\documentclass{jktrX}
\newcommand{\embfig}[3]{
\medskip\centerline{
         \includegraphics[width=.#1\textwidth]{#2}
}\centerline{#3}}

\newcommand{\Abb}{\embfig{90}{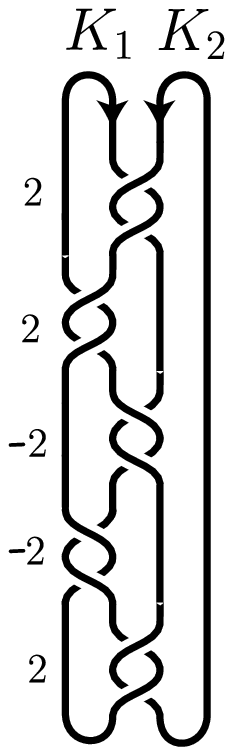}{Figure 2.1:
$S(34,21)$
}}

\newcommand{\Acc}{\embfig{55}{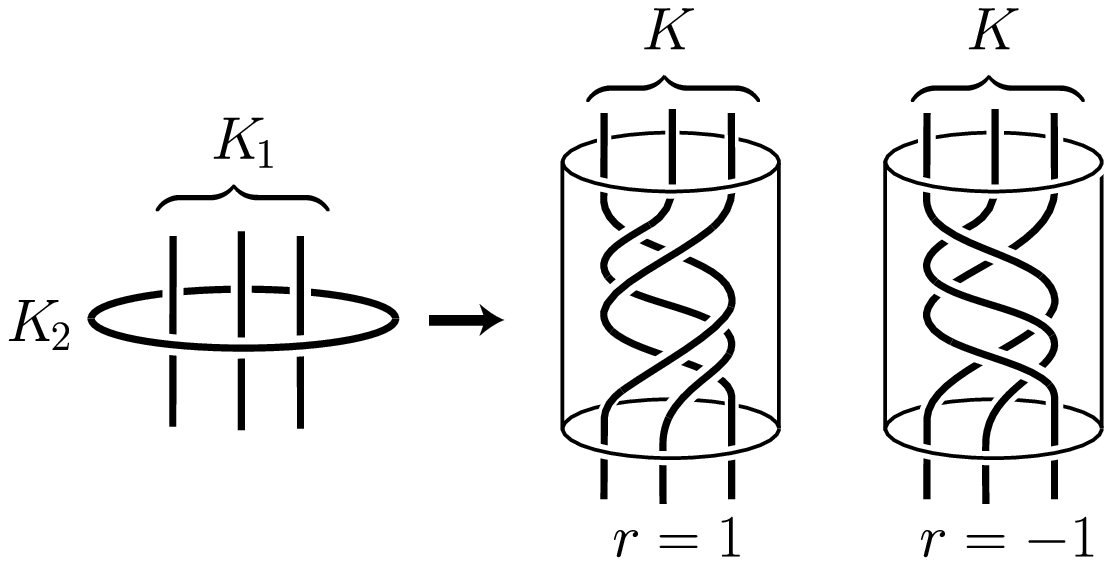}{Figure 2.2:
Dehn twists along $K_2$
}}

\newcommand{\Baa}{\embfig{40}{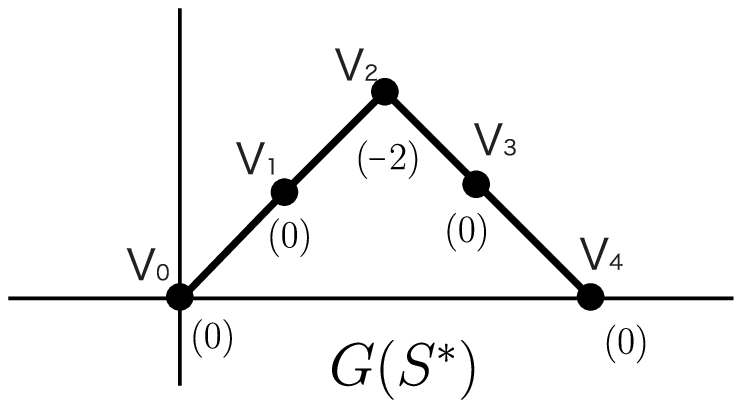}{Figure 3.1:
Graph $G(S^{*})$ for $S^{*}=[2,0,2,-2,-2,0,-2]$
}}
\newcommand{\Bb}{\embfig{75}{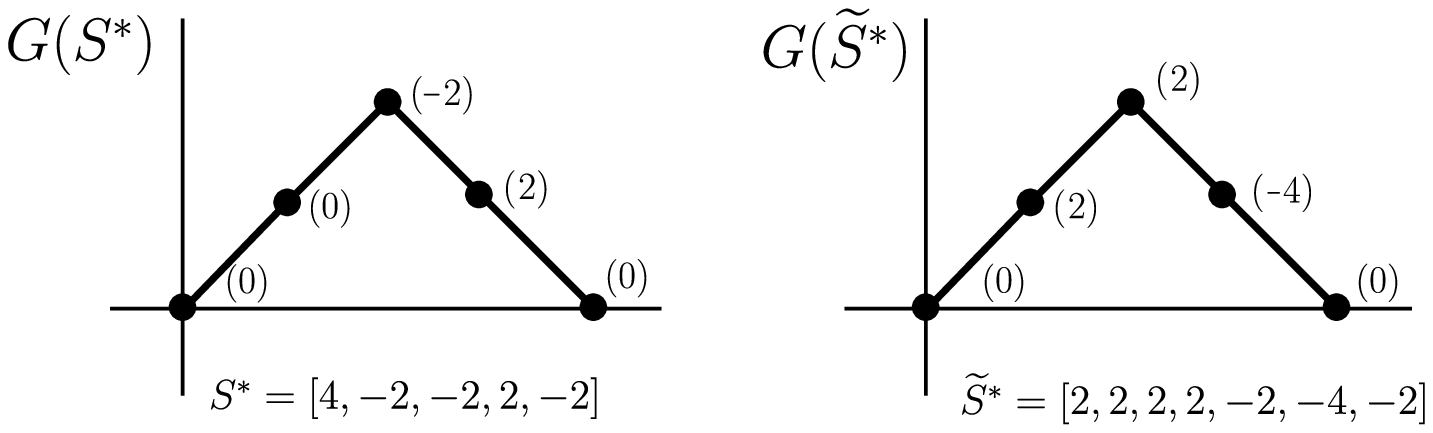}{Figure 3.2:
The graphs for $S^{*}$ and 
$\widetilde{S}^{*}$
}}

\newcommand{\foA}{\embfig{80}{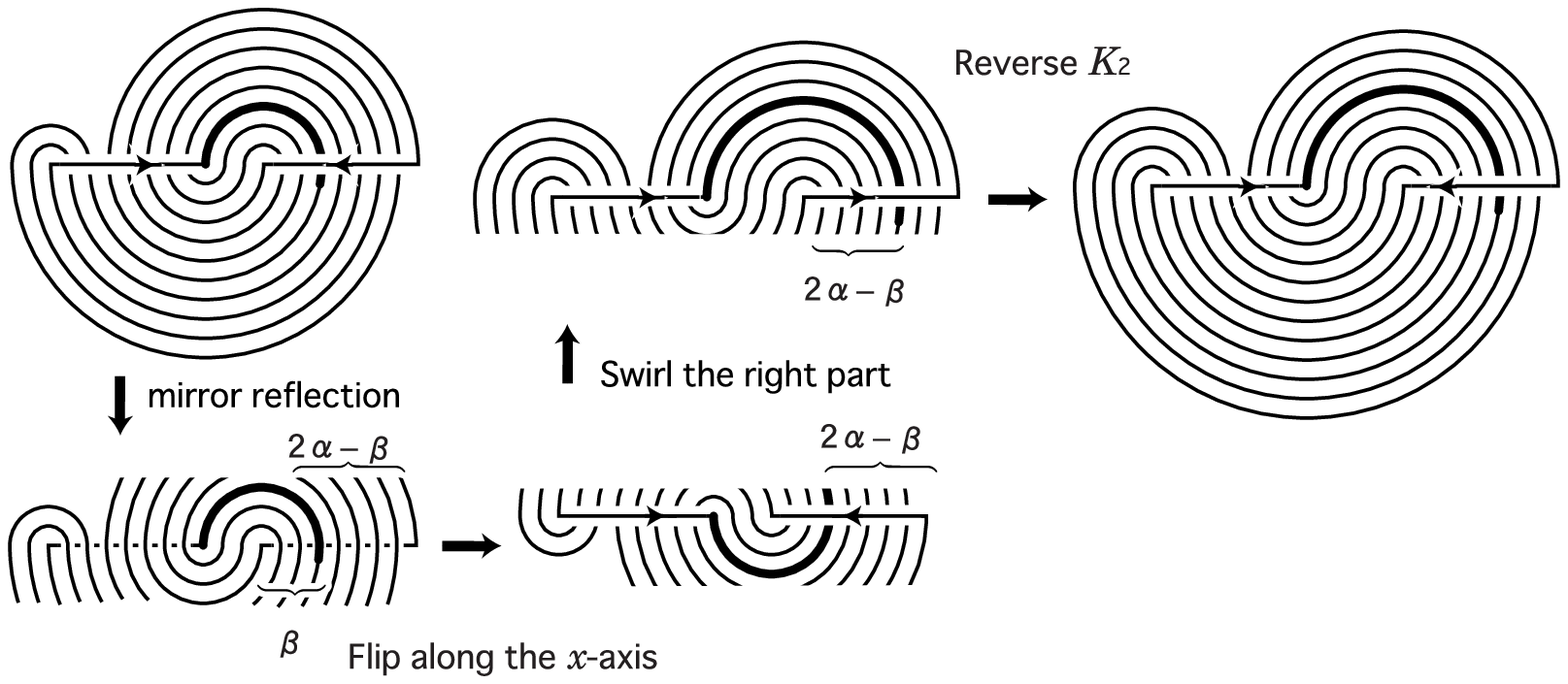}{Figure 3.3:
Deformation from $B(2\alpha,\beta)$ to 
$B(2\alpha, 2\alpha-\beta)$, e.g. $B(8,3)$
to $B(8,5)$
}}

\newcommand{\foC}{
\embfig{95}{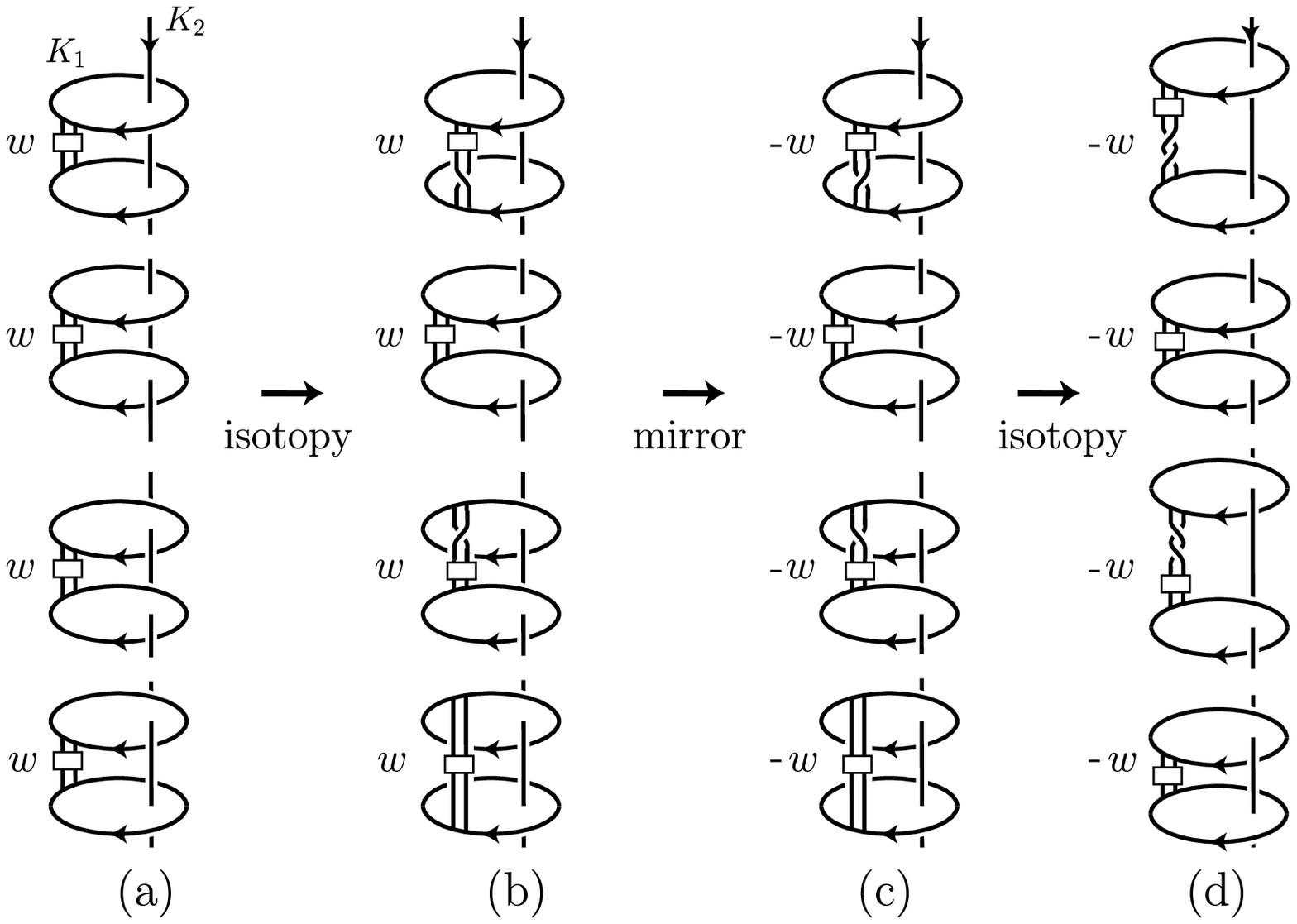}{Figure 3.5:
Deformation from $B(2\alpha,\beta)$ to 
$B(2\alpha, 2\alpha-\beta)$.
}}

\newcommand{\foD}{\embfig{50}{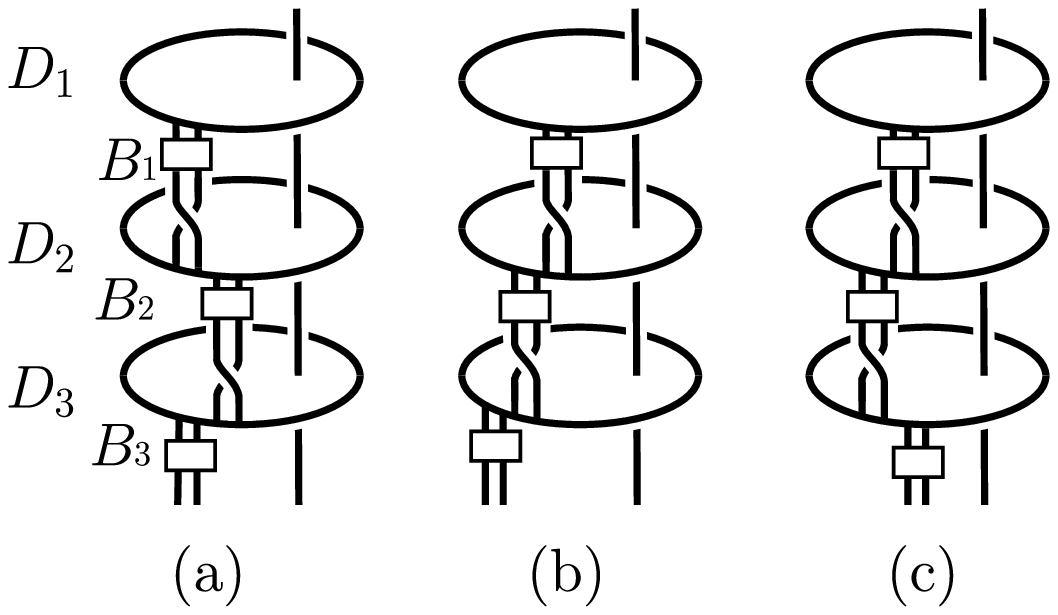}{Figure 3.6:
Sliding bands to change relative positions
}}

\newcommand{\foE}{\embfig{95}{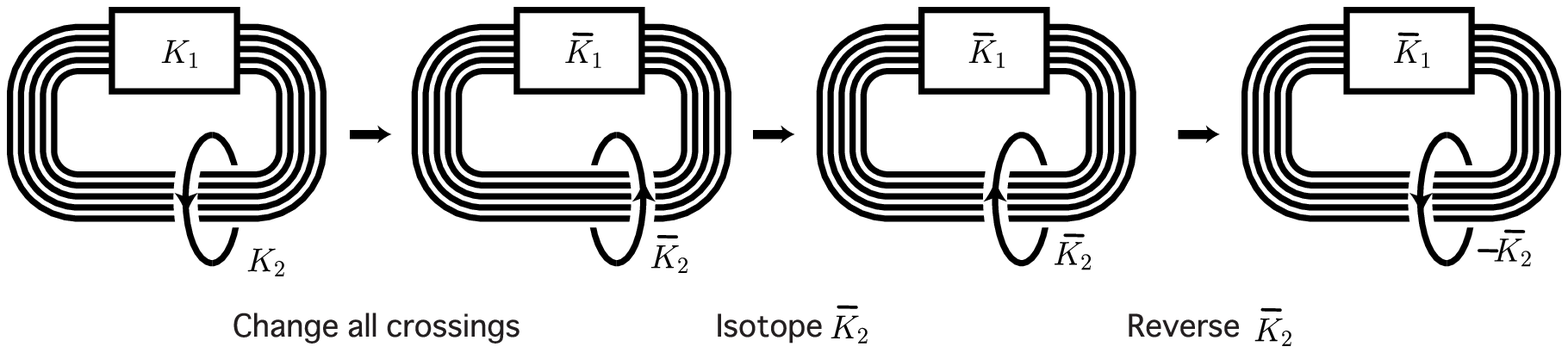}{Figure 3.7:
Reflect $K_1\cup K_2$ and reverse $K_2$
}}

\newcommand{\sixa}{\embfig{70}{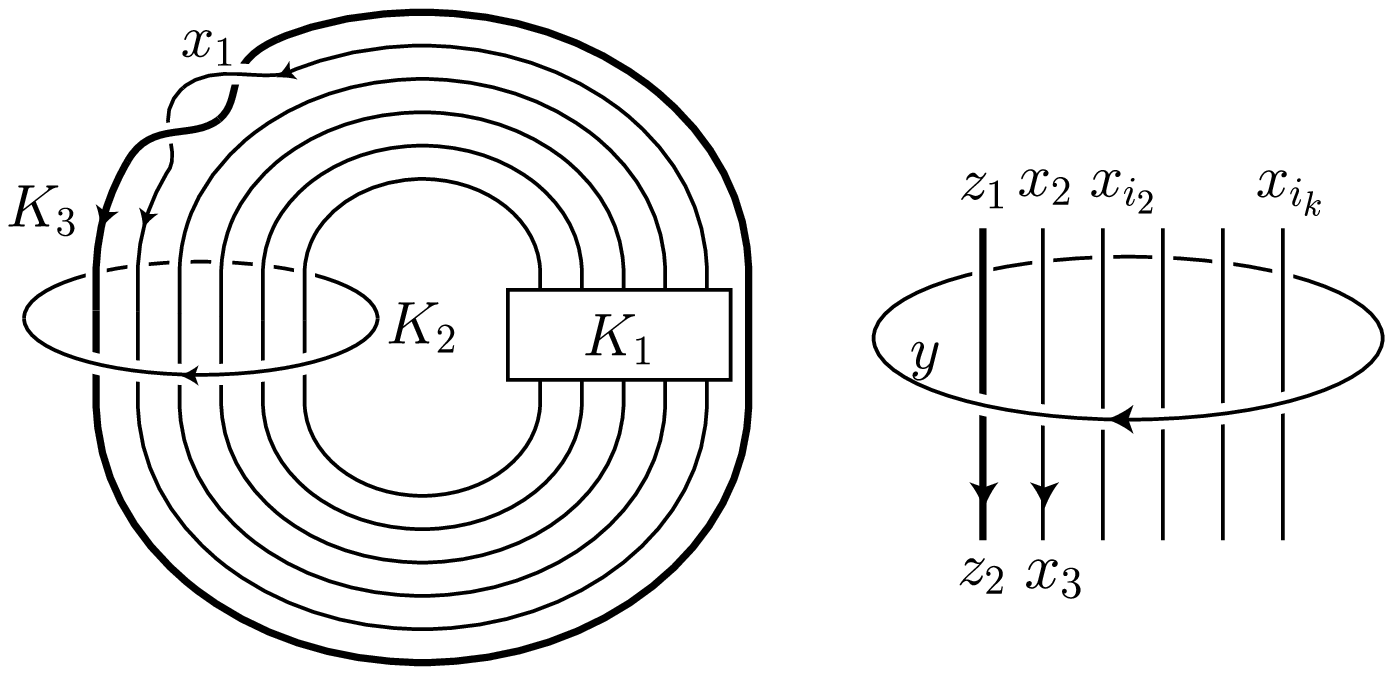}{Figure 4.1: 
A diagram of 
$\widetilde{L}=K_1\cup K_2\cup K_3$
}}

\newcommand{\tena}{\embfig{82}{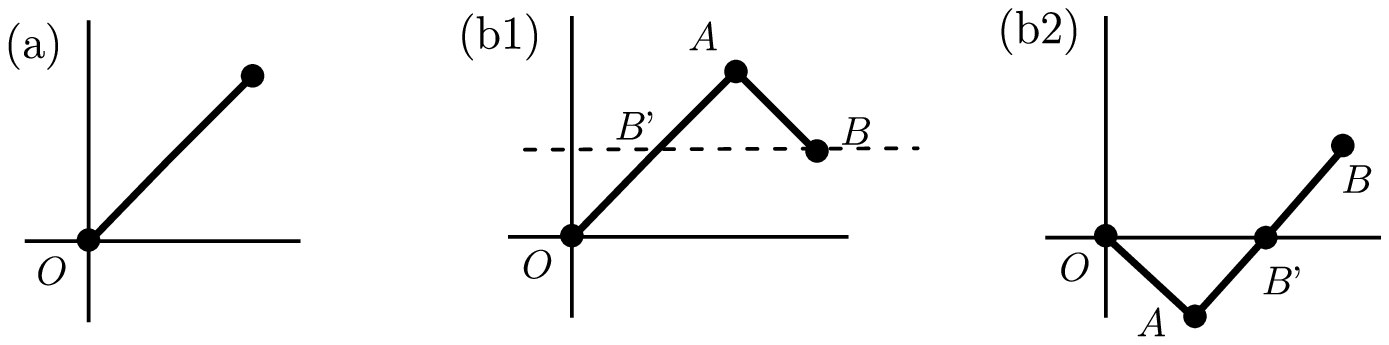}{Figure 7.1:
Graphs $G$ for $\nu(G)=2$ and $3$
}}

\newcommand{\tenc}{\embfig{96}{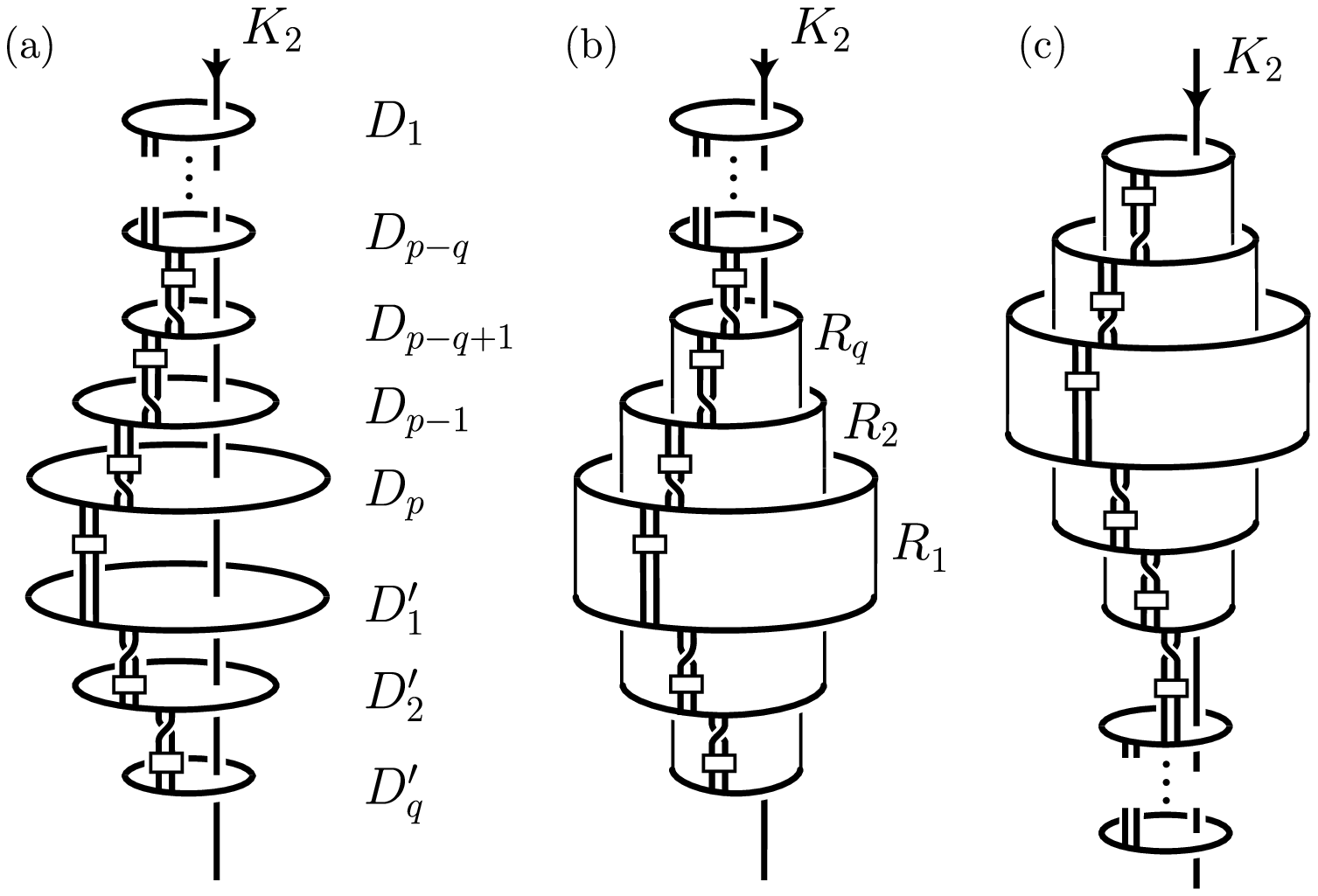}{Figure 7.2:
Construction of a canonical surface
}}

\newcommand{\tend}{\embfig{98}{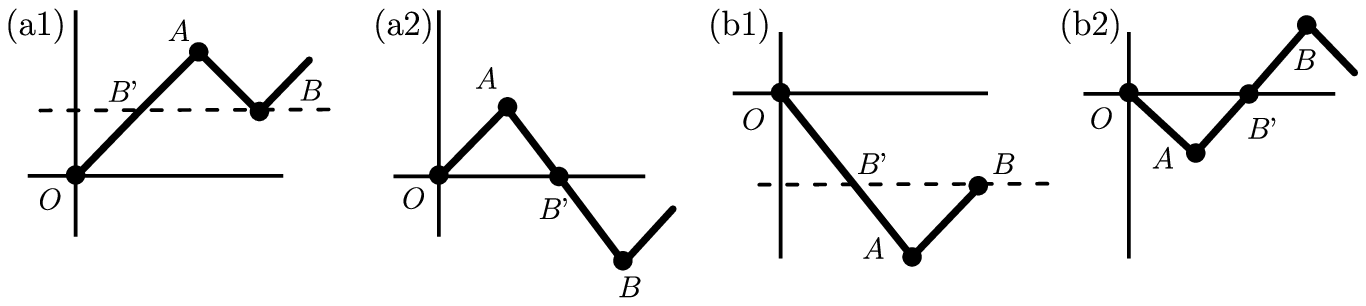}{Figure 7.3:
Graph $G$ for $\nu(G)\ge 4$,
where $O$ is a local minimal, or maximal
}}

\newcommand{\tenf}{\embfig{66}{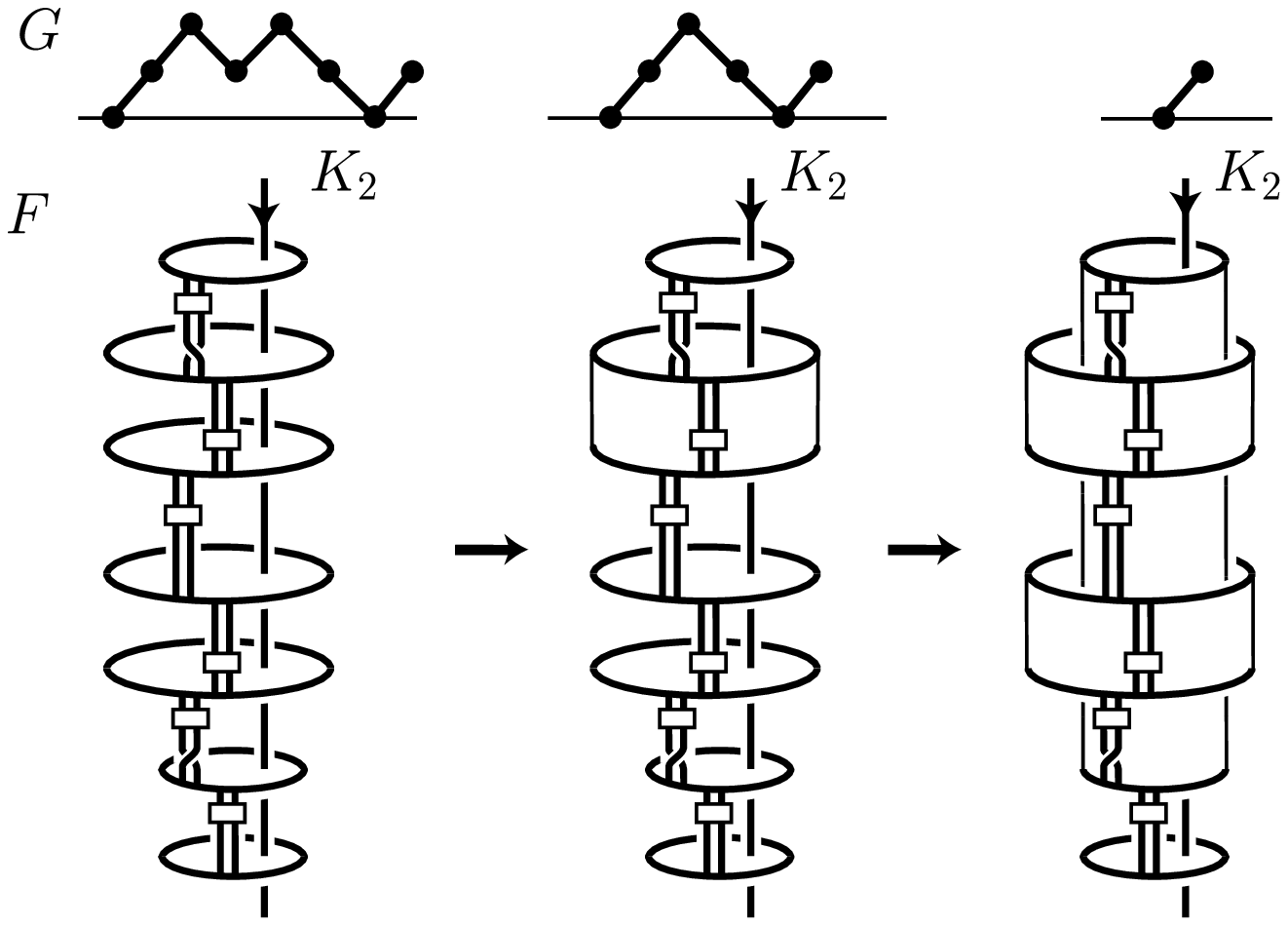}{Figure 7.4:
Construction of a canonical surface from a graph $G$
(I)}}

\newcommand{\teng}{\embfig{97}{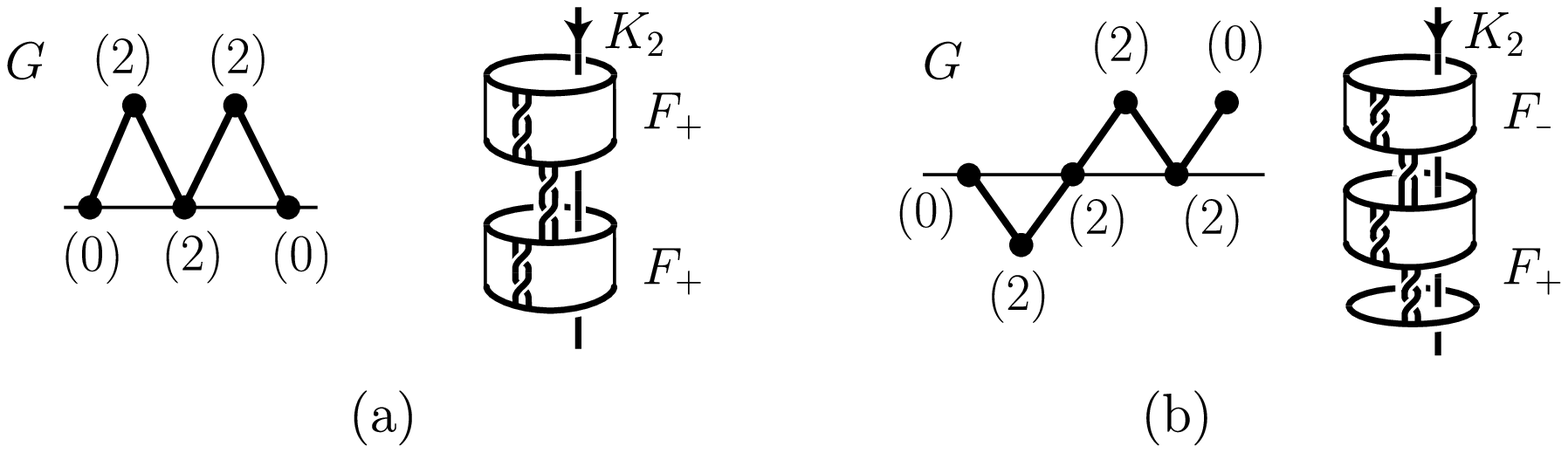}{Figure 7.5:
Construction of a canonical surface from $G$ (II)
}}
\newcommand{\tenh}{\embfig{96}{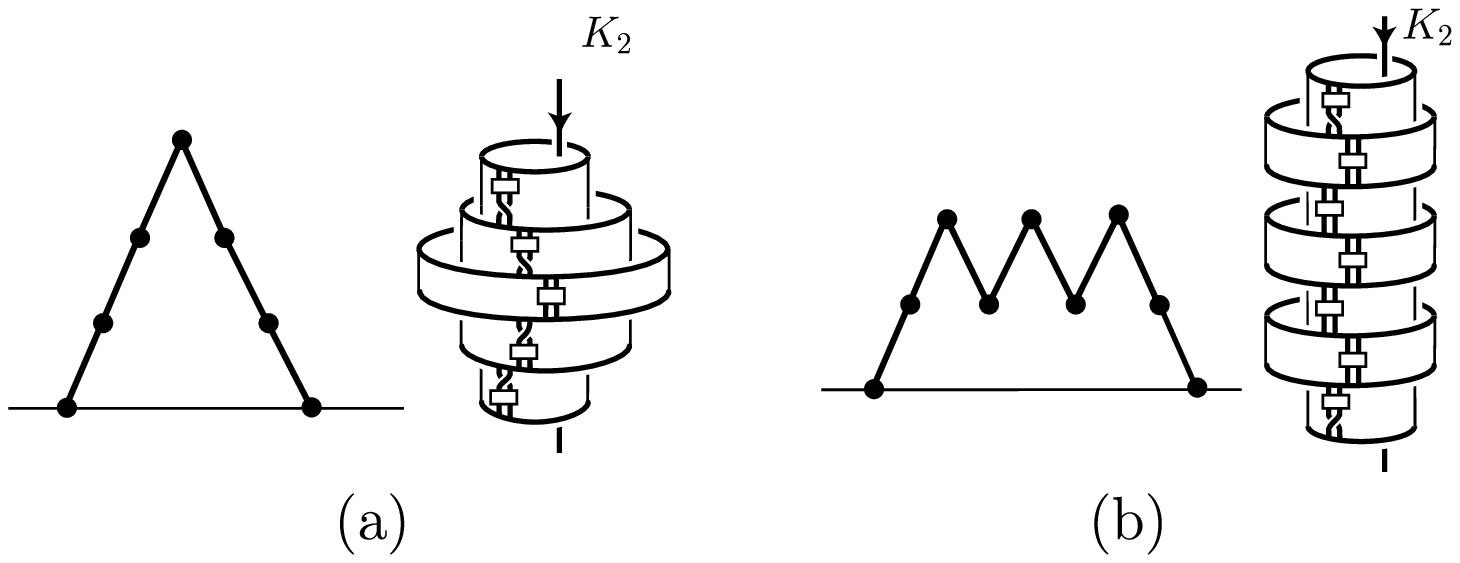}{Figure 7.6:
Construction of a canonical surface from $G$ (III)
}}

\newcommand{\fotnA}{\embfig{75}{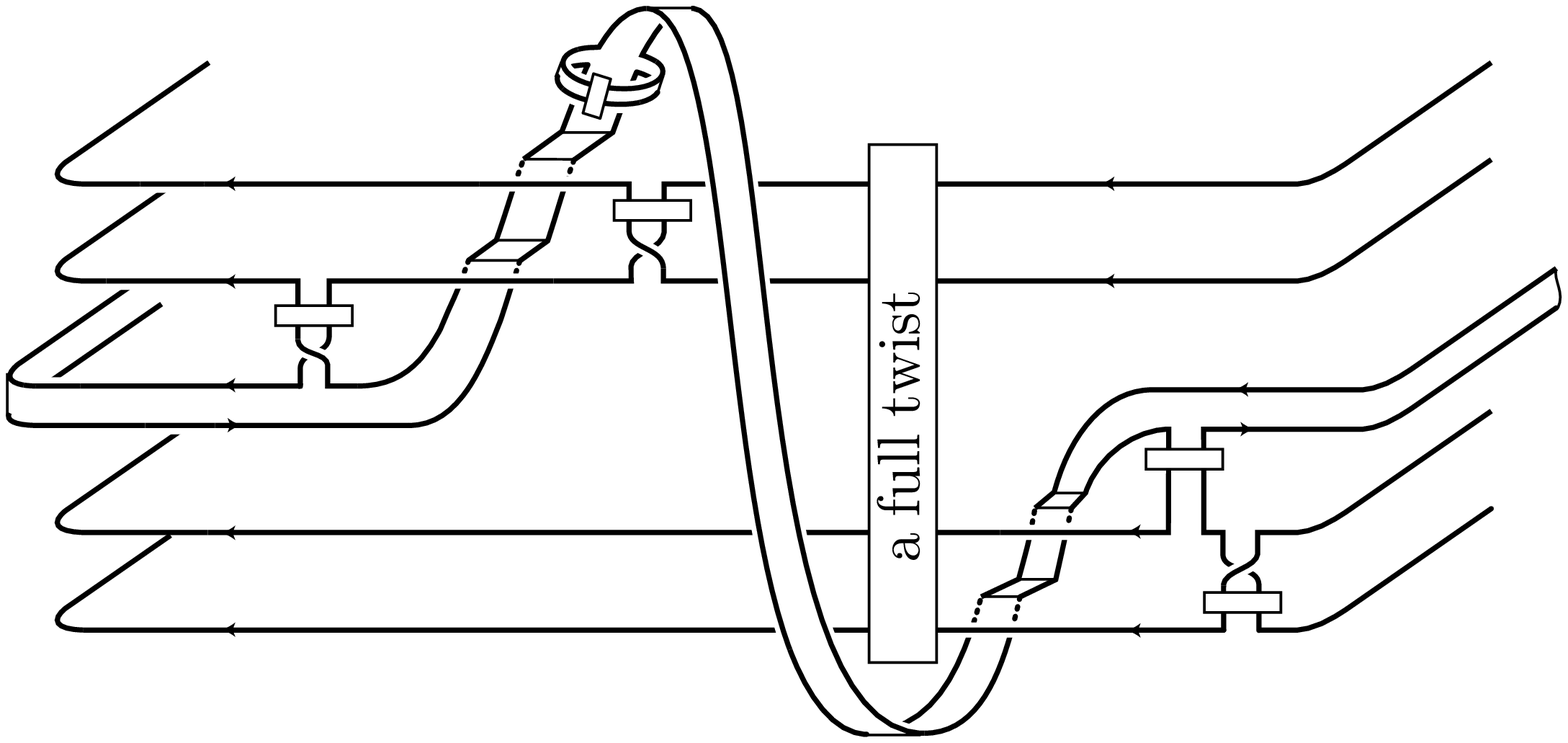}{Figure 7.7:
A canonical surface $F(r)$, where $r=1$}}

\newcommand{\constAa}{\embfig{72}{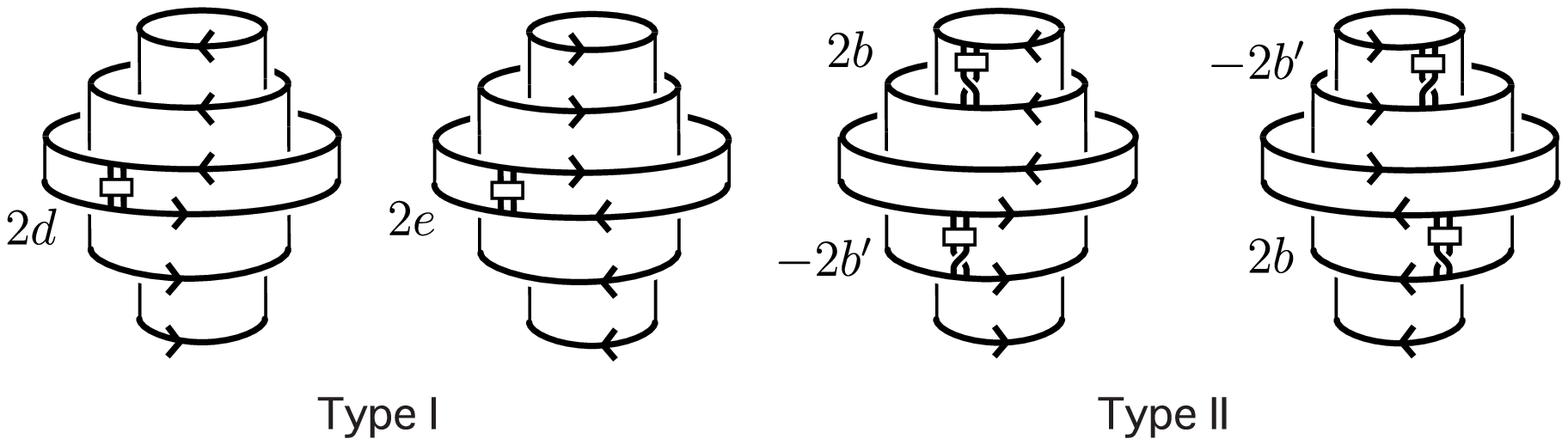}{Figure 7.8.1:
Bands of Types I and II}}

\newcommand{\constAb}{\embfig{60}{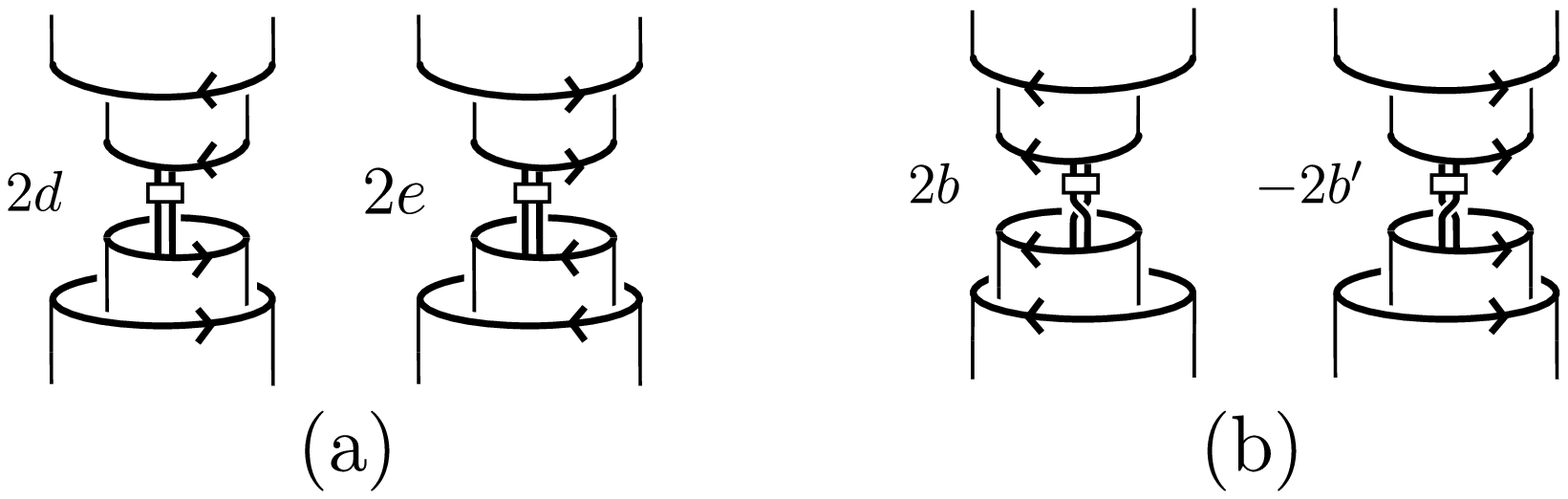}
{Figure 7.8.2:
Bands of Type III}}

\newcommand{\constAc}{\embfig{60}{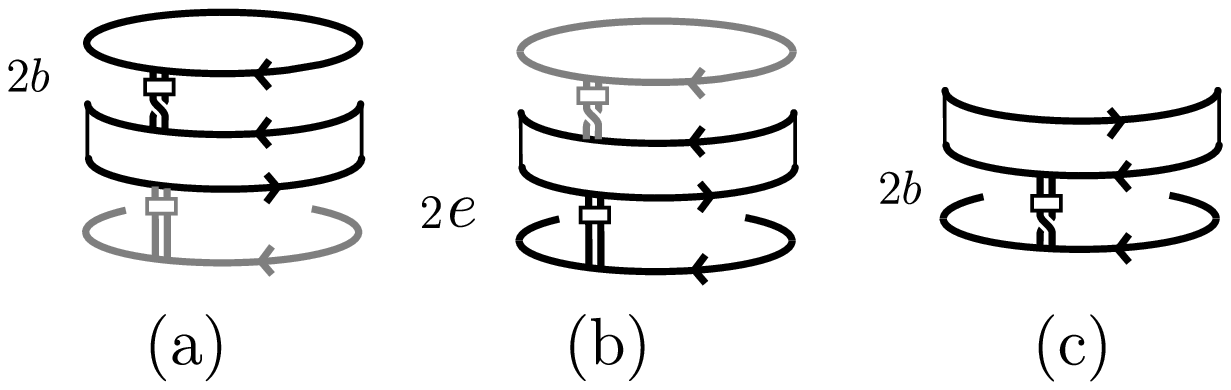}{Figure 7.8.3:
Bands of Type IV}}

\newcommand{\constAd}{\embfig{75}{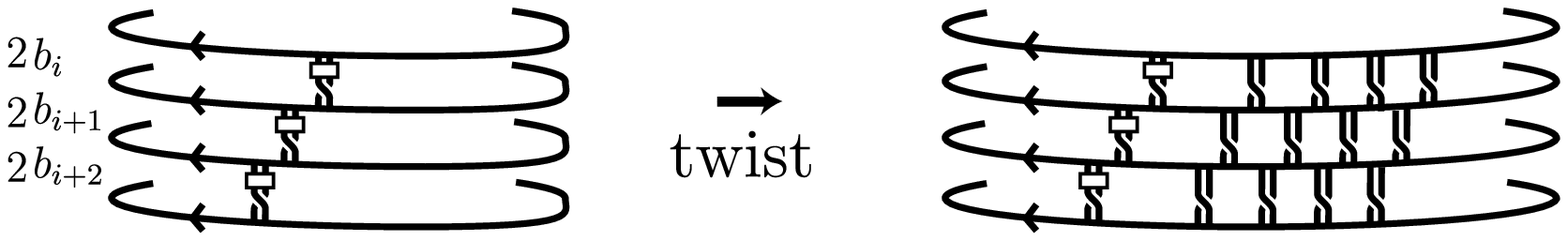}{Figure 7.8.4:
Bands of Type V}}

\newcommand{\constB}{\embfig{60}{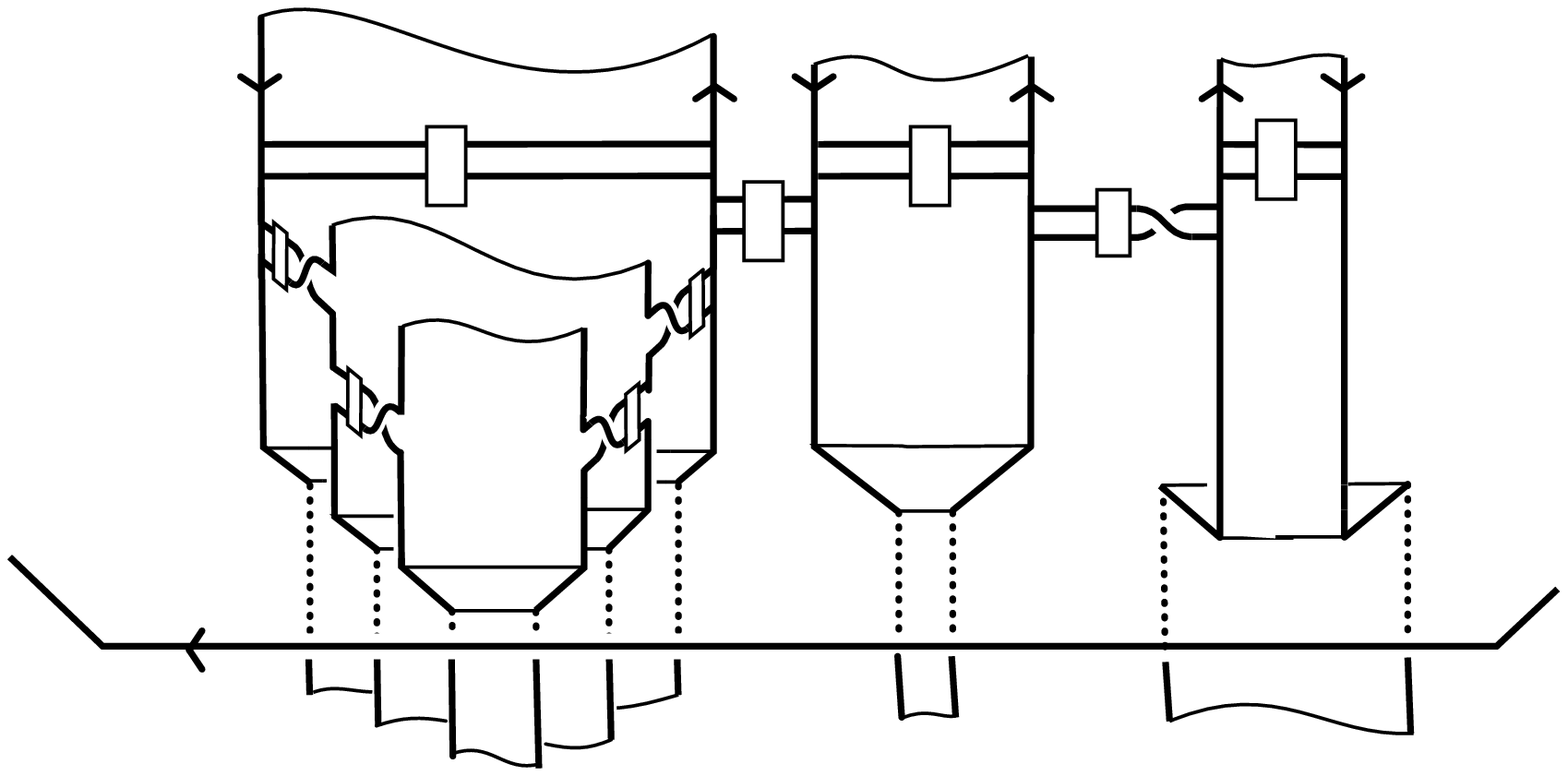}{Figure 7.9:
Placement of the bands of Types 
I, II and III}}

\newcommand{\fotnC}{\embfig{70}{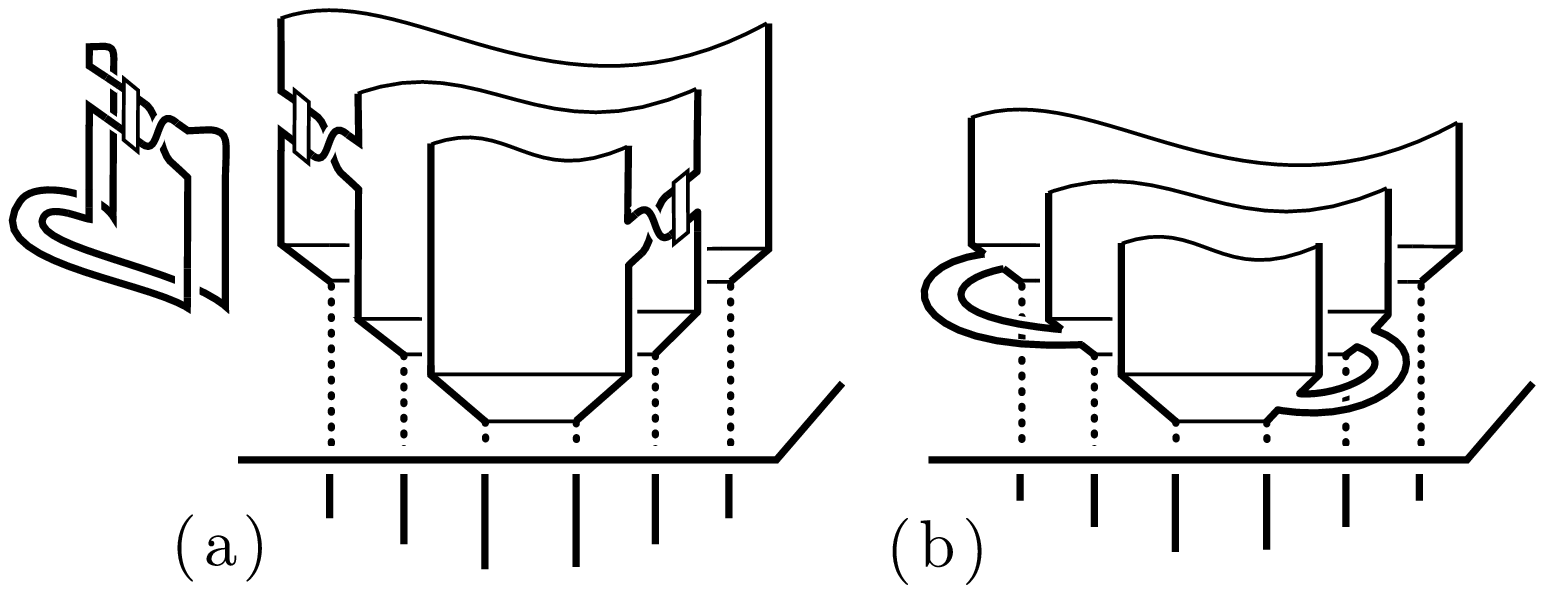}{Figure 8.1:
Deplumbing and compression for Type II}}

\newcommand{\fotnD}{\embfig{50}{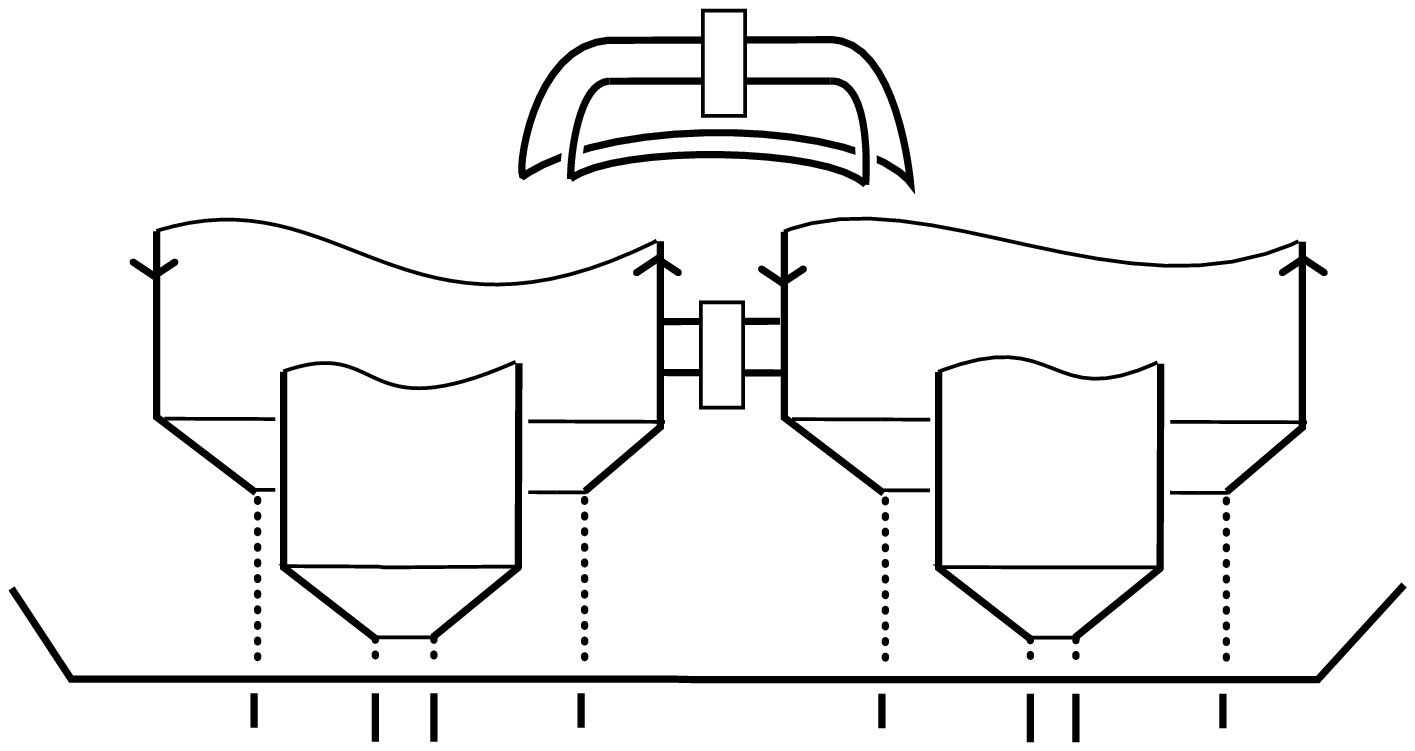}{Figure 8.2:
Deplumbing an annulus for Type III (a)}}

\newcommand{\fotnE}{\embfig{85}{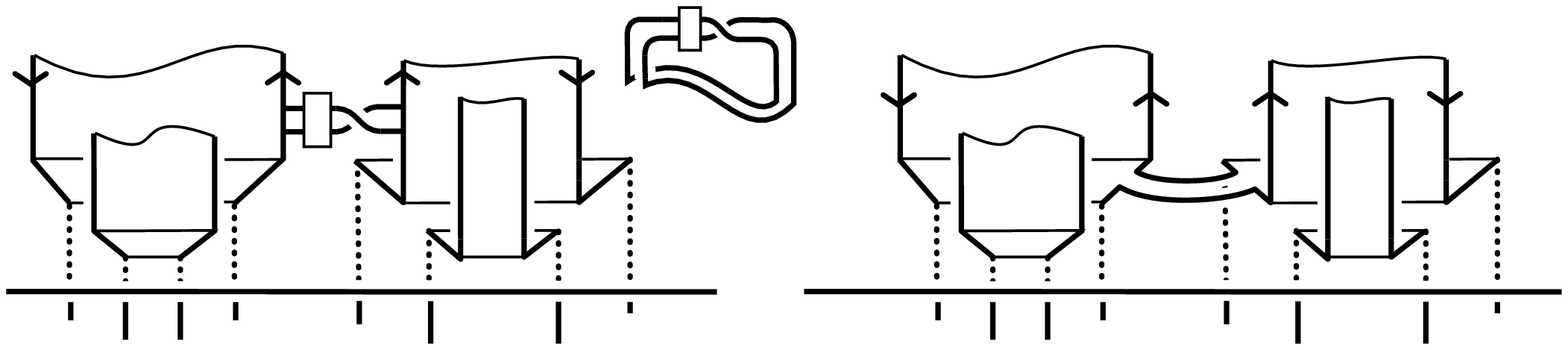}{Figure 8.3:
Deplumbing and compression for Type III (b)}}

\newcommand{\fotnF}{\embfig{60}{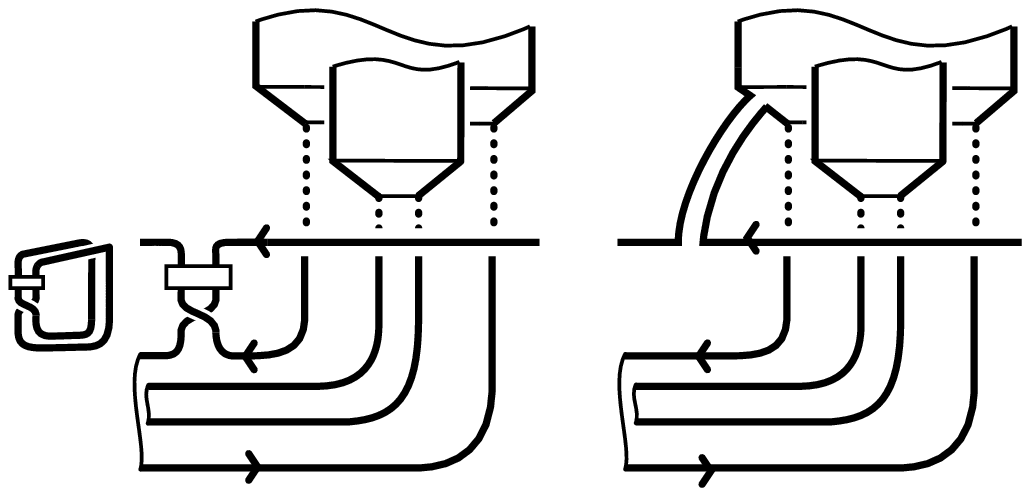}{Figure 8.4:
Deplumbing and compression for Type IV (a)}}

\newcommand{\fotnG}{\embfig{70}{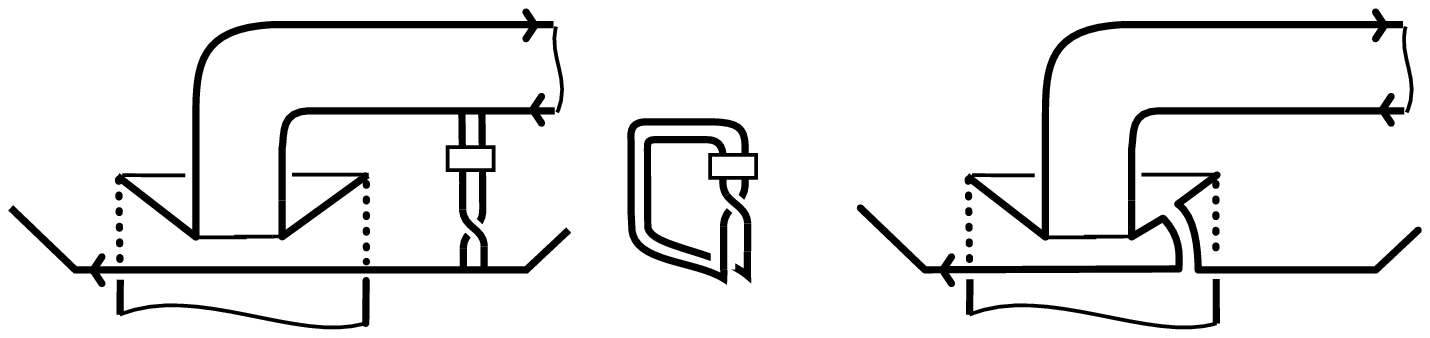}{Figure 8.5:
Deplumbing and compression for Type IV (c)}}

\newcommand{\fotnH}{\embfig{80}{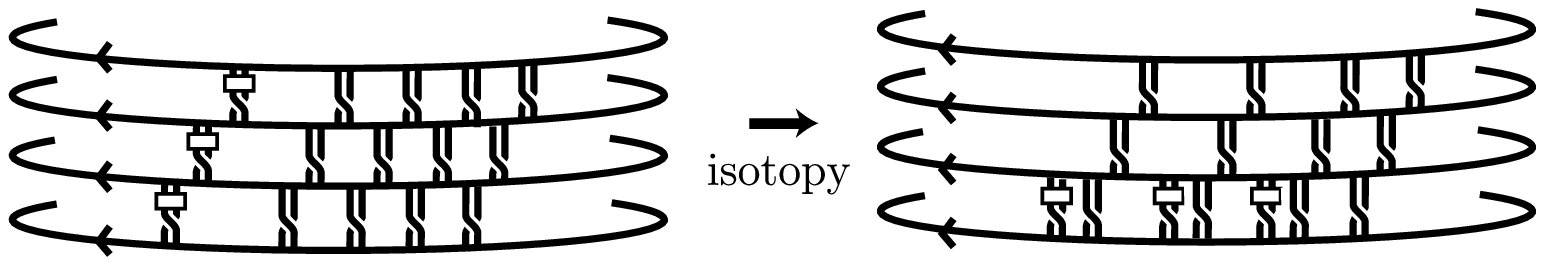}{Figure 8.6:
Deplumbing and compression for Type V}}
\newcommand{\fif}{\embfig{40}{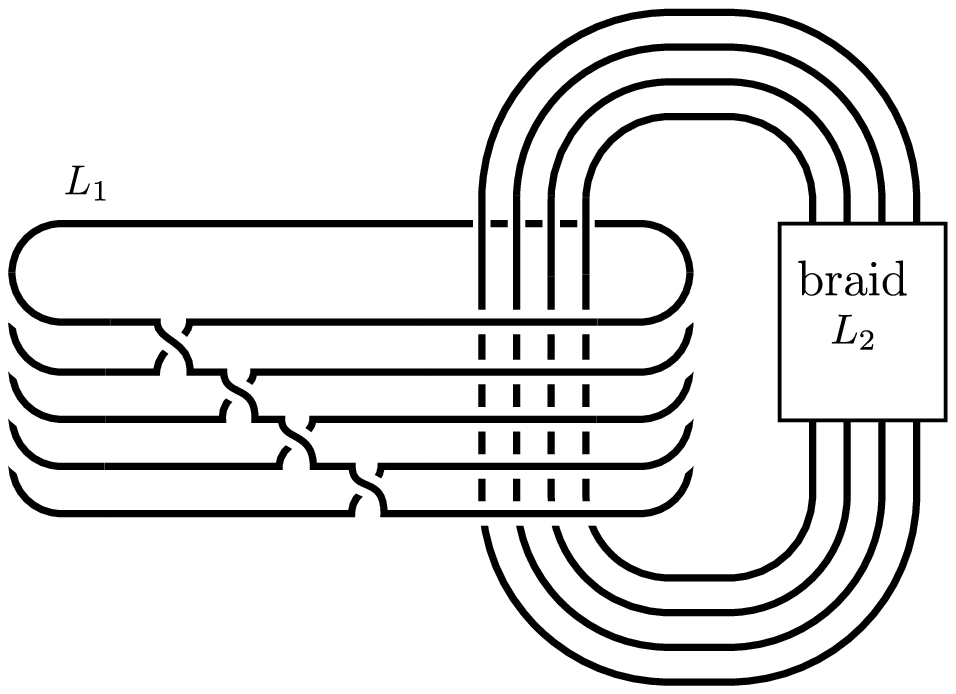}{Figure 9.1:
A braid penetrating a disk
}}

\newcommand{\stnA}{\embfig{58}{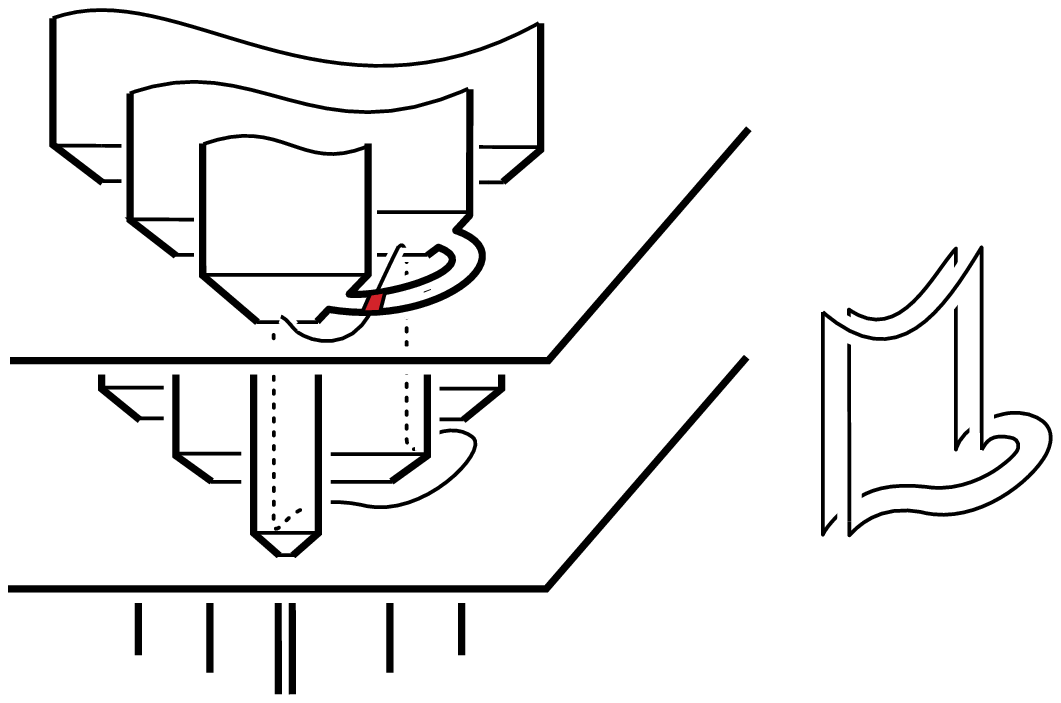}{Figure 9.2:
Plumbing a Hopf band for Type II
}}
\newcommand{\stnB}{\embfig{58}{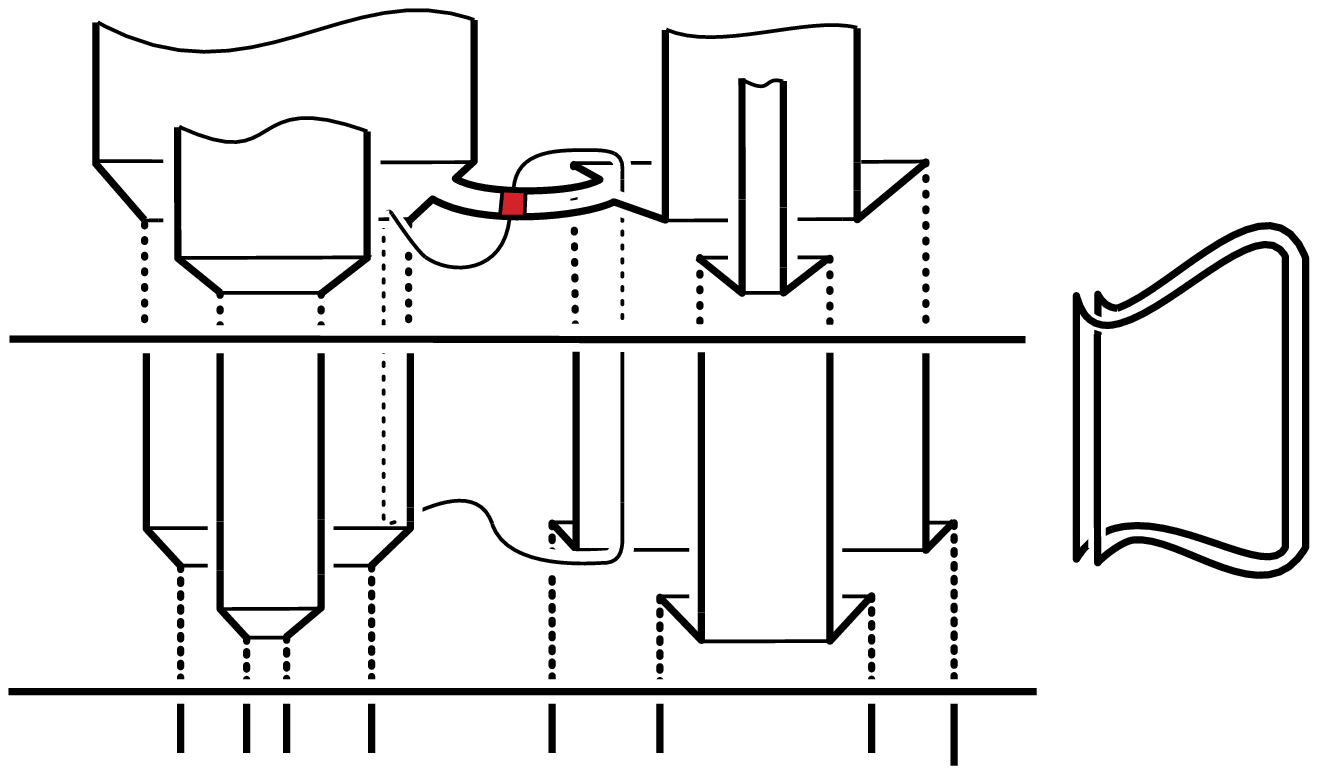}{Figure 9.3:
Plumbing a Hopf band for Type III
}}

\newcommand{\stnC}{\embfig{98}{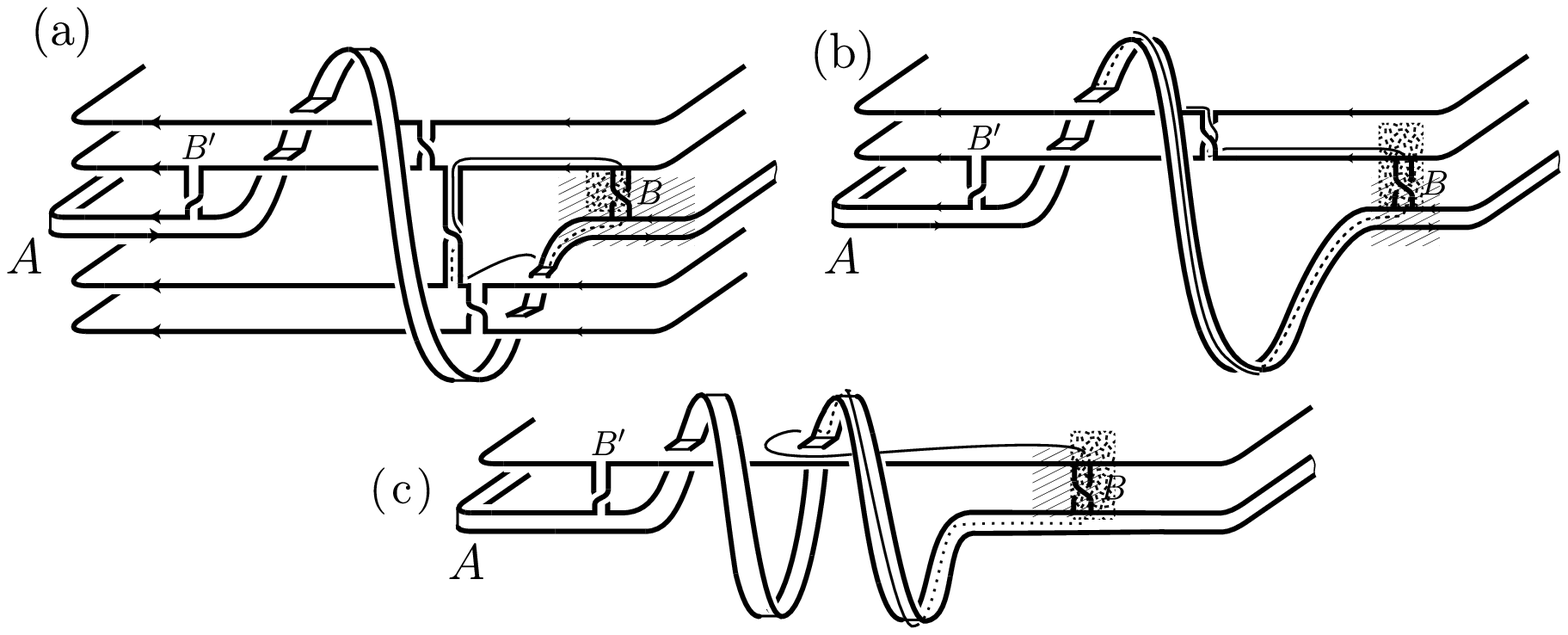}{Figure 9.4:
Plumbing a Hopf band for Type IV}}
\newcommand{\elebb}{\embfig{65}{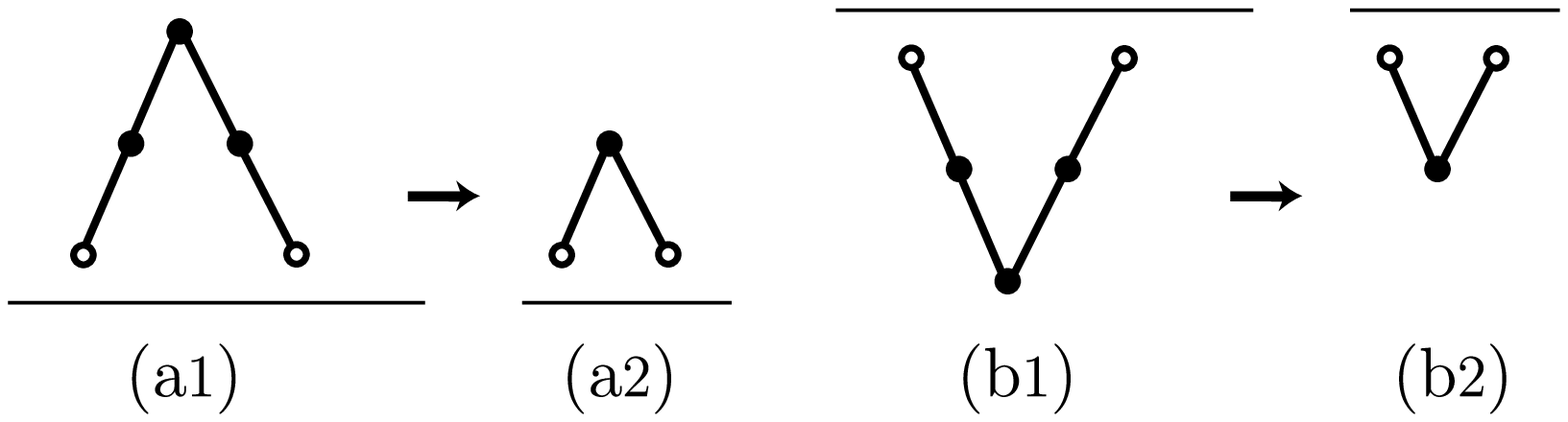}
{Figure 10.1:
Reduction on the graph (I)
}}

\newcommand{\eleb}{\embfig{96}{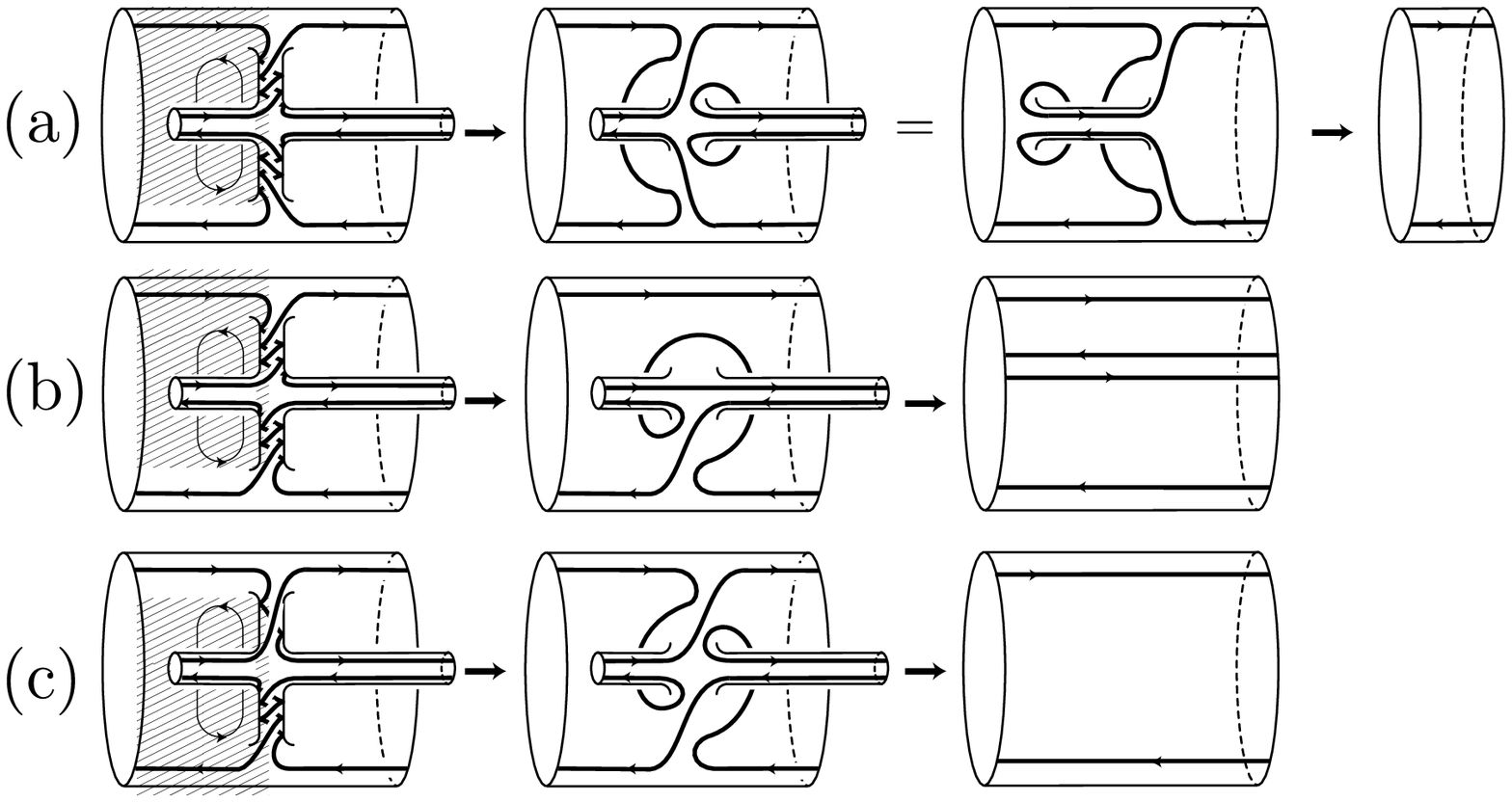}
{Figure 10.2:
Reduction by complementary disk decompositions
}}

\newcommand{\elecc}{\embfig{95}{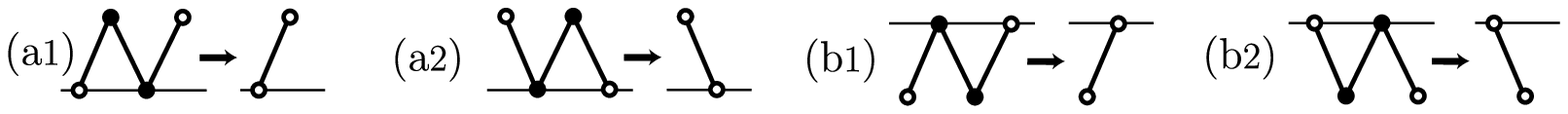}
{Figure 10.3:
Reduction on the graph (II)
}}
\newcommand{\eledd}{\embfig{26}{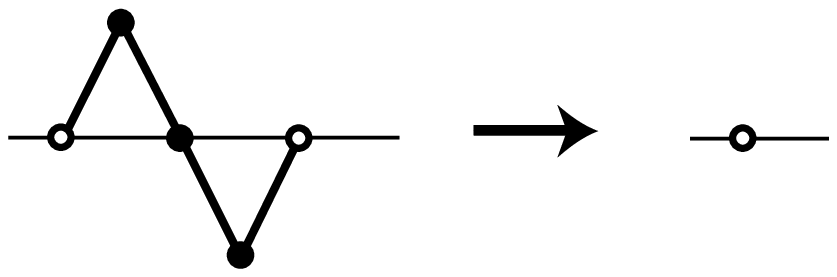}
{Figure 10.4:
Reduction on the graph (III)
}}

\newcommand{\eleee}{\embfig{77}{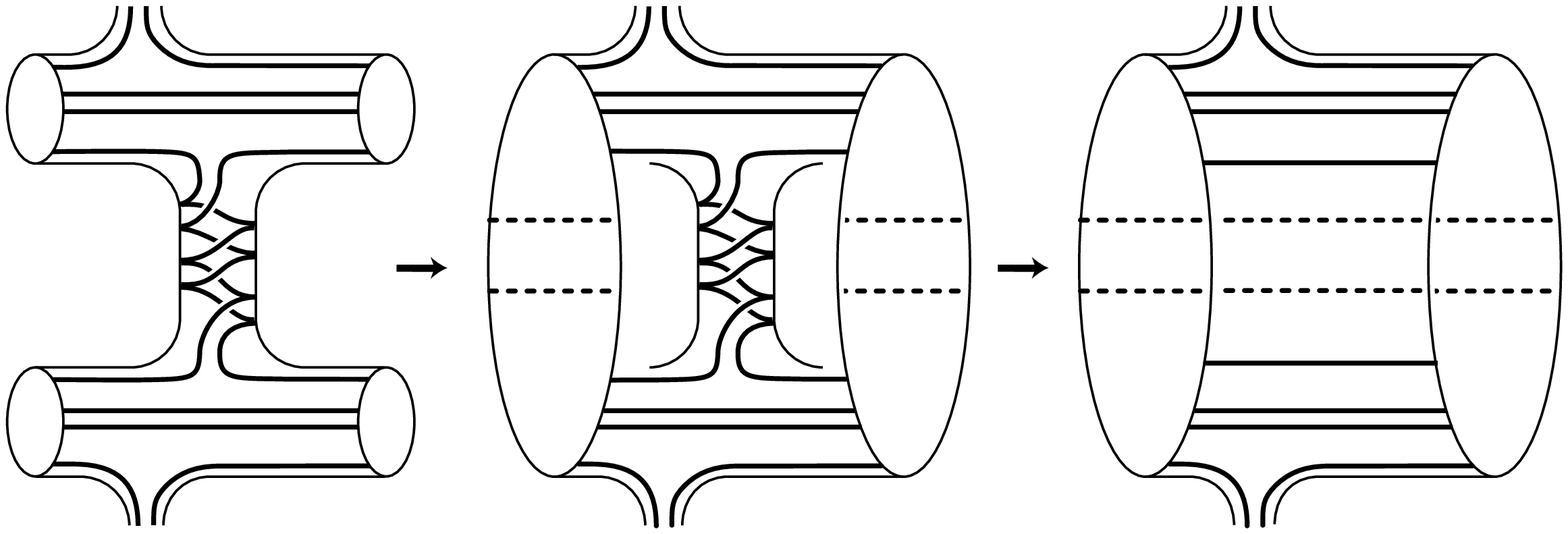}
{Figure 10.5:
Amalgamation of tori (I)
}}
\newcommand{\eleff}{\embfig{70}{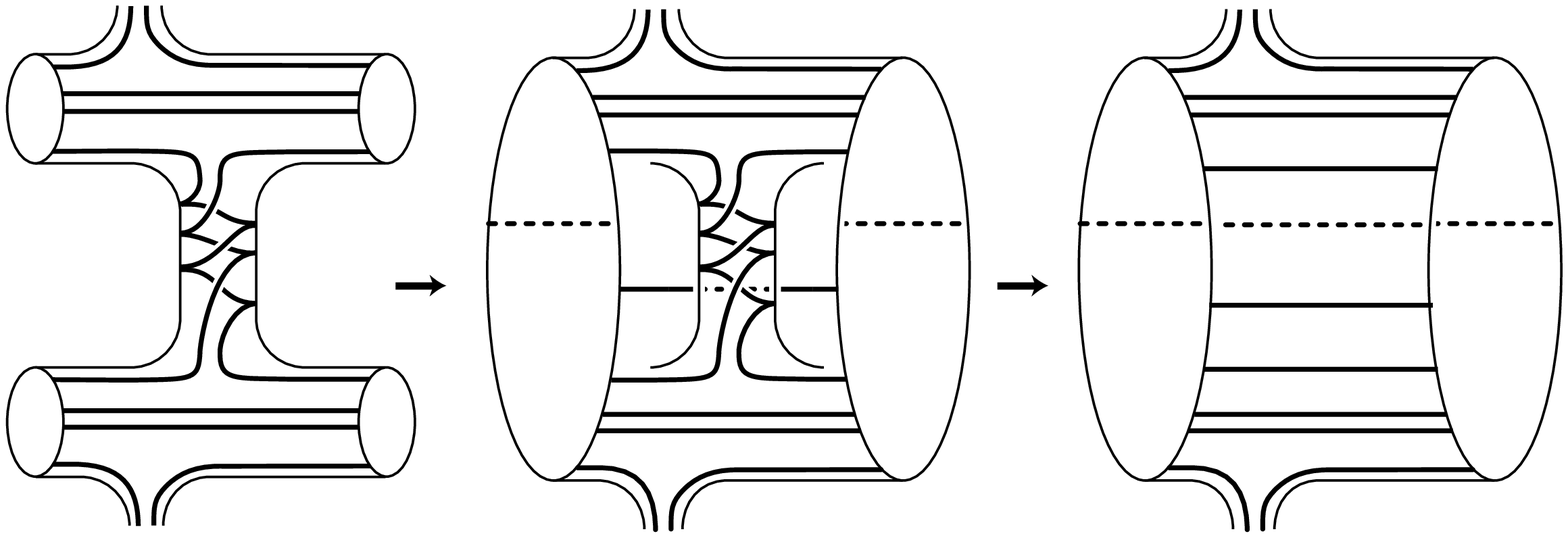}
{Figure 10.6:
Amalgamation of tori (II)
}}

\newcommand{\thir}{\embfig{85}{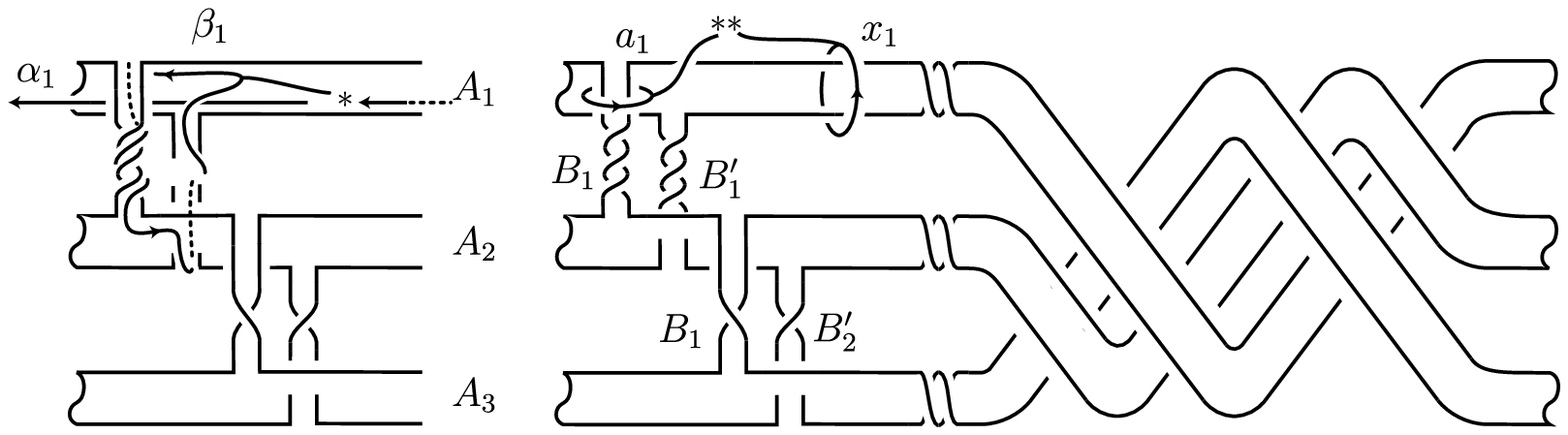}{Figure 11.1:
Generators of $\pi_1(S^3-\hat{F})$ and $\pi_1(\hat{F})$
}}

\newcommand{\etna}{\embfig{85}{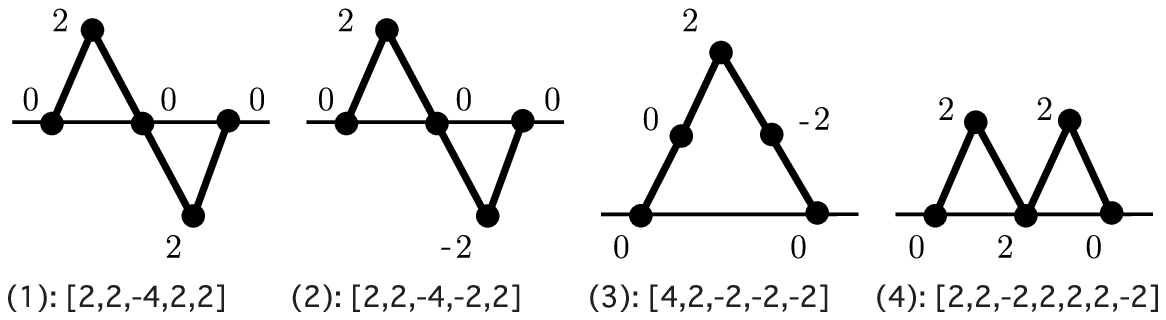}{Figure 12.1:
The graphs for $2$-bridge links with $\ellk=0$
}}

\newcommand{\etnb}{\embfig{90}{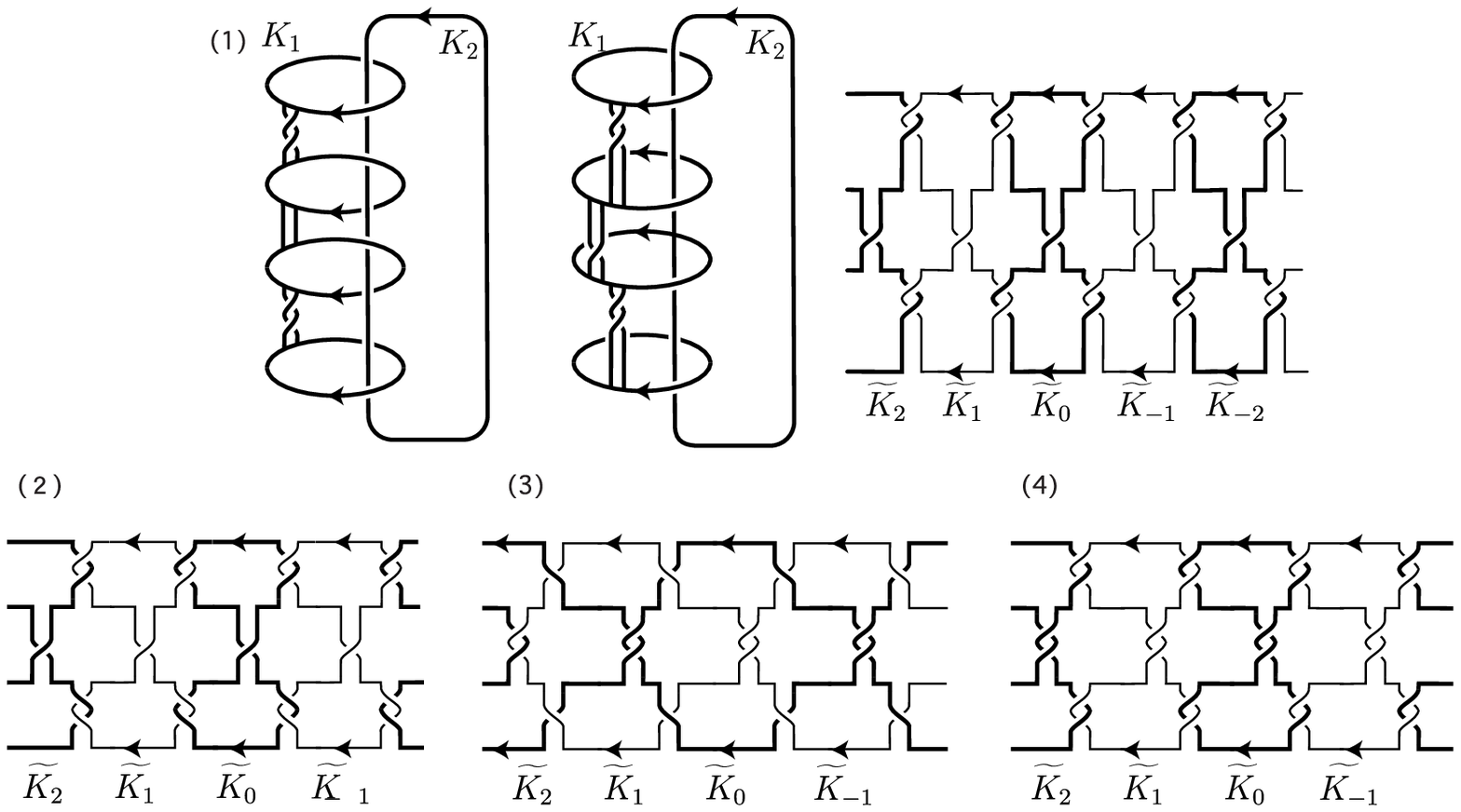}{Figure 12.2:
Lifts of $K_1$ in $M^3$
}}

\newcommand{\etnc}{\embfig{80}{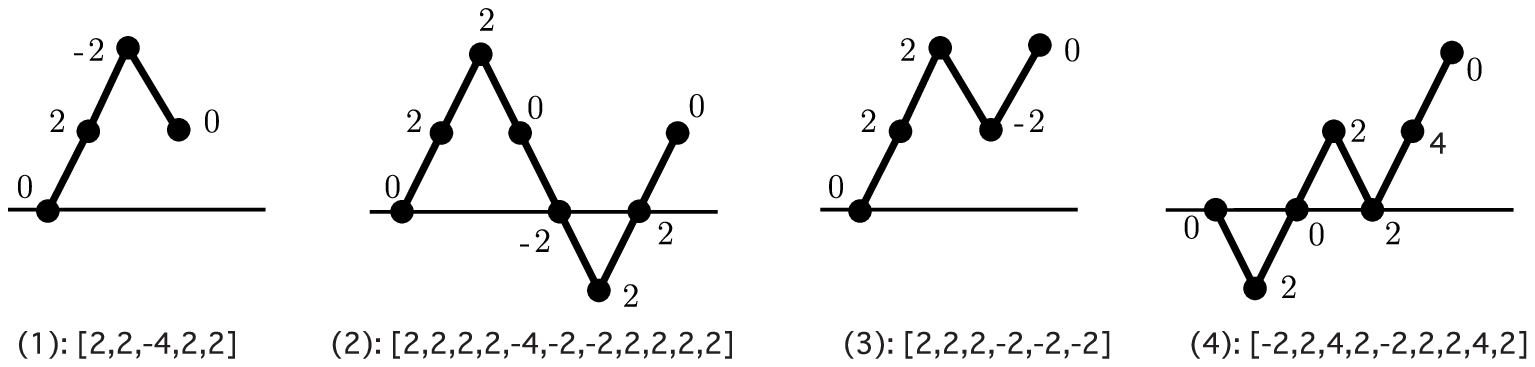}{Figure 12.3:
More graphs for $2$-bridge links with $\ellk\neq 0$
}}

\newcommand{\satel}{\embfig{80}{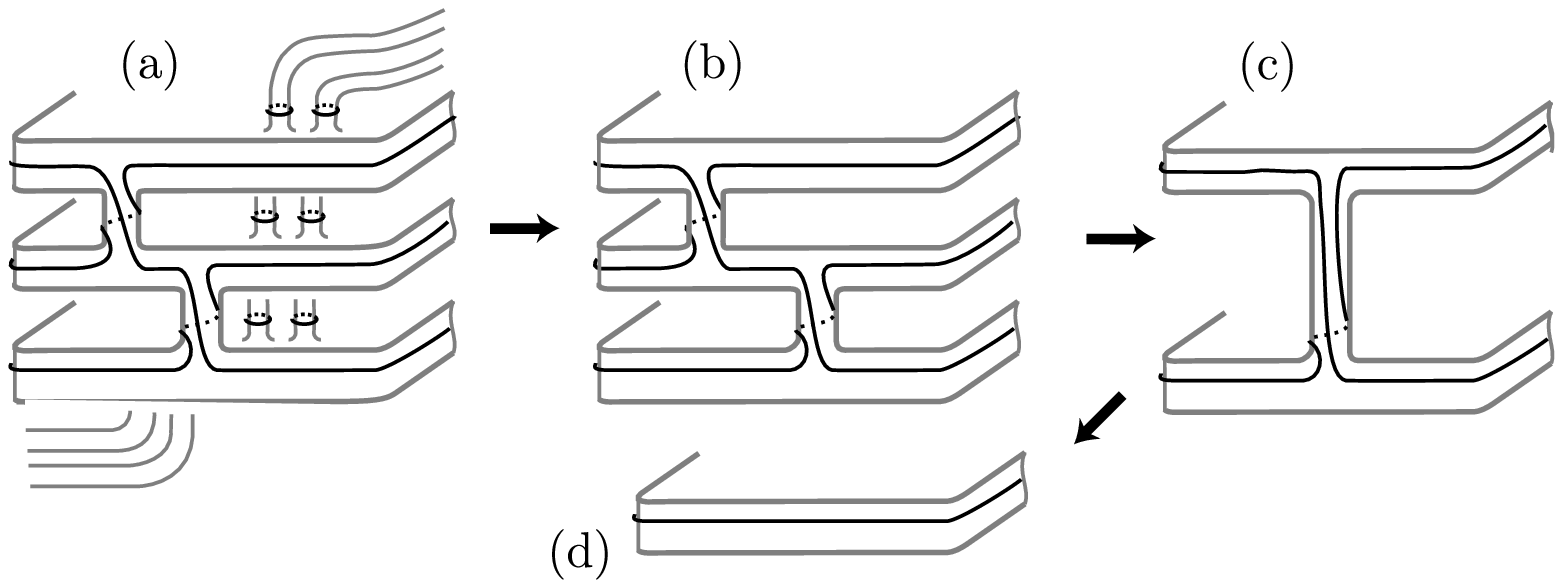}{Figure 13.1:
$C$-product decompositions
}}

\newcommand{\eighteen}{\embfig{80}{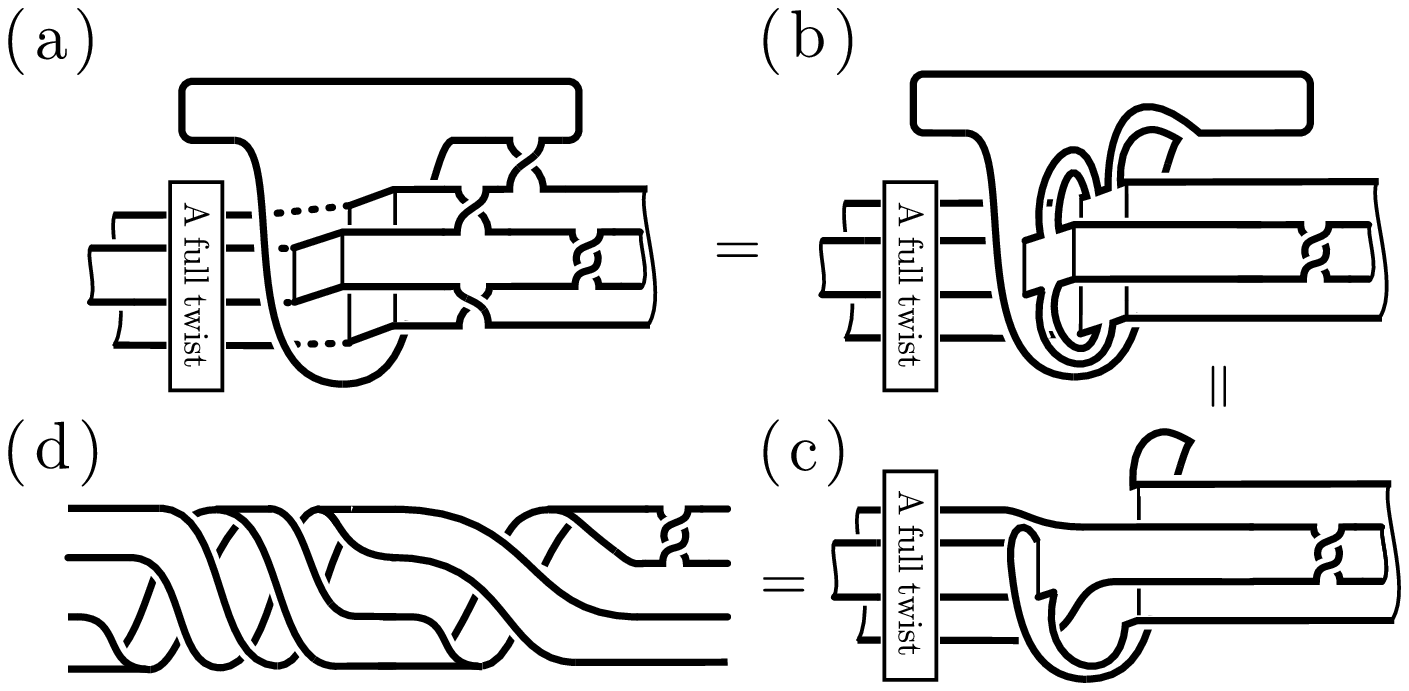}
{Figure 14.1:
$g(K(46,39|1))=1$
}}

\def\draftnote{ }
\def\today{}

\usepackage{amsfonts}
\usepackage{amssymb}
\usepackage{amsmath}
\usepackage{graphicx}

\def\qed{\hfill \fbox{}}
\renewcommand{\hat}[1]{\widehat{#1}}
\def\bysame{\leavevmode\hbox to3em{\hrulefill}\thinspace}

\newcommand{\secti}[1]{
\addtocounter{section}{1}
\setcounter{thm}{0}
\setcounter{equation}{0}
\noindent
{\bf \thesection.\ #1 }
\medskip
}%

\def\ellkB{\mbox{$\ellk B(2\alpha,\beta)$}}

\def \sone{\raisebox{-8pt}
{\begin{picture}(26,16)(-2,0) \thicklines  
\put(0,25){\line(1,-2){5}} 
\put(8,10){\line(1,-2){5}}  
\put(0,0){\line(1,2){12.5}} 
\put(22,0){\line(0,1){25}} 
\end{picture}}
}

\def \stwo{\raisebox{-8pt}
{\begin{picture}(26,16)(-2,0) \thicklines  
\put(10,25){\line(1,-2){5}} 
\put(18,10){\line(1,-2){5}}  
\put(10,0){\line(1,2){12.5}} 
\put(0,0){\line(0,1){25}} 
\end{picture}}
}

\newcommand{\wti}[1]{\widetilde{#1}}
\newcommand{\tif}{\tilde{f}}
\newcommand{\tiL}[1]{\widetilde{L}}
\def\delti{\widetilde{\Delta}}
\def\es{\varepsilon}
\def\estil{\widetilde{S}}
\newcommand{\ellk}{\mbox{$\ell k$}}

\renewcommand{\bi}[2]{
\mbox{$
\left(
\begin{array}{@{}c@{}}
#1\\
#2
\end{array}
\right)$}}

\def\Theorem{Theorem  }

\def\Proposition{Proposition  }

\def\Corollary{Corollary  }
\def\Lemma{Lemma  }
\def\Figure{Figure  }
\def\Figures{Figures  }

\newtheorem{thm}{Theorem}[section]
\newtheorem{lem}[thm]{Lemma}%%%%Feb27

\newtheorem{prop}[thm]{Proposition}
\newtheorem{cor}[thm]{Corollary}
 \newtheorem{dfn}[thm]{Definition}
 \newtheorem{rem}[thm]{Remark}
 \newtheorem{ex}[thm]{Example}

 \title{Fibred torti-rational knots}
\author{Mikami Hirasawa}
\address{Department of Mathematics,
Nagoya Institute of Technology,
Nagoya Aichi 466-8555 Japan.\\
{\it E-mail: hirasawa.mikami@nitech.ac.jp}
}

\author{Kunio Murasugi}
\address{Department of Mathematics,
University of Toronto
Toronto, ON M5S2E4 Canada.\\
{\it E-mail: murasugi@math.toronto.edu}
}

\begin{document}
\maketitle

\begin{abstract}
A torti-rational knot, 
denoted by $K(2\alpha,\beta|r)$, 
is a knot obtained from the $2$-bridge link $B(2\alpha,\beta)$ 
by applying Dehn twists an arbitrary number of times, $r$, 
along one component of $B(2\alpha,\beta)$. 
We determine the genus of $K(2\alpha,\beta|r)$ and solve a 
question of when $K(2\alpha,\beta|r)$ is fibred.  
In most cases, the Alexander polynomials 
determine the genus and fibredness of these knots.
We develop both
algebraic and geometric techniques to describe the genus
and fibredness by means of continued fraction expansions
of $\beta/2\alpha$. 
Then, we explicitly construct minimal genus Seifert surfaces.
As an application, we solve the same question for 
the satellite knots of tunnel number one.
\end{abstract}


\section{Introduction}
A torti-rational knot \footnote[1]{
This naming is due to Lee Rudolph.}, denoted by $K(2\alpha,\beta|r)$, 
is a knot obtained from the $2$-bridge link $B(2\alpha,\beta)$ 
by applying Dehn twists an arbitrary number of times, $r$, 
along one component of $B(2\alpha,\beta)$. 
(For the precise definition, see Section 2.)

Torti-rational knots have occasionally
appeared in literatures of knot theory.  
For example, twist knots are 
torti-rational knots. By \cite{KMS}, we know when
$K(2\alpha, \beta|r)$ is unknotted (see Proposition \ref{prop:8.6++}).
A torti-rational knot is a $g1$-$b1$ knot 
(i.e., admits a genus-one bridge-one decomposition), and hence 
has tunnel number one.
In 1991, Morimoto and Sakuma \cite{MS}
proved that for a satellite knot of tunnel number
one, the companion knot is a torus knot $T(p,q)$
and the pattern knot is a torti-rational knot 
$K(2\alpha, \beta|pq)$, for some $p,q$, and $\alpha, \beta$.
Then,
Goda and Teragaito \cite{GT} determined which of such satellite knots
of tunnel number one are of genus one. 

In this paper, we study torti-rational knots systematically and
completely determine the genus of 
any torti-rational knot and solve a question of when it is fibred.  
%%If $r=0$, then $K(2\alpha, \beta|0)$ is off course a trivial knot,
%%we assume $r \ne 0$ unless specified otherwise.  

In fact, we prove the following:

%Theorem 1.1. 
%\begin{thm}\label{thm:torti}
{\bf Theorem A}. (Theorems \ref{thm:C} and \ref{thm:A})
{\it 
Let $B(2\alpha, \beta)$ be an oriented $2$-bridge link,
with linking number $\ell$. 
Let $K=K(2\alpha, \beta|r)$ be a torti-rational knot.  
  Suppose $\ell \ne 0$. \\
%   %Then we have:\\   
(1)	The genus of $K$ is exactly half of the degree of 
the Alexander polynomial $\Delta_K(t)$.\\
(2)	$K$ is fibred if (and only if) $\Delta_K(t)$
is monic (i.e, the leading coefficient 
%of $\Delta_K(t)$ 
is $\pm1$).
}%
%\end{thm}

See
Theorems \ref{thm:8.4} %Thm 6.1
and
\ref{thm:5.1} %Thm 9.1
for a practical method for determination.
 
\medskip

If $\ell = 0$, Theorem A does not hold true.  
For this case, the genus and the characterization 
of a fibred torti-rational
knot are stated as follows:

%Later we will give the complete characterization of fibred 
%knots in terms of 
%continued fractions. See \Theorem \ref{thm:5.1}.
%\begin{thm}\label{thm:D} %2.4

{\bf Theorem B.} (Theorems \ref{thm:D} and \ref{thm:B})
{\it
Suppose $\ell=0$,
Let
$[2c_1, 2c_2, \ldots, 2c_m]$
be the continued fraction of $\frac{\beta}{2\alpha}$.\\
(1) For any $r\neq 0$,
${\displaystyle
g(K(2\alpha, \beta|r))=
\dfrac{1}{2}\sum_{i:\ {\rm odd}} |c_i|
}$.\\
(2)
(a) If $|r|>1$, then $K(2\alpha, \beta|r)$ is not
fibred.
(b) Suppose $r=\pm 1$. Then
$K(2 \alpha, \beta|r)$ is fibred if and only if 
$\frac{\beta}{2\alpha}$ has the continued fraction of
the following special form:
$\frac{\beta}{2\alpha}=
\pm[2 a_1, 2, 2a_2, 2, \ldots, 2a_p, \pm 2, -2a'_1, -2,
-2a'_2, -2, \ldots, -2a'_q]$,
where $2\alpha>\beta>-2\alpha$, 
$a_i, a'_j>~0$ and $\sum_{i=1}^pa_i
=\sum_{j=1}^q a'_j$. 
}%\it

See Section 2, for our convention of continued fractions.
\medskip

To prove these theorems, we construct explicitly a 
minimal genus Seifert surface for 
$K$ and determine whether or 
not it is a fibre surface for $K$. 
Proofs of these theorems will be given in Sections 10 and 11.

%If $\ell = 0$, Theorem A does not hold true.  
%For this case, the genus and the characterization 
%of a fibred torti-rational
%knot will be stated in 
%Theorem \ref{thm:D} and Theorem \ref{thm:B} in Section 2.  

This work is a part of our project to determine
the genus and fibredness of double torus knots
(i.e., knots embedded in a standard closed surface of genus $2$).
See \cite{HirM} for a relevant work.
Double torus knots are classified into five types
(Type $(1,1), (1,2), (1,3), (2,3)$ and $(3,3)$).
In \cite{HirM}, we settled the problem for all double torus
knots of type $(1,1)$.
A $g1$-$b1$ knot can be presented as a double torus knot
of type either $(1, 2), (2, 2)$ or $(2,3)$. In this paper, we
settle the problem for the $g1$-$b1$ knots presented as of 
type $(1,2)$.
As an application of our study,
we determine the genus and the fibredness problem 
for the satellite knots of tunnel number one. 
In fact, we show that a similar theorem to 
Theorem A
holds true for satellite knots in a slightly wider class.  
The precise statements can be found in Section 13.

Recently, Goda, Hayashi and Song \cite{GHS}
study torti-rational knots
with a different motivation.
Their approach is completely different from ours.
%We note that the method we develop in this paper is
%also useful to prove some of their results 
%(Theorems 1.3 and 1.6), though our and their approaches
%are completely different.

This paper is organized as follows.
In Section 2, we give precise statements of our main
theorems (Theorems \ref{thm:C} - \ref{thm:B}).
In Section 3, we first introduce several
notions needed in this paper,
such as {\it graphs of continued fractions}, 
{\it dual graphs}. Then we prove that for our
study of $K(2\alpha, \beta|r)$, we may assume 
$\ell\ge0$ and $r>0$, where $\ell$ is the linking number
of $B(2\alpha, \beta)$. This restriction simplifies considerably
the proofs of our main theorems.
At the end of Section 3, we construct a 
spanning disk of a nice form 
for one component of the $2$-bridge link.
In Sections 4 and 5, 
we study the Alexander 
polynomial of various knots:
In Section 4, we determine the Alexander polynomials
of $K(2\alpha, \beta|r)$.
In Section 5, 
we prove one 
subtle property of the (2-variable) 
Alexander polynomial of
$B(2\alpha,\beta)$.  
The determination of the degree of the Alexander
polynomials of  $K(2\alpha, \beta|r)$ depends on
this property.
Sections 6 is devoted to
characterizing the monic 
Alexander polynomial of a knot $K(2\alpha,\beta|r)$:
First, we deal with knots $K(2\alpha, \beta|r)$
for the case $\ell >0$, and 
characterize the monic Alexander polynomials
in terms of a continued fraction of $\beta/2\alpha$
using the formulae given in Section 5.
In particular, we give an equivalent algebraic condition
for Theorem \ref{thm:A} (Theorem \ref{thm:8.4}).
However, for the case $\ell=0$, the monic Alexander 
polynomials of $K(2\alpha, \beta|r)$ cannot be characterized
by the continued fractions.
This case is considered in the rest of Section 6,
and the characterization will be done 
using a geometric interpretation of the Alexander polynomials
of $K(2\alpha, \beta|r)$.
In Section 7, we construct explicitly 
a Seifert surface for $K(2\alpha,\beta|r)$, 
which in most cases is of minimal genus.
In Section 8, we prove Theorem \ref{thm:C}.
In this case, some of the surfaces constructed in Section 7
are not of minimal genus, but we obtain  minimal
genus surfaces after explicitly compressing them.
In Section 9, we prove Theorem \ref{thm:A}.
In Sections 10 and 11, we prove Theorems \ref{thm:D}
and \ref{thm:B}. 
Various examples that illustrate our main theorems are 
discussed  in Section 12.
In section 13, we consider the satellite knot with fibred 
companion and 
$K(2\alpha,\beta|r)$, $r \ne 0$, as a pattern, and  
prove an analogous theorem to Theorem A.
In the final section, Section 14, 
we determine the genus one knots in our 
family of knots $K(2\alpha,\beta|r)$.
In particular, we find satellite knots among them, and
hence give a negative answer to the problem posed in
\cite{AM}.

%After our manuscripts are completed, we learned that
%Goda, Hayashi and Song also study torti-rational knots
%with completely different motivation in their paper \cite{GHS}.
%We note that some of their main results
%(Theorems 1.3 and 1.6) can be reproved by the method we developed
%in this paper, which are completely different from their approach.

After is this paper was completed, 
D. Silver pointed out that \Theorem 5.5 in this paper
makes it possible to prove Theorem A algebraically 
using Brown's graphs in \cite{brown}
(without explicit construction of minimal genus Seifert 
surfaces).
However, Brown's method does not work for proving Theorem B.

\section{Statements of main theorems}

We begin with an (even) continued fraction 
of a rational number $\frac{\beta}{2 \alpha},
0<\beta<2 \alpha$, and gcd$(2\alpha, \beta)=1$.
The (unique) continued fraction of
%\begin{minipage}{5cm}
$$\frac{\beta}{2 \alpha}
      =\cfrac{1}{2c_1
          -\cfrac{1}{2c_2
          -\cfrac{1}{2c_3
          -\cfrac{1}{\ddots
          -\cfrac{1}{2c_{m-1}
          -\cfrac{1}{2c_m}}}}}}, $$
%\end{minipage}
%\hfil
%\begin{minipage}{5cm}
where $c_i \neq 0$, is denoted by
$\frac{\beta}{2\alpha}=[2c_1, 2c_2,  \cdots, 2 c_m]$ 
or $[[c_1, c_2,  \cdots, c_m]]$.
Note that $m$ is odd.
Throughout this paper, we consider only even continued fraction
expansions, and hence omit the word \lq even\rq.
Now, using the continued fraction of 
$\frac{\beta}{2 \alpha}$,
we can obtain a diagram of
an oriented $2$-bridge link
$B(2\alpha, \beta)$ as
follows. 
%\end{minipage}
%\hfil
%\begin{minipage}{2cm}

%\Abb
%\end{minipage}
%\medskip

\begin{minipage}{8cm}
Let
$\sigma_1=$ \sone and $\sigma_2=$\stwo
be Artin's generators of the
$3$-braid group. First construct a $3$-braid
$\gamma=\sigma_2^{2c_1}
\sigma_1^{2c_2}
\sigma_2^{2c_3}\cdots \sigma_2^{2c_m}$.
Close $\gamma$ by joining the first and second strings 
(at the both ends)
and then join the top and bottom of the third string by
a simple arc as in \Figure 2.1.
We give downward orientation to the second and third
strings.
\Figure 2.1 shows the (oriented) $2$-bridge link obtained from
the continued fraction $\frac{21}{34}=
[2, 2, -2, -2, 2]=[[1,1,-1,-1,1]]$.
Now an oriented $2$-bridge link $B(2\alpha, \beta)$
consists of two unknotted knots $K_1$ and $K_2$,
where $K_2$ is formed from the third
and fourth strings.
\end{minipage}
\hfil
\begin{minipage}{2cm}
\Abb
\end{minipage}
\medskip

%%
%%
%%
%%Let
%%$\sigma_1=$ \sone and $\sigma_2=$\stwo
%%be Artin's generators of the
%%$3$-braid group. First construct a $3$-braid
%%$\gamma=\sigma_2^{2c_1}
%%\sigma_1^{2c_2}
%%\sigma_2^{2c_3}\cdots \sigma_2^{2c_m}$.
%%Close $\gamma$ by joining the first and second strings 
%%(at the both ends)
%%and then join the top and bottom of the third string by
%%a simple arc as in \Figure 2.1.
%%We give downward orientation to the second and third
%%strings.
%%\Figure 2.1 shows the (oriented) $2$-bridge link obtained from
%%the continued fraction $\frac{21}{34}=
%%[2, 2, -2, -2, 2]=[[1,1,-1,-1,1]]$.
%%Now an oriented $2$-bridge link $B(2\alpha, \beta)$
%%consists of two unknotted knots $K_1$ and $K_2$,
%%where $K_2$ is formed from the third
%%and fourth strings.

\medskip
{\bf Note.}
Our convention for the orientation
of a $2$-bridge link is not standard, but is used
for the convenience in utilizing the $2$-variable
Alexander polynomials. (Usually we reverse the orientation
of one component so that the $2$-bridge link is fibred if
and only if all the entries of the even continued fraction
are~$\pm~2$.)
\medskip

Since $K_2$ is unknotted, $K_1$ can be considered as 
a knot in an unknotted solid tours $V$ and $K_2$ 
is a meridian of $V$.
Then by applying Dehn twists along
$K_2$ in an arbitrary number of times , say $r$,
we obtain a new knot $K$ from $K_1$.
We denote this knot $K$ by $K(2\alpha, \beta|r)$,
or simply $K(r)$.

More precisely, one Dehn twist along $K_2$
is the operation that replaces the part of $K_1$ in a cylinder
by the braid
$(\sigma_1\sigma_2\cdots\sigma_{k-1})^{k}$, 
where $k$ is the 
wrapping number. See \Figure 2.2.
(Since $B(2\alpha, \beta)$ is symmetric, $K_1$ and $K_2$
can be interchanged, and hence this notation is justified.)

\Acc
%\Fi{2.3}

We note that if $r=0$, then $K(2\alpha, \beta|r)$
is unknotted for any $\alpha, \beta$, and henceforth we 
assume
$r\neq 0$ unless otherwise specified.

Now, given an oriented $2$-bridge link $B(2\alpha, \beta)$,
let $\ell=\ellk(K_1, K_2)$ be the linking number between
$K_1$ and $K_2$
which, for simplicity, is denoted by $\ellk B(2\alpha, \beta)$.

Our main theorems in this paper are the following four theorems:

\begin{thm}\label{thm:C} %2.3
Suppose $\ell\neq 0$. %Then for any $r\neq 0$, we have:
Then the genus of $K=K(2\alpha, \beta|r)$ is half of
the degree of its Alexander polynomial $\Delta_{K}(t)$. 
Namely we have
$g(K)=
\dfrac{1}{2}{\rm deg}\Delta_{K}(t)$.
\end{thm}

\begin{thm} \label{thm:A} %2.1
Suppose $\ell\neq0$. Then 
$K=K(2\alpha, \beta|r)$
is a fibred knot if (and only if) $\Delta_K(t)$
is monic, i.e., $\Delta_K(0)=\pm 1$.
\end{thm}

Theorem \ref{thm:A} is divided into two parts:
Theorem \ref{thm:8.4} is the algebraic part, 
where we determine when
$\Delta_K(t)$ is monic in terms of continued fractions,
and Theorem \ref{thm:5.1} is the geometric part, 
where we show the fibredness,
by actually constructing fibre surfaces using the
continued fractions.

%Later we will give the complete characterization of fibred 
%knots in terms of 
%continued fractions. See \Theorem \ref{thm:5.1}.
\begin{thm}\label{thm:D} %2.4
Suppose $\ell=0$,
Let
$[2c_1, 2c_2, \ldots, 2c_m]$
be the continued fraction of $\frac{\beta}{2\alpha}$.
Then for any $r\neq 0$,
${\displaystyle
g(K(2\alpha, \beta|r))=
\dfrac{1}{2}\sum_{i:\ {\rm odd}} |c_i|
}$.
\end{thm}

\begin{thm}\label{thm:B} %2.2
Suppose $\ell=0$. 
(a) If $|r|>1$, then $K(2\alpha, \beta|r)$ is not
fibred.
(b) Suppose $r=\pm 1$. Then
$K(2 \alpha, \beta|r)$ is fibred if and only if 
$\frac{\beta}{2\alpha}$ has the continued fraction of
the following special form:\\
$\frac{\beta}{2\alpha}=
\pm[2 a_1, 2, 2a_2, 2, \ldots, 2a_p, \pm 2, -2a'_1, -2,
-2a'_2, -2, \ldots, -2a'_q]$,
where $2\alpha>\beta>-2\alpha$, 
$a_i, a'_j>~0$ and $\sum_{i=1}^pa_i
=\sum_{j=1}^q a'_j$. 
\end{thm}

In Theorem \ref{thm:B}, we have non-fibred knots $K$
such that $\Delta_K(t)$ are 
monic and $\deg \Delta_K(t)=2 g(K)$.

\section{Preliminaries}%3
In this section, we first introduce two fundamental 
concepts, a {\it graph of a continued fraction}
and the {\it dual graph}, which
play a key role throughout this paper.
Next, in Subsection 3.4, we show that we can assume
$\ellkB\ge 0$ and $r>0$ without loss of
generality. This assumption is very important
to simplify the proof of our main theorem.
In the last Subsection (Subsection 3.5), 
we introduce the concept of the primitive
spanning disk for $K_1$, which is the 
first step of constructing a minimal genus 
Seifert surface for $K(2\alpha, \beta|r)$.

{\bf 3.1. Modified continued fractions and their graphs.}

Let $S=[[c_1, c_2, \ldots, c_{2k+1}]]$ be
the continued fraction of $\frac{\beta}{2\alpha}$,
where $-2\alpha<\beta<2\alpha$ and
gcd$(2\alpha, \beta)=1$.
The {\it length} of $S$ is defined as $2k+1$.
To define the dual of $S$, we need to extend $S$ slightly to
$S^*$, called the modified form of $S$.
We will see that $S$ and $S^*$ correspond to 
the same rational number.

\begin{dfn}\label{dfn:3.1}
Let $S=[[c_1, c_2, \ldots, c_{2k+1}]]$.
Then we obtain a continued fraction $S^*$
by thoroughly repeating the following and
call it {\it modified form of} $S$.\\
(1) If $c_{2i+1}>1, 0\le i\le k$, then $c_{2i+1}$ is replaced by
the  new sequence of length 
\hspace*{5mm} 
$2c_{2i+1}-1$,
$(1, 0, 1, 0, \ldots, 0,1)$, and\\
(2) if $c_{2i+1}<-1, 0\le i\le k$, then $c_{2i+1}$ is replaced
by the new sequence of length 
\hspace*{5mm}
$2|c_{i+1}|-1$,
$(-1, 0, -1, 0, \ldots, 0, -1)$.
Note that the length of $S^*$ is \\
\centerline{${\displaystyle
\sum_{i=0}^k (2|c_{2i+1}|-2)+2k+1
=\sum_{i=0}^k 2|c_{2i+1}|-1}$.
}\\
The original continued fraction $S$ may be called
the {\it standard} continued fraction of $\beta/2\alpha$,
which does not contain entries $0$.
\end{dfn}

The modified form of $\beta/2\alpha$ is of the form:
%$(3.1)
\begin{equation}
[2u_1, 2v_1, 2u_2, 2v_2, \ldots, 2u_d, 2v_d, 
2u_{d+1}],
\end{equation} 
where $u_i=+1$ or $-1$, for $1\le i\le d+1$, and
$v_i, (1\le i\le d)$ are arbitrary, including $0$.

Now, given the continued fraction $S$ of $\beta/2\alpha$,
consider the modified form $S^*$ for $S$ of the form (3.1).

\begin{dfn}\label{dfn:3.2}
The {\it graph} $G(S^*)$ of $S^*$,
(or the graph $G(S)$ of $S$), is a plane graph in 
${\mathbb R}^2$, consisting of $d+2$ vertices 
$V_0, V_1, \ldots, V_{d+1}$ and
$d+1$ line segments $E_k, (1\le k\le d+1)$ joining
two vertices $V_{k-1}$ and $V_k$, where
$V_0 = (0, 0)$ and $V_i=(i, \sum_{j=1}^i u_j)$, for
$1\le i\le d+1$.
The graph is a {\it weighted graph}, when the weight
of $V_i, (1\le i \le d)$ is defined as $2 v_i$. 
The weights of both $V_0$ and $V_{d+1}$ are $0$.
\end{dfn}
 
\Baa
%\Fi{3.1}

\begin{ex}\label{ex:3.1}
Let $S^*=[2, 0, 2, -2, -2, 0, -2]$.
Then $G(S^*)$ is depicted in \Figure 3.1.
The weight of $V_i$ is denoted by $(m)$ near $V_i$.
\end{ex}

The following is
immediate from the diagram of $B(2\alpha, \beta)$
(Figure 2.1).

\begin{prop}\label{prop:3.3}
The $y$-coordinate of the last vertex
$V_{d+1}$ gives the linking number $\ell=
\ellkB$. Namely, 
${\displaystyle \ell=\sum_{i=1}^{d+1} u_i}$.
\end{prop}

{\bf 3.2. Dual graphs and dual continued fractions.}

For a continued fraction $S$,
we define the dual $\wti{S}$ to $S$
and 
the dual graph $\wti{G}$ to a graph $G(S)$.
Then we have the following theorem, 
proved in Subsection 3.5:

\begin{thm}\label{thm:3.6} %thm 3.9 (Apr23) %3.5 Feb23 2006.
Let $S$ be the continued fraction of $\beta/2\alpha$.
Then the dual $\estil$ of $S$ is the continued fraction of 
$(2\alpha-\beta)/2\alpha$
(resp. $(-2\alpha-\beta)/2\alpha$),
if $\beta>0$ (resp. $\beta<0$).
\end{thm}

\begin{dfn}\label{dfn:3.4}
The {\it dual} $\wti{G}$ to the graph $G(S)$
is defined as follows. The underlying graph of 
$\wti{G}$ is exactly the same as that of $G(S)$,
but the weight $\wti{w}(V_i)$ is given as follows;\\
(1) If $V_i$ is a local maximal or local minimal vertex
(including the ends of $G(S)$), then 
$\wti{w}(V_i)=-w(V_i)$, and\\
%(2) if $V_i$ is an end vertex, then 
%$\wti{w}(V_i) = 0$, and 
(2) for the other vertices, $\wti{w}(V_i)= 2 \es_i-w(V_i)$,
where $\es_i$ is the sign of $u_i$, i.e., $\es_i=u_i/|u_i|$

The dual $\wti{S}^*$ to the modified form $S^*$ is defined to be
the modified form of the continued fraction represented
by the dual graph $\wti{G}$.
The dual $\wti{S}$ of $S$ is the standard continued fraction
obtained from $\wti{S}^*$.
\end{dfn}

\Bb
%\Fi{3.2}

\begin{ex}\label{ex:3.2} %ex3.6 (Apr23)
Let $S=[[2,-1,-1,1,-1]]$. \\
Then
$S^*=[[1,0,1,-1,-1,1,-1]]=[4,-2,-2,2,-2]$.\\
Thus 
$\widetilde{S}^*=[[1,1,1,1,-1,-2,-1]]=
[2,2,2,2,-2,-4,-2]=\wti{S}$.
\end{ex}

In the following, we give an alternative formulation
of $\widetilde{S}$, the dual of $S$.
Given the continued fraction of $\beta/2\alpha$,
%$(3.2) 
\begin{equation}
[[c_1,c_2, \ldots, c_{2d+1}]], c_i\neq 0,
\end{equation}
consider the
partial sequence of (3.2):
%$(3.3)
\begin{equation}
\{c_1, c_3, c_5, \ldots, c_{2d+1}\},
\end{equation}
consisting of only $c_{2i+1}, 0\le i\le d$.

In this sequence, for convenience, write $-c_i$,
where $c_i<0$, so that we may assume that
if $i$ is odd, then $c_i$ is always positive.
Thus, the sequence (3.3) is divided into several
\lq positive' or \lq negative' sub-sequences:
Therefore we can write, 
$[[c_1, c_2, \ldots, c_{2d+1}]]=$
$[[a_1, b_1, a_2,$
$ \ldots, a_p, b_p,$
$ -a_{p+1}, -b_{p+1},
\ldots, -a_{r}, b_r, a_{r+1}, b_{r+1}, \ldots]]$,
where $a_i>0$ for all $i$.
Note that $a_i=c_{2i-1}$ or $-c_{2i-1}, 
i=1, 2, \ldots$ and $b_j=c_{2j}, j=1, 2, \ldots$

We call a sequence of the form
$[[a_1, b_1, a_2, b_2, \ldots, a_k, b_k, a_{k+1}]]$\\
(or $[[-a_1, -b_1, -a_2, -b_2, \ldots, -a_k, -b_k, -a_{k+1}]]$)
a {\it positive} (or {\it negative}) sequence, where
$a_i>0,
1\le i \le k+1$, but $b_j, (1\le j\le k)$ are arbitrary ($\neq 0)$.
We denote by $P_i$ (resp. $Q_i$) a positive (resp. negative)
subsequence.

\begin{ex}\label{ex:3.3}
\begin{align*}
&\ \ \ [[1, 1, 2, -1, 1, -1, -2, 1, -2, -1, 2, 1, 2, 1, -2, -1, -2]]\\
&=[[\underbrace{1, 1, 2, -1, 1}_{P_1}, -1 ,
\underbrace{-2, -(-1), -2}_{Q_1}, -1,
\underbrace{2, 1, 2}_{P_2}, 1,
\underbrace{-2, -1, -2}_{Q_2}]]\\
&=[[P_1,-1,Q_1,-1,P_2,1,Q_2]]
\end{align*}
\end{ex}

Thus, the sequence (3.2) can be written as
%(3.4)
\begin{equation}
[[c_1, c_2, c_3, \ldots, c_{2d+1}]]=
\{P_1, d_1, Q_1, e_1, P_2, d_2, Q_2, e_2,P_3, \ldots\}, 
\end{equation} 
where $d_i,e_j$ are some $c_{2k}$.

This form (3.4) is called the {\it canonical decomposition} 
of the continued fraction of $\beta/2\alpha$.

\begin{rem}\label{rem:3.5}
If $\beta>0$, then the first entry $c_1>0$, but if $\beta<0$, 
then $c_1<0$, and hence, the canonical decomposition begins 
with $Q_1$ (not a positive sequence $P_1$ and $d_1$ is missing).
However, since this does not change our argument, 
we may assume in general that $c_1>0$.
\end{rem}

Now, the {\it dual continued fraction} of (3.2) 
is reformulated as follows.

Let $S=\{P_1, d_1, Q_1, e_1, P_2, \ldots\}$
 be the canonical decomposition of 
 $[[c_1, c_2, \ldots, c_{2d+1}]]$.

First the {\it dual of a positive sequence} $P$ 
is obtained as follows.

Given $P=[[a_1, b_1, a_2, b_2, \ldots a_{m+1}]], a_j>0$, 
consider the modified form $P^*$ of $P$
%(3.5)\hfil
\begin{equation}
P^*=[[a^*_1, b_1^*, a_2^*, b_2^*, \ldots, 
a_k^*, b_k^*, a^*_{k+1}]],
\end{equation}
where $a^*_j=1 (1\le j\le k+1)$ and $b^*_j$'s $(1\le j\le k)$ 
are arbitrary including $0$.

Then the dual of $P^*$, denoted by $\wti{P^*}$, 
is $\wti{P^*}=
[[\wti{a}_1^*, \wti{b}_1^*, \wti{a}_2^*, \wti{b}_2^*, 
\ldots,\wti{a}_{k+1}^*]$,
where $\wti{a}_j^*=a_j^*=1$, for 
$1\le j\le k+1$ and $\wti{b}_j^*=1-b_j^*$ for $1\le j\le k$.\\
The {\it dual $\wti{P}$ of $P$} 
is the standard form obtained
from $\wti{P}^*$.

For the negative sequence $Q$, 
apply the same operation for the positive sequence $-Q$ 
to obtain the dual $\wti{-Q}$ of $-Q$.
Then the dual $\wti{Q}$ of $Q$ 
is the negative sequence $-(\wti{-Q})$.

Finally, the {\it dual} $\wti{S}$ of $S$ is
$\{\wti{P}_1,-d_1,\wti{Q}_1,-e_1,\wti{P}_2,-d_2,
\wti{Q}_2,-e_2, \ldots\}$.

{\bf Example} \ref{ex:3.3} {\bf (continued)}

(1)Since $P_1=[[1,1,2,-1,1]]$,  $P^*_1=[[1,1,1,0,1,-1,1]]$,
and hence\\
$\wti{P}^*_1=[[1,0,1,1,1,2,1]]$ and
$\wti{P}_1=[[2,1,1,2,1]]$.\\
(2) Since $P_2=[[2,1,2]]$,
$P^*_2=[[1,0,1,1,1,0,1]]$,
and hence\\
$\wti{P}^*_2=[[1,1,1,0,1,1,1]]$ and
$\wti{P}_2=[[1,1,2,1,1]]$.\\
(3)
Since $Q_1=[[-2,1,-2]], -Q_1=[[2,-1,2]]$ and hence\\
$(-Q_1)^*=[[1,0,1,-1,1,0,1]]$,\\
$\wti{-Q_1^*}=[[1,1,1,2,1,1,1]]=\wti{-Q_1}$, so
$\wti{Q}_1=[[-1,-1,-1,-2,-1,-1,-1]]$\\
(4)
Since $Q_2=[[-2,-1,-2]], -Q_2=[[2,1,2]]$ and hence\\
$\wti{-Q_2^*}=[[1,1,2,1,1]]$, so
$\wti{Q}_2=[[-1,-1,-2,-1,-1]]$.
Thus, 
$\estil=[[2,1,1,2,1,1,-1,-1,-1,-2,-1,-1,-1,1,1,1,2,1,1,-1,
-1,-1,-2,-1,-1]]$

{\bf 3.3. Applications.}

In this subsection, we study
some of the invariants of a $2$-bridge link 
$B(2\alpha, \beta)$ deduced from $S$ or its dual $\estil$.

The following three propositions show that $S$ or $\estil$ 
determines the degree of the Alexander polynomial 
$\Delta_{B(2\alpha, \beta)} (x, y)$ of $B(2\alpha, \beta)$.

Let $S=\{P_1, d_1, Q_1, e_1, P_2, \ldots, P_m, d_m, Q_m\}$ 
be the canonical decomposition of the continued fraction of 
$\beta/2\alpha$. 
We write more precisely:
%(3.6)  
\begin{align}
P_i&=[[a_{i,1}, b_{i,1}, a_{i,2}, b_{i,2}, \ldots, 
b_{i,s_i}, a_{i, s_{i}+1}]],\ {\rm and} \nonumber\\
Q_j&=[[-a'_{j,1},-b'_{j,1},-a'_{j,2},-b'_{j,2}, \ldots,
-b'_{j,q_j},-a'_{j,q_{j}+1}]]
\end{align}

\begin{dfn}\label{dfn:3.8}
For $1\le i,j \le m$, define
%(3.7)\hfil 
\begin{align}
\rho_i&=|\{b_{i,\ell}|b_{i, \ell}=1, 1\le \ell \le s_i\}|,\nonumber\\
\rho'_j&=|\{b'_{j,\ell}|b'_{j,\ell}=1,1\le \ell \le q_j\}|,\  {\rm and}\nonumber\\
\rho&=\rho(\beta/2\alpha)=\sum_{i=1}^{m} 
\rho_i +\sum_{j=1}^{m}\rho_{j'}.
\end{align}

We call this $\rho$ the {\it deficiency} (see Theorem \ref{thm:7.4}
and Sections 6 and 8).

Further, we define;
%(3.8) \hfil 
\begin{align}
\lambda_i&=\sum_{\ell=1}^{s_{i}+1}a_{i, \ell}, 
1\le i\le m \nonumber\\
\lambda'_j&=\sum_{\ell=1}^{q_{j}+1}a'_{j,\ell}, 1\le j\le m,\ {\rm and}\nonumber\\
\lambda&=\sum_{i=1}^{m}\lambda_i +\sum_{j=1}^{m}\lambda'_j.
\end{align}
\end{dfn}

Note that $\lambda$ equals the number of edges in $G(S)$,
which also equals the number of disks in Figure 3.4.
This number is neatly evaluated by Kanenobu as follows :

%\begin{rem}\label{kanenobu}
%The wrapping number of $B(2\alpha, \beta)=K_1\cup K_2$, denoted by 
%$w(2\alpha, \beta)$, is defined as the minimal number of 
%intersection points of $K_1$ and a disk spanned by $K_2$. 
%This number is neatly evaluated by Kanenobu \cite[(4.10)]{K} as follows :
%
%The wrapping number $w(2\alpha, \beta)$ is equal to the 
%number of edges of $G(S)$.
%\end{rem}

\begin{prop}\label{prop:3.9} %[KM] 
\cite[(4.10)]{K}
Write $\Delta_{B(2\alpha, \beta)}(x,y)=
{\displaystyle \sum_{0\le i, j}} c_{i,j} x^i y^j
\in {\mathbb Z}[x,y]$
in such a way that $\min y$-deg $\Delta_{B(2\alpha, \beta)}
(x,y) = \min\{j|c_{i,j}\neq 0\}=0$.
Then $\max y$-deg $\Delta_{B(2\alpha, \beta)}(x,y)=
\max\{j|c_{i,j}\neq 0\}=\lambda -1$.
\end{prop}

The following proposition shows that $\lambda$ and $\rho$
are related to the dual of $S$.

\begin{prop}\label{prop:3.10}
Let $\estil$ be the dual of $S$. Then the length of $\estil$
is $2(\lambda-\rho)-1$.
\end{prop}

{\it Proof.}
First consider the positive sequence $P_i$.

Let $P_i=[
2a_{i,1},2b_{i,1},2a_{i,2},\ldots,2a_{i,s_i},2b_{i,s_i},
2a_{i,s_{i}+1}]$.
Then the length of $P_i$ is $2s_{i}+1$.
Now to get the dual, consider the modified form 
$P^*_i$ of $P_i$
that is of the form:
%\begin{equation*}
\hfil
$P^*_i=[\underbrace{2, 0, 2, \ldots, 0, 2}_{2a_{i,1}-1}, 
2b_{i, 1}, \underbrace{2, 0, 2, \ldots, 0, 2}_{
2a_{i,2}-1}, 2b_{i, 2}, \ldots].
$
%\end{equation*}

Then to obtain $\wti{P}_i$, replace $0$ in $P^*_i$
by $2$ and $b_{i,r}$ by $1-b_{i,r}$.
Therefore, in $\wti{P}_i$, $0$ occurs exactly $\rho_i$ times.
Since the length of $P^*_i$ is $\sum^{s_i+1}_{r=1} 
(2a_{i, r}-1)+s_i$, the length of the dual $\wti{P}_i$ is
\begin{align*}
\sum^{s_i+1}_{r=1} (2a_{i, r}-1)+s_i -2 \rho_i
&= \sum^{s_i+1}_{r=1} 2a_{i, r}-(s_i +1)+ s_i -2\rho_i\\
&= 2\lambda_i-2\rho_i-1.
\end{align*}

By the same reasoning, the length of the dual 
$\wti{Q}_j$ of $Q_j$ is equal to $2\lambda'_j-2\rho'_j-1$.

Therefore, the length of the dual $\estil$ of $S$ is 
\begin{equation*}
\sum^m_{i=1} (2\lambda_i-2\rho_i-1)+\sum^m_{j=1} 
(2\lambda'_j-2\rho'_j-1)+(2m-1)=2\lambda-2\rho-1.
\end{equation*}
\qed

Note that \Proposition \ref{prop:3.10} 
holds if $Q_m=\phi$ (and hence $d_m$ is missing) or $P_1=\phi$ 
(and hence $d_1$ is missing).

 Combining \Proposition 2 in \cite{KM} with 
 \Proposition \ref{prop:3.10} and
 \Proposition \ref{prop:4.4} (in the next subsection), 
 we obtain:

\begin{prop}\label{prop:3.11}%3.14 
Let $\Delta_{B(2\alpha, \beta)} (x, y)$ be the 
Alexander polynomial of a 2-bridge link $B(2\alpha, \beta)$. 
Then $\Delta_{B(2\alpha, \beta)} (t, t)$ is a polynomial of 
degree $2(\lambda-\rho-1)$.
\end{prop}

{\it Proof.}
Let $\Delta_{B}(t)$ be the reduced Alexander polynomial of 
$B(2\alpha, \beta)$. 
Then $\Delta_{B}(t)=\Delta_{B(2\alpha, \beta)} (t, t)(1-t)$.
Apply \Proposition 2 in \cite{KM}.
\qed
 
%%\section{Reduction} %4
\medskip
{\bf 3.4. Reduction.}

In this subsection, we justify
the assumptions $\ell\ge 0$ and $r>0$.
This restriction drastically 
simplifies the proofs of
the main theorems.

For a knot $K$, we denote by $\overline{K}$
the mirror image of $K$.

\begin{thm}\label{thm:1} %4.1
In studying the genera and fibredness of 
$K(2\alpha, \beta|r)$
with $r\neq0$,
we may assume:\\
(1) $-2\alpha<\beta<2\alpha$,
(2) $\ellk B(2\alpha, \beta)\ge 0$, and
(3) $r>0$.

More precisely, we have the following:
Suppose $0<\beta<2\alpha$. Then, for any $r\neq 0$\\
(I) \ 
$K(2\alpha, \beta|r)=\overline{K(2\alpha, -\beta|-r)}$.\\
(II)
$K(2\alpha, \beta|r)=
\overline{K(2\alpha, 2\alpha-\beta|-r)}$ and
$K(2\alpha, -\beta|r)=
\overline{K(2\alpha, -2\alpha+\beta|-r)}$.
\end{thm}

We remark the following:
(i) $\ellkB=
-\ellk B(2\alpha, -\beta)$,
(ii) $\ellkB=
\ellk B(2\alpha, \pm 2\alpha-\beta)$,
(iii) if $\ell B(2\alpha, \beta)=0$, then we may assume
$r>0$ and $\beta>0$.

\begin{ex} %\label{ex:4.8}
(1) $K(4,-1|3) = \overline{K(4,1|-3)} = K(4,3|3)$.\\
(2) $K(14,9|2)= \overline{K(14,-9|-2)}=K(14,-5|2)$.\\
(3) $K(14,9|-2)=\overline{K(14,-9|2)}$.\\
(4) $K(14,-9|-2)=\overline{K(14,-5|2)}$.\\
Note that $\ellk B(4,-1)=-2$ and $\ellk B(14,9)=-1$. 
So, we first take the mirror
image to make the linking numbers positive,
at the expense of changing the sign of~$r$.
\end{ex}

{\it Proof of \Theorem \ref{thm:1}.}
First, we prove the former part, namely we show;

\begin{prop}
We may assume (1), (2) and (3) of \Theorem \ref{thm:1}.
\end{prop}

{\it Proof.}
Take a $2$-bridge link $B(2\alpha, \beta)=K_1\cup K_2$.
Without loss of generality, we may assume 
$-2\alpha<\beta<2\alpha$.
If $\ellk(K_1,K_2)<0$, take
 $B(2\alpha, -\beta)$. %:=\bar{K}_1\cup \bar{K}_2$. 
 This corresponds to
 taking the mirror image of $B(2\alpha, \beta)$
 while preserving the orientation
 of the components. 
 Therefore, $\ellk B(2\alpha, -\beta)>0$,
 and hence, from now on, 
 we assume $2$-bridge links always have
 a non-negative linking number. Note that we still have
 $-2\alpha<-\beta<2\alpha$.
 Take $K(2\alpha, \beta|r)$. 
 Suppose $r<0$.
 Let  $B(2\alpha',\beta')=K'_1 \cup K'_2$
 be the link obtained by
taking the mirror
 image of $B(2\alpha, \beta)$
 while reversing the orientation of
 $K_2$. 
 Now the linking number is preserved,
 i.e., $\ellkB=\ellk
 B(2\alpha',\beta')$.
 Recall that $K(2\alpha, \beta|r)$ is obtained
 by twisting $K_1$ by $K_2$, $r$ times. 
 This does not depend on
 the orientation of $K_2$, 
 and hence
the knot obtained by twisting $K_1$ along $K_2$ $r$ times
is the mirror image of the knot obtained by twisting 
$K_{1}'$ along $K_{2}'$ $-r$ times.
 Therefore, we see 
 $K(2\alpha, \beta|r)=\overline{K(2\alpha', \beta'|-r)}$
\qed

\medskip

%Next we determine the  $2$-bridge link $L'$
%that is 
%the mirror image of $L$ with one 
%component reversed.
Next, we show the following to prove the latter half
of Theorem \ref{thm:1}.

\begin{prop}\label{prop:4.4}
Let $L=B(2\alpha,\beta)$ be a $2$-bridge link,
where $-2\alpha<\beta<2\alpha$.
Let $L'$ be obtained by taking the mirror image of $L$
while reversing the orientation of one component.
Then, we have:
$L'= 
\begin{cases}
\begin{array}{ll}
B(2\alpha, 2\alpha-\beta) & {\rm if\ } \beta>0,\\
B(2\alpha, -2\alpha-\beta) & {\rm if\ } \beta<0.
\end{array}
\end{cases}
$
\end{prop}

{\it Proof.}
Since the other case is similar,
we only deal with the case $\beta>0$. 
Consider the Schubert normal form of $B(2\alpha, \beta)$.
See \Figure 3.3.

\foA

First, take the mirror image by changing all crossings
simultaneously. Flip the figure by the horizontal axis.
Now we have the Schubert normal form of 
$B(2\alpha, -\beta)$. Rotate the right over-bridge 
clockwise by $\pi$. Change the orientation of
the component containing the right over-bridge.
This gives the Schubert normal form of
$B(2\alpha, 2\alpha -\beta)$.
\qed

%\embfig{80}{4-1.eps}{Figure 4.1:
%Deformation from $B(2\alpha,\beta)$ to 
%$B(2\alpha, 2\alpha-\beta)$, e.g. $B(8,3)$
%to $B(8,5)$
%}

By two propositions above, we have
\Theorem \ref{thm:1}.
\qed

\medskip

\noindent
\begin{minipage}{7.5cm}
{\bf 3.5. Primitive spanning disk for $K_1$.}

In this subsection, we introduce the notion of
{\it primitive spanning disk}
 for $K_1$, which locally look like Figure 3.5 (b).
This surface is the first step to
construct a minimal genus Seifert surface
for $K(2\alpha, \beta|r)$.
Let $D$ be
a diagram obtained from the continued fraction
$S$ of $B(2\alpha, \beta)$ as in Figure 2.1.
By a slight modification of $D$ corresponding to
the modification of $S$ to $S^*$, as in \Figure 3.4,
construct a spanning disk for $K_1$,
which consists of horizontal disks and vertical bands,
whose interiors are mutually disjoint.
In \Figure 3.4, each box contains an even number of
twists (including $0$). 
%%We call this Seifert surface a {\it primitive disk}
%%for $K_1$. 
Note that
in \Figure 3.4, all disks are showing the same side,
though $K_2$ may penetrate them from various sides.
The set of horizontal disks is divided into several 
families so that each member of a family
meets $K_2$ from the same side as its neighbouring
member(s).
This corresponds to the canonical
decomposition $\{P_1, d_1, Q_1, e_1, \cdots\}$
of $S$.
For simplicity,
the disks belonging to the family corresponding
to $P_i$'s (resp. $Q_i$'s) are called
{\it positive disks} (resp. {\it negative disks}),
and a band connecting  two positive (resp. negative)
disks is called a 
{\it positive} (resp. {\it negative}) band.
The other bands are called {\it connecting bands}.

\end{minipage}
%\foB
\hfil
\begin{minipage}{4cm}
\medskip
\medskip
\medskip
\embfig{80}{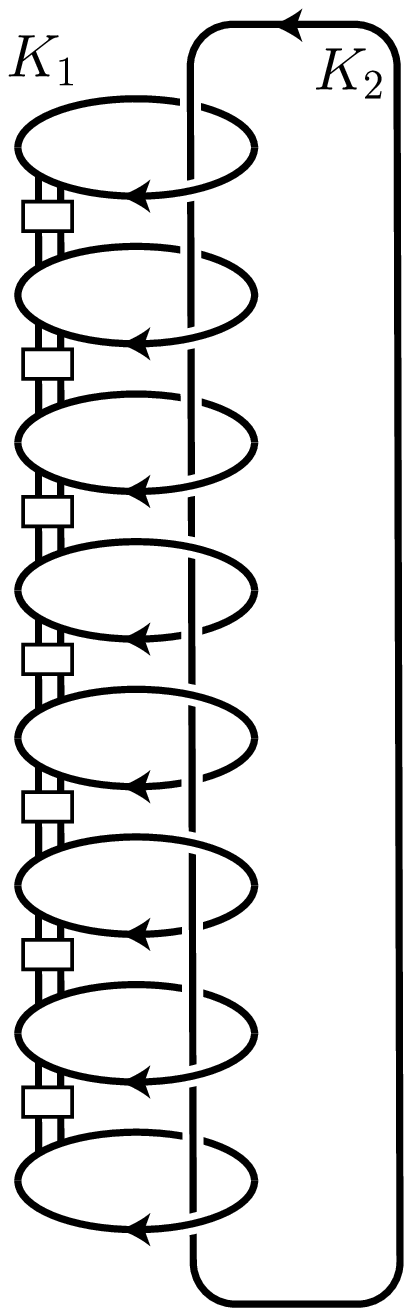}{Figure 3.4:
$B(2\alpha, \beta)$
}
\end{minipage}

\begin{rem}\label{rem:4.1}
 Disks and bands in the spanning disk for $K_1$
 correspond to edges and vertices, respectively, 
 of the graph $G(S)$ of $S$ as follows,
 except for the end vertices of $G(S)$.\\
%\begin{tabular}{ll}
(1) Positive/negative disks correspond to 
edges with positive/negative slope.\\
(2) Positive/negative bands correspond to vertices
between positive/negative edges. \\
(3) Connecting bands correspond to local maximal or 
minimal vertices\\
(4) The number of twists of a band corresponds to
the weight of a vertex.
%\end{tabular}
\end{rem}

\begin{dfn}\label{dfn:primitivedisk}
Slide each band so that both of its ends are 
attached to the front edge of each small disk
as in \Figure 3.5 (b).
The {\it primitive spanning disk} for $K_1$ is
the union of all the small disks 
together with all bands arranged this way.
See Figures 7.2 (a) and 7.4 left.
\end{dfn}

We remark that in the process of sliding a band,
another band may stand in the way. However,
as shown in the following proposition, we
can always arrange the bands so that each of them
appears as in \Figure 3.5 (b).

\foC
%\embfig{99}{4-3.eps}{Figure 4.3:
%Deformation from $B(2\alpha,\beta)$ to 
%$B(2\alpha, 2\alpha-\beta)$.
%}\medskip

\begin{prop}\label{prop:4.2}
 A relative position
 of the bands in a primitive spanning disk can be
 arbitrary.
 \end{prop}

{\it Proof.}
Examining the case locally suffices.
See \Figure 3.6, where disks, say $D_1, D_2, D_3$
and bands $B_1, B_2, B_3$ are depicted.
To change from (a) to (b), fix $D_2$
and everything
lying above $D_2$, and simultaneously 
turn around everything
that hangs below $D_2$.
Similarly we can change (b) to (c).
\qed
%To change from (b) to (c), fix $D_3$
%and everything 
%lying above $D_3$ and simultaneously move everything
%that hang below $D_3$ clockwise.

\foD

Now, we demonstrate the process of replacing
$B(2\alpha, \beta) =K_1 \cup K_2$ by
$B(2\alpha, 2\alpha-\beta)$, that is
to take the mirror image and reverse
(the orientation of) $K_2$:
Since $K_2$ consecutively penetrates the disks
transversely, we have a diagram of $K_1$ as in 
Figure 3.7, where (i) all the disks are concentric,
(ii) the higher disk appears smaller
and (iii) the only crossings are in the twists of bands.
 \Figure 3.7 shows the process of taking the mirror
image and reversing the mirror image of $K_2$.

\foE
%\embfig{80}{4-5.eps}{Figure 4.5:
%Reflect $K_1\cup K_2$ and reverse $K_2$
%}\medskip

Then we notice that the effect of the process is 
simply replacing each of the bands in $K_1$ by 
its mirror image.
Therefore, the process can be depicted as in
\Figure 3.5 (b) to (c).
Now reversing the operation of (a) to (b), we obtain
the standard diagram of the  $2$-bridge link
$B(2\alpha, 2\alpha-\beta)$.

Finally, \Theorem \ref{thm:3.6} is now almost immediate.

{\it Proof of Theorem} \ref{thm:3.6}. 
It is easy to see that the final diagram \Figure 3.5(d) 
is the primitive disk obtained from the dual 
$\widetilde{S}$ of $S$. 
Therefore \Theorem \ref{thm:3.6} follows from 
\Proposition \ref{prop:4.4}.
\qed

\section{Alexander polynomials (I)}%6 Feb 25  Sec 4
In this section, we determine the Alexander polynomial 
$\Delta_{K(r)}(t)$ for $K(2\alpha, \beta|r)$. 
In fact, we prove the following

\begin{prop}\label{prop:6.1}
 Let $\Delta_{B(2\alpha, \beta)} (x, y)$ be the 
 Alexander polynomial of an (oriented) $2$-bridge link 
 $B(2\alpha, \beta)$.
Let $\Delta_{K(r)}(t)$ be the Alexander polynomial of 
$K(2\alpha, \beta|r), r>~0$.\\
(1)\cite{Ki} If $\ellkB=\ell\ne0$, then
$\Delta_{K(r)}(t)=\dfrac{1-t}{1-t^\ell} 
\Delta_{B(2\alpha, \beta)} (t,t^{\ell r})$\\
(2) If $\ell=0$, then, for some $a=\pm 1$ and $b$,
\begin{equation*}
\Delta_{K(r)}(t)=r \biggl[
   \dfrac{\Delta_{B(2\alpha, \beta)}(x, y)}
   {1-y}\biggr]_{x=t, y=1} (1-t)+a t^b
\end{equation*}
\end{prop}

In Subsection 6.2, %%%
we will give a geometric interpretation of 
$\Delta_{K(r)}(t)$ when $\ell=0$. (See also \cite{Gon}
 or \cite{MM})

Now \Proposition \ref{prop:6.1} (1) 
follows from a general result due  
to Kidwell \cite{Ki}, and hence we omit the proof. 
However, part(2) was not proved in \cite{Ki}.
In this section, we prove the following more general
result suggested by M. Kidwell.

\begin{prop}\label{prop:6.2}
 Let $K_1$ be an oriented knot embedded in an 
 (unknotted) solid torus $V$. 
 Suppose $\ellk (K_1, K_2)=0$, where
 $K_2$ is an oriented meridian of $\partial V$.
 Denote by $K_1(r)$ the knot obtained from $K_1$
 by applying Dehn twists $r$ times along $K_2$ ($r>0$).
Let $L=K_1\cup K_2$. Then, for some $a=\pm 1$ and $b$, we have:
%
%(4.1) %%%
\begin{equation}
\Delta_{K_1(r)}(t)=r(1-t) \biggl[
\dfrac{\Delta_{L}(x, y)}{1-y}\biggr]_{
\genfrac{}{}{0pt}{}{x=t}{y=1}}+
a t^b   
\Delta_{K_1}(t).
\end{equation} 
\end{prop}

\Proposition \ref{prop:6.1} (2) 
follows from \Proposition \ref{prop:6.2} immediately, 
since $\Delta_{K_1}(t)=1$.

{\it Proof of \Proposition \ref{prop:6.2}.}
First, consider the link $L=K_1\cup K_2$. 
We add one trivial knot $K_3$ to $L$ such that 
$\ellk(K_1, K_3)=\ellk(K_2, K_3)=1$
as in \Figure 4.1. %%%
Let $\wti{L}=K_1\cup K_2\cup K_3$
be the $3$-component link. 

\sixa %\Fi{6.1}  %%4.1

Using this diagram, we obtain the following Wirtinger 
presentation of the link group 
$G(\wti{L})$ of $\wti{L}$.
\begin{equation*}
G(\wti{L})=\langle
x_1, x_2, \ldots, x_m, y, z_1, z_2|r_1, \ldots, r_m, s, t_1, t_2
\rangle,\ {\rm where}
\end{equation*}
\hspace*{4mm}
$
\begin{array}{ll}
r_1=z_1  x_1  z_1^{-1}  x_2^{-1},
&
\hspace*{8mm}
s=(x_{i_{k}}^{\varepsilon_k}\cdots 
x_{i_2}^{\varepsilon_2} x_2  z_1) y (z_1^{-1}  x_2^{-1}  
x_{i_2}^{-\varepsilon_2}\cdots x_{i_k}^{-\varepsilon_k})y^{-1},\\ 
r_2=y  x_2 y^{-1}  x_3^{-1},
&
\hspace*{8mm} t_1=y z_1 y^{-1} z_2^{-1},\\
r_3=w_3  x_3  w_3^{-1}  x_4^{-1},
&
\hspace*{8mm} t_2=x_1 z_2 x_1^{-1} z_1^{-1}\\
\hspace*{8mm}\vdots & \\                                 
r_m=w_m  x_m  w_m^{-1}  x_1^{-1}. & 
\end{array}$\\
Here, $w_i$ is a word in $x_i$ and/or $y$.

We note that $\varepsilon_k+\cdots+\varepsilon_2+1=0$, 
since $\ell k(K_1, K_2)=0$.

Now using $t_1$ and $t_2$, we can eliminate $z_2$ and 
obtain a new presentation.
For simplicity, we write $z=z_1$. Then
\begin{equation*}
G(\wti{L})=\langle
x_1, x_2, \ldots, x_m, y, z|r_1, \ldots, r_m, s, t'\rangle,\ {\rm where}\ 
t'=x_1 y z y^{-1} x_1^{-1} z^{-1}.
\end{equation*}

From this presentation, 
we obtain the Alexander matrix 
$M(\wti{L}$) for $\wti{L}$. 
The matrix $M(\wti{L})$ is an 
$(m+2)\times(m+2)$ matrix. 
A simple calculation shows 

%(4.2)\begin{tabular}{ll} %%%
\begin{align}
&(1)\ 
(\dfrac{\partial r_i}{\partial y})^\phi=
\delta_i(1-x),\  {\rm where}\ \delta_i=0, 1\ {\rm  or}\ -y^{-1},\nonumber\\
&(2)\
(\dfrac{\partial r_1}{\partial z})^\phi=1-x, 
(\dfrac{\partial r_i}{\partial z})^\phi=0,\ {\rm for}\ i \ne1,
\end{align}
where $\partial$  indicates Fox's 
free derivative and $\phi$ is the induced homomorphism 
from $G(\wti{L})$ to the free abelian group 
$G(\wti{L})/[G(\wti{L}),G(\wti{L})]$, where 
$x_i^\phi=x, y^\phi=y$ and $z^\phi=z$.
Let $U=x_{i_k}^{\varepsilon_k} 
\cdots x_{i_2}^{\varepsilon_2} x_{2}z$.
Then $s=U y U^{-1} y^{-1}$, and
%(4.3)  %%%
\begin{equation}(\dfrac{\partial s}{\partial x_i})^\phi=
(1-y) (\dfrac{\partial U}{\partial x_i})^\phi.
\end{equation}

Since $U$ does not involve $y$ and $\varepsilon_k+ 
\cdots +\varepsilon_2+1=0$,
we see
%(4.4) %%%
\begin{align}
&(1)\ (\dfrac{\partial s}{\partial y})^\phi=z-1\nonumber\\
&(2)\ (\dfrac{\partial s}{\partial z})^\phi=1-y
\end{align}
Furthermore, we have:
%(4.5) %%%
%\hfil
\begin{align}
&(1)\ (\dfrac{\partial t'}{\partial x_1})^\phi=1-z,  
\
(\dfrac{\partial t'}{\partial x_i})^\phi=0, {\rm for}\ i \ne1,\nonumber\\
&(2)\ (\dfrac{\partial t'}{\partial y})^\phi=x(1-z),\nonumber\\
&(3)\ (\dfrac{\partial t'}{\partial z})^\phi=xy-1.
\end{align}
Now, the Alexander polynomial 
 $\Delta_{\wti{L}}(x, y, z)$ of $\wti{L}$ is obtained as follows.

 Denote by $\hat{M}(\wti{L})$ 
 the $(m+1)\times(m+2)$ matrix obtained from 
 $M(\wti{L})$ by striking out the $m^{\rm th}$ row:
$\biggl(
(\dfrac{\partial r_m}{\partial x_1})^\phi,
 \cdots, (\dfrac{\partial r_m}{\partial x_m})^\phi, 
 (\dfrac{\partial r_m}{\partial y})^\phi,
 (\dfrac{\partial r_m}{\partial z})^\phi\biggr)$

Further,  $\hat{M}(\wti{L})_\nu$ denotes the 
$(m+1)\times(m+1)$ matrix obtained from$\hat{M}(\wti{L})$ 
by striking out the column corresponding to the generator 
$\nu$. (For instance, to get
$\hat{M}(\wti{L})_z$, eliminate the last column of 
$\hat{M}(\wti{L})$.)
 Then the following is known:
%(4.6) %%%
\begin{equation}
\Delta_{\wti{L}}(x, y, z)\doteq
\dfrac{\det  \hat{M}(\wti{L})_z}{1-z}.
\end{equation}

Since the last row of $\hat{M}(\wti{L})_z$ is divisible by 
$1-z$, we have:
%(4.7) %%%
\begin{equation}
\Delta_{\wti{L}}(x, y, z)=\det
\left[
\begin{array}{c|c}
\bigl(\dfrac{\partial r_i}{\partial x_j}\bigr)^\phi&
\delta_{i}(1-x)\\
\noalign{\vskip3pt\hrule \vskip3pt}
(1-y)\bigl(\dfrac{\partial U}{\partial x_j}\bigr)^\phi&
z-1\\
\noalign{\vskip3pt\hrule \vskip3pt}
1\ \ 0\cdots  0&x
\end{array}
\right].
\end{equation}

Let $\hat{L}(r)$ be the link obtained from 
$K_1\cup K_3$ by applying Dehn twists $r(>0)$ times along $K_2$.

Since $\ell k(K_2,K_1)=0$ and $\ell k(K_2, K_3)=1$, 
by Kidwell's theorem \cite[Corollary 3.2]{Ki} we have:
%$(4.8)\hfil %%%
\begin{equation}
\Delta_{\hat{L}}(x, z)=
\dfrac{1}{1-z} \det
\left[
\begin{array}{c|c}
\bigl(\dfrac{\partial r_i}{\partial x_j}\bigr)^\phi_{y=z^r}&
\delta_{i}(1-x)\\
\noalign{\vskip3pt\hrule \vskip3pt}
(1-z^r)\bigl(\dfrac{\partial U}{\partial x_j}\bigr)^\phi&
z-1\\
\noalign{\vskip3pt\hrule \vskip3pt}
1\ \  0\ 0\ \cdots \ 0&x
\end{array}
\right]
\end{equation}.

Further, our knot $K_1(r)$ is obtained from $\hat{L}$ 
by eliminating $K_3$, and hence, by 
Torres' \Theorem \cite{T}, 
noting $\ell k(K_3, K_1(r))=1$, we have:
%(4.9) %%%
\begin{align}
&\Delta_{K_1(r)} (x)=\Delta_{\hat{L}} (x, 1),
\ {\rm and\ hence},\\
%(4.10) %%%
&\Delta_{K_1(r)} (x)=\det
\left[
\begin{array}{c|c}
\bigl(\dfrac{\partial r_i}{\partial x_j}\bigr)^\phi&
\delta_{i}(1-x)\\
\noalign{\vskip3pt\hrule \vskip3pt}
r\bigl(\dfrac{\partial U}{\partial x_j}\bigr)^\phi&
-1\\
\noalign{\vskip3pt\hrule \vskip3pt}
1\ \  0\ 0\ \cdots 0&x
\end{array}
\right]_{y=z=1}
\end{align}

We evaluate $\Delta_{K_1(r)}(x)$ 
by expanding it along the last row, and hence
%(4.11)\\ %%%
%\hspace*{1cm}
%\hfil
\begin{equation}
\Delta_{K_1(r)} (x)\doteq\det
\left[
\begin{array}{c|c}
\bigl(\dfrac{\partial r_i}{\partial x_j}
\bigr)^\phi_{\genfrac{}{}{0pt}{}{i\ge 1}{j\ge 2}}&
\delta_{i}(1-x)\\
\noalign{\vskip3pt\hrule \vskip3pt}
r\bigl(\dfrac{\partial U}{\partial x_j}\bigr)^\phi_{ j\ge2 }&
-1
\end{array}
\right]_{y=z=1}
+(-1)^m x \det 
\left[
\begin{array}{c}
\bigl(\dfrac{\partial r_i}{\partial x_j}
\bigr)^\phi_{\genfrac{}{}{0pt}{}{i\ge 1}{j\ge 1} }\\
\noalign{\vskip3pt}
r\bigl(\dfrac{\partial U}{\partial x_j}\bigr)^\phi_{j\ge 1}
\end{array}
\right]_{y=z=1}
\end{equation}

First we claim:

\begin{lem}\label{lem:6.3}
$\det
\left[
\begin{array}{c}
\bigl(\dfrac{\partial r_i}{\partial x_j}
\bigr)^\phi_{\genfrac{}{}{0pt}{}{i\ge 1}{j\ge  1}}\\
\noalign{\vskip3pt}
r\bigl(\dfrac{\partial U}{\partial x_j}\bigr)^\phi_{j\ge1}\\
\end{array}
\right]_{y=z=1}
=0$ 
\end{lem}

{\it Proof.}
Since $y=1$ and $\es_k+\cdots+\es_2+1=0$, we have
$\sum_{j=1}^m(\frac{\partial r_i}{\partial x_j})_{y=1}^\phi =0$
and $\sum_{j=1}^m\bigl(\frac{\partial U}
{\partial x_j}\bigr)_{y=1}^\phi =0$,
and hence \Lemma \ref{lem:6.3} follows.
\qed

Now we return to the proof of \Proposition \ref{prop:6.2}.
From (4.11) %%%
and \Lemma \ref{lem:6.3}, we see the following:
%(4.12) %%%
\begin{equation}
\Delta_{K_1(r)} (x)\doteq\det
\left[
\begin{array}{c|c}
\bigl(\dfrac{\partial r_i}{\partial x_j}
\bigr)^\phi_{\genfrac{}{}{0pt}{}{i\ge 1}{j\ge 2}}&
\delta_{i}(1-x)\\
\noalign{\vskip3pt\hrule \vskip3pt}
r\bigl(\dfrac{\partial U}{\partial x_j}\bigr)^\phi_{j\ge2}&
-1
\end{array}
\right]_{y=1}
\end{equation}

The determinant is decomposed into two
terms as follows:
%
%(4.13) %%%
%%\hfil
\begin{equation}
\Delta_{K_1(r)}(x) \doteq \det
\left[
\begin{array}{c|c}
\bigl(\dfrac{\partial r_i}{\partial x_j}
\bigr)^\phi_{\genfrac{}{}{0pt}{}{i\ge 1}{j\ge 2}}&
\delta_{i}(1-x)\\
\noalign{\vskip3pt\hrule \vskip3pt}
r\bigl(\dfrac{\partial U}{\partial x_j}\bigr)^\phi_{j\ge2}&
0
\end{array}
\right]_{y=1}
+
\det
\left[
\begin{array}{c|c}
\bigl(\dfrac{\partial r_i}{\partial x_j}
\bigr)^\phi_{\genfrac{}{}{0pt}{}{i\ge 1}{j\ge 2}}&
0\\
\noalign{\vskip3pt\hrule \vskip3pt}
r\bigl(\dfrac{\partial U}{\partial x_j}\bigr)^\phi_{j\ge2}&
-1
\end{array}
\right]_{y=1}
\end{equation}

The second term is equivalent to

$\det\Bigl[\bigl(\dfrac{\partial r_i}{\partial x_j}\bigr)^\phi_
{\genfrac{}{}{0pt}{}{1\le i \le m-1}{2\le j\le m}}\Bigr]_{y=1}
$
that is equal to $\Delta_{K_1}(x)$ (up to $\pm x^k$).
Therefore, the final step is to show that

%(4.14) %%%
%\hfil
\begin{equation}
\det
\left[
\begin{array}{c|c}
\bigl(\dfrac{\partial r_i}{\partial x_j}
\bigr)^\phi_{\genfrac{}{}{0pt}{}{i\ge 1}{j\ge 2}}&
\delta_{i}\\
\noalign{\vskip3pt\hrule \vskip3pt}
\bigl(\dfrac{\partial U}{\partial x_j}\bigr)^\phi_{j\ge2}&
0
\end{array}
\right]_{y=1}
\doteq
\left[\dfrac{\Delta_{B(2\alpha, \beta)}(x,y)}{1-y}\right]_{y=1}.
\end{equation}

To show (4.14) %%%
we go back to $M(\wti{L})$ and compute 
$\Delta_{\wti{L}} (x, y, z)$ in a different way. 
We use the following formula:
%(4.15)%%%
\begin{equation}
\Delta_{\wti{L}} (x, y, z)=
\dfrac{\det\hat{M}(\wti{L})_y}{1-y}
\end{equation}

Then the row 
$(\dfrac{\partial s}{\partial x_1},
\dfrac{\partial s}{\partial x_2}, 
\cdots, \dfrac{\partial s}{\partial x_m},
\dfrac{\partial s}{\partial z})^\phi$ 
is divisible by $1-y$, and hence, we have:
%(4.16) %%%
\begin{equation}
\Delta_{\wti{L}} (x, y, z)=\det
\left[
\begin{array}{c|c}
\bigl(\dfrac{\partial r_i}{\partial x_j}\bigr)^\phi_{j\ge1}&
\begin{array}{c}
1-x\\
0\\
\vdots\\
0\\
\end{array}\\
\noalign{\vskip3pt\hrule \vskip3pt}
\bigl(\dfrac{\partial U}{\partial x_j}\bigr)^\phi_{j\ge1}&
1\\
\noalign{\vskip3pt\hrule \vskip3pt}
1\hspace{-1.5mm}-\hspace{-1.5mm}z\ 0 \cdots 0&
xy-1
\end{array}
\right]
\end{equation}

Now we try to find $\Delta_L (x, y)$ from 
$\Delta_{\wti{L}} (x,y,z)$.

To do this, we eliminate $K_3$ from $\wti{L}=K_1 \cup K_2 \cup K_3$.
Then, since $\ell k (K_3, K_1)=\ell k (K_3, K_2)=1$, 
Torres' \Theorem \cite{T} implies:
%(4.17)\hfil %%%
\begin{equation}
\Delta_L (x, y)=\dfrac{\Delta_{\wti{L}}(x, y, 1)}{xy-1},
\end{equation}

that is, from (4.16),%%%
%(4.18)\hfil %%%
\begin{equation}
\Delta_L (x, y)=\det
\left[
\begin{array}{c}
\bigl(\dfrac{\partial r_i}{\partial x_j}\bigr)^\phi_{j\ge1}\\
\noalign{\vskip3pt}
\bigl(\dfrac{\partial U}{\partial x_j}\bigr)^\phi_{j\ge1}
\end{array}
\right]_{z=1}
=N
\end{equation}

We describe $N$ precisely.
First we note:
%(4.19) 
\begin{align}
&(1)\  %%%
{\displaystyle
\sum_{j=1}^m (\frac{\partial r_i}{\partial x_j})^\phi
=
\begin{cases}
y^{\es}-1, {\rm if\ } r_i {\rm \ is\  of\  the\  from}: 
y^{\es} x_i 
y^{-\es} x_{i+1}^{-1}, \es=\pm1
\\
0, {\rm otherwise}.
\end{cases}
}\nonumber\\
&(2)\  
{\displaystyle
\sum_{j=1}^m (\frac{\partial U}{\partial x_j})^\phi=0
}.
\end{align}

Therefore, if we add all columns of $N$ to the first column 
to get $N_1$, then the first column of $N_1$ is 
divisible by $1-y$.  Further,

%(4.20) %%%
%\hfil
\begin{equation}
{\displaystyle
\sum_{j=1}^m (\frac{\partial r_i}{\partial x_j}
)^\phi=\es y^{\frac{\es -1}{2}} (1-y)}.
\end{equation}
Since $\es y^{\frac{\es -1}{2}} =\delta_i$,
we have:
%(4.21) \hfil %%%
\begin{equation}
\dfrac{N_1}{1-y}=(-1)^m \det     
\left[
\begin{array}{c|c}
\bigl(\dfrac{\partial r_i}{\partial x_j}\bigr)^\phi_{j\ge2}&
\delta_{i}\\
\noalign{\vskip3pt\hrule \vskip3pt}
\bigl(\dfrac{\partial U}{\partial x_j}\bigr)^\phi_{j\ge2}&
0
\end{array}
\right]
=\dfrac{\Delta_L (x, y)}{1-y}.
\end{equation}

Evaluations of both polynomials at $y=1$ give (4.14). %%%

The proof of \Proposition \ref{prop:6.2} is now completed.
\qed

\section{Alexander polynomials (II)}%7 %%Feb25 5

We have established some relationships between 
the Alexander polynomial of 
$K(2\alpha, \beta|r)$ 
and that of the 2-bridge link $B(2\alpha, \beta)$. 
However, these relations are not sufficient to our purpose.
Therefore, in this section, we prove some subtle properties of 
$\Delta_{B(2\alpha, \beta)}(x, y)$. 
These properties are indispensable to study 
the Alexander polynomial of our knot 
$K(2\alpha, \beta|r)$. See \Theorem \ref{thm:7.4}.

Let $S=\{P_1,d_1, Q_1, e_1, P_2, \cdots, P_m, d_m, Q_m\}$
be the canonical decomposition of the continued fraction of 
$\beta/2\alpha$.
Let $\rho_i, \rho_j, \rho, \lambda_i, \lambda_j$ and 
$\lambda$ be integers as defined in Definition \ref{dfn:3.8}.
Now, by \Proposition \ref{prop:3.9}, we can write
%(5.1) 
\begin{equation}
\Delta_{B (2\alpha, \beta)}(x, y)=
f_{\lambda-1}(x)
y^{\lambda-1}+f_{\lambda-2}(x) y^{\lambda-2}+\cdots+f_0(x),
\end{equation}
 where $f_i(x), 0 \le i \le \lambda-1$, 
are integer polynomials in $x$ of degree at most $\lambda-1$.
 
 Our purpose is to determine these polynomials $f_i(x)$, 
in particular, $f_{\lambda-1}(x)$.
 
\medskip
{\bf 5.1 Skein relation}
 
Let $[[u_1, v_1, u_2, v_2, \cdots, $
$u_s, v_s, u_{s+1}]]$ 
be the continued fraction of $\beta/2\alpha$. 
Then it is shown in \cite[Theorem 2 (4.2)]{K} that
%(5.2) \hfil
\begin{align}
&\ \ \ \Delta_{B (2\alpha, \beta)}(x,y)\nonumber\\
&=v_s (x-1)(y-1) F_{u_{s+1}} (x, y)
\Delta[[u_1, v_1, \cdots, u_s]]
-\Delta[[u_1, v_1, 
\cdots, v_{s-1}, u_s +u_{s+1}]],
\end{align}
where 
$\Delta [[c_1, \cdots, c_k]]$ is the Alexander polynomial of 
the 2-bridge link associated to the continued fraction 
$[[c_1, c_2, \cdots, c_k]]$, and $F_{n}(x,y)$ is defined below:
%(5.3) 
\begin{align}
&(1)\ F_0 (x,y)=0.\nonumber\\
&(2)\ {\rm For}\ n>0,\nonumber\\
&\ \ \ (a)\   F_n (x,y)=1+xy+\cdots+(xy)^{n-1}=
\frac{(xy)^n-1}{xy-1},\nonumber\\
&\ \ \ (b)\  F_{-n} (x,y)=-\{(xy)^{-1}+\cdots+(xy)^{-n}\}
=\frac{-1}{(xy)^n} F_n(x,y).
\end{align}

%\begin{tabular}{lll}
%(5.3) & (1) & $F_0 (x,y)=0$.\\
%          & (2)  & For $n>0$, \\
%& & (a) $F_n (x,y)=1+xy+\cdots+(xy)^{n-1}=
%\frac{(xy)^n-1}{xy-1}$,\\
%& & (b) $F_{-n} (x,y)=-\{(xy)^{-1}+\cdots+(xy)^{-n}\}
%=\frac{-1}{(xy)^n} F_n(x,y)$.
%\end{tabular}

Note that $F_c(x,y)=\Delta[[c]]$.

Formula (5.2) is obtained by applying crossing changes 
and smoothing
at $v_s$,i.e.,
at the crossings
corresponding to $v_s$.
 
We should note that (5.2) is slightly different from the 
original formula given in \cite[(4.2)]{K}, 
since we use a different notation.
 
By applying (5.2) on all $v_j, j=1, 2, \cdots, s$, 
we obtain $\Delta_{B(2\alpha, \beta)}(x, y)$ 
in terms of various $
\Delta [[c]]$, where $c$ is written as the sum of $u_i$.
 
The following example illustrates a calculation.

\begin{ex}\label{ex:7.1}
Write $\frac{\beta}{2\alpha}=[[u_1, v_1, u_2, v_2, u_3]]$.
Then,
\begin{align*}
&\ \ \ \ \Delta_{B(2\alpha, \beta)}\\
&=v_2 (x-1)(y-1) F_{u_3}(x, y) 
\Delta [[u_1, v_1, u_2]]-\Delta[[u_1,v_1, u_2+u_3]]\\
&=v_2 (x-1)(y-1) F_{u_3}(x, y)\{v_1(x-1)(y-1) 
F_{u_2}(x, y) F_{u_1}(x, y)-\Delta[[u_1+u_2]]\}\\
&\ \ \ \ -\{v_1(x-1)(y-1) F_{u_2+u_3}(x, y)
 \Delta[[u_1]] 
 -\Delta[[u_1+u_2+u_3]]\}\\
&=v_1 v_2 (x-1)^2 (y-1)^2 F_{u_1} F_{u_2} F_{u_3} 
-(x-1)(y-1)\{v_1 F_{u_1} F_{u_2+u_3} +
v_2 F_{u_1+u_2} F_{u_3}\}\\
&\ \ \ \ +F_{u_1+u_2+u_3}
\end{align*}
\end{ex}

As is illustrated in Example \ref{ex:7.1}, we see that
$\Delta_{B(2\alpha, \beta)}(x, y)$ is of the following form:
%(5.4) 
\begin{equation}
\Delta_{B(2\alpha,\beta)}(x,y)=
\sum_{\genfrac{}{}{0pt}{}{0\le k \le s}{
1\le i_1<i_2<\cdots<i_k\le s}}
(-1)^k v_{i_1}v_{i_2}\cdots v_{i_k}(x-1)^k (y-1)^k
F_{\mu_1}F_{\mu_2}\cdots F_{\mu_k},
\end{equation}
where the summation is taken over all indices 
$i_j$ such that 
$1\le i_1<i_2<\cdots<i_k\le s$, and 
$\mu_j$ is of the form: 
$\mu_j=u_{j_1}+u_{j_1+1} +\cdots+u_{j_1+p}$ and  
$\mu_1 + \mu_2 + \cdots +\mu_{k+1} 
= u_1 +u_2 + \cdots+ u_{s+1}$.

For convenience, we denote by $\Lambda_{p,r}$ 
the set of all $p$ indices $i_1, \cdots, i_p$
such that $1 \le i_1< \cdots <i_p\le r$.
Since $F_c(x,y)$ is a rational function, we replace 
$F_c(x,y)$ by a polynomial $\widetilde{F}_c(x,y)$ below.
%%\begin{tabular}{ll}
%(5.5)
\begin{align} 
&{\it For}\ n>0,\nonumber\\
&(1)\ \wti{F}_n(x,y) = (xy-1) F_n(x,y) =(xy)^n -1.\nonumber\\
&(2)\ \wti{F}_{-n}(x,y) = (xy)^n(xy-1) 
F_{-n}(x,y)= (-1)\wti{F}_n(x,y)=(-1)\{(xy)^n -1\}.
\end{align}
Using these polynomials, we obtain an integer polynomial   
$\delti_{B(2 \alpha,\beta)}(x,y)$ from
$\Delta_{B(2 \alpha,\beta)}(x,y)$:
%
%(5.6)\hfil
\begin{equation}
\delti_{B(2 \alpha,\beta)}(x,y) =
(xy)^{\sum_{j=1}^m \lambda'_j}
(xy-1)^{\sum_i(s_i+1)+\sum_j (q_j+1)}
\Delta_{B(2 \alpha,\beta)}(x,y).
\end{equation}
Therefore we have:
%
%(5.7)
\begin{align}  
(1)\ &\max \mbox{$y$-}\deg\delti_{B(2 \alpha,\beta)}(x,y)\nonumber\\
&=\max \mbox{$y$-}\deg \Delta_{B(2 \alpha,\beta)}(x,y)
+{\displaystyle\sum_{i=1}^{m}s_i +\sum_{j=1}^{m}q_j +2m
}\nonumber\\
&=  \lambda +\sum_{i=1}^{m}s_i +\sum_{j=1}^{m}q_j + 2m-1,\nonumber\\
(2)\ &  \min \mbox{$y$-}deg \delti_{B(2 \alpha,\beta)}(x,y) = 0.
\end{align}
Now we can write
%(5.8) \hfil
\begin{align}
\delti_{B(2 \alpha,\beta)}(x,y)
&= \tif_\nu (x) y^\nu  + \cdots +\tif_0(x),\ {\rm  and}\nonumber\\
\tif_\nu (x)&=f_{\lambda -1}(x) x^{\sum_1^m s_i + 
\sum_1^m q_j +2m},
\end{align}
where $\nu= \lambda +\sum_{i=1}^{m}s_i + 
\sum_{j=1}^{m}q_j + 2m - 1$.
 
First we show
%(5.9)
\begin{equation}
\deg \tif_\nu(x) = \lambda +
\sum_1^m s_i +   \sum_1^m q_j -
\rho
+2m - 1.
\end{equation}

 \medskip
{\bf 5.2 Proof of (5.9) (I)}

We consider two special cases.

Case 1.  All $u_i >0$.

Consider  $\beta/2\alpha = [[u_1, v_1,u_2,v_2, 
\cdots, u_s,v_s,u_{s+1}]]$.\\
Then
$\displaystyle
\Delta_{B(2 \alpha,\beta)}(x,y)
=\sum_{\Lambda_{k,s},0\le k\le s}(-1)^k v_{i_1} v_{i_2} 
\cdots v_{i_k} (x-1)^k(y-1)^k 
F_{\mu_1} F_{\mu_2} \cdots F_{\mu_{k+1}}$,\\
where  $\mu_i>0, 1\le i\le k+1,{\rm and}\
 \lambda =\sum_{i=1}^{k+1} \mu_i,  \lambda' =0$.
Therefore,
%(5.10)\hfil
\begin{align}
&\ \ \ \ \delti_{B(2 \alpha,\beta)}(x,y)\nonumber\\
&=
(xy-1)^{s+1} \Delta_{B(2 \alpha,\beta)}(x,y)\nonumber\\
&=\sum_{\genfrac{}{}{0pt}{}
 {\Lambda_{k,s}}{ 0\le k\le s}}
 (-1)^k v_{i_1} v_{i_2} 
         \cdots v_{i_k} (x-1)^k
 (y-1)^k \wti{F}_{\mu_1} \wti{F}_{\mu_2} \cdots 
 \wti{F}_{\mu_{k+1}}(xy-1)^{s-k}.
 \end{align}

Case 2.   All  $-u_i <0$.

Consider   $\beta/2\alpha= [[-u_1, -v_1,-u_2,-v_2, \cdots, 
-u_q,-v_q,-u_{q+1}]]$. Then
 $\Delta_{B(2 \alpha,\beta)}(x,y) =    (-1)^k(-v_{i_1})(-v_{i_2}) 
 \cdots (-v_{i_k}) 
 (x-1)^k(y-1)^k F_{-\mu_1} F_{-\mu_2} \cdots 
 F_{-\mu_{k+1}}$,\\
and hence  $\lambda'  = u_1 + u_2 +\cdots+u_{q+1}= 
\mu_1 +\mu_2+\cdots+\mu_{k+1}$. 
Therefore:
%(5.11) 
\begin{align}
&\ \ \ \ 
\delti_{B(2\alpha,\beta)}(x,y) \nonumber\\
&=(xy)^{\lambda'}
(xy-1)^{q+1}
\Delta_{B(2 \alpha,\beta)}(x,y)\nonumber\\
&= \sum_{\Lambda_{k,q}, 0\le k \le q}
 v_{i_1}
 v_{i_2} \cdots 
v_{i_k} (-1)^{k+1}(x-1)^k(y-1)^k 
\wti{F}_{\mu_1} \wti{F}_{\mu_2} 
\cdots \wti{F}_{\mu_{k+1}}(xy-1)^{q-k}\nonumber\\
&=-\sum_{\Lambda_{k,q},0\le k \le q}
 (-1)^k v_{i_1} v_{i_2} 
\cdots v_{i_k} (x-1)^k(y-1)^k 
\wti{F}_{\mu_1} \wti{F}_{\mu_2} \cdots \wti{F}_{\mu_{k+1}}
(xy-1)^{q-k}.
\end{align}
 
Note that (5.10) and (5.11) are of the same form.
 
Now consider the general case.
Let $\{P_1, d_1, Q_1, e_1, P_2,\cdots, P_m, d_m, Q_m\}$ 
be the canonical decomposition of the continued fraction of 
$\beta/2\alpha$.\\
Denote 
 $P_i=[[a_{i,1},b_{i,1},a_{i,2},b_{i,2}, \cdots, a_{i,s_i}, 
 b_{i,s_i},a_{i,s_i +1}]]$, $1\le i \le m$, and\\
 $Q_j= [[-a'_{j,1},-b'_{j,1},-a'_{j,2},-b'_{j,2}, 
 \cdots,-a'_{j,q_j},-b'_{j,q_j},-a'_{j,q_j+1}]]$, $1\le j\le m$,
 where $a_{i,p}>0$ and $a'_{j,p} >0$, but 
 $b_{i,q}, b'_{j,q}$ are arbitrary.
Then by \Proposition \ref{prop:3.9},
 
max $y$-deg $\Delta_{B(2 \alpha,\beta)}(x,y) = 
{\displaystyle
\sum_{i=1}^{m}\sum_{k=1}^{s_{i}+1}
a_{i,k} + \sum_{j=1}^{m}\sum_{k=1}^{q_{j}+1} a'_{j,k} - 1 
}
= \lambda - 1$.
 
First we try to find the term with the max $y$-degree 
in $\delti_{B(2 \alpha,\beta)}(x,y)$.
 
Denote by $\Delta_{P_i}(x,y)$  (resp. $\Delta_{Q_j}(x,y))$ 
the Alexander polynomial of the 2-bridge link associated 
to $P_i$  (resp. $Q_j$). 
Then, as we did above, we obtain\\
$\Delta_{P_i} (x,y) = 
{\displaystyle
\sum_{\Lambda_{k,s_i},0\le k \le s_i}
}(-1)^k b_{i,p_1} b_{i,p_2} 
\cdots b_{i,p_k} (x-1)^k(y-1)^k F_{\mu_1} 
F_{\mu_2} 
\cdots F_{\mu_{k+1}}$,\\
where $\mu_1+\mu_2 + \cdots +\mu_{k+1}= \lambda_i$,
and hence,
%(5.12)\hfil
\begin{align}
&\ \ \ \  \delti_{P_i} (x,y) \nonumber\\
&= 
(xy-1)^{s_i+1}
\Delta_{P_i} (x,y)\nonumber\\
&=\sum_{\Lambda_{k,s_i}}
(-1)^k b_{i,p_1} b_{i,p_2} \cdots b_{i,p_k} 
(x-1)^k(y-1)^k \wti{F}_{\mu_1} \wti{F}_{\mu_2} 
\cdots \wti{F}_{\mu_{k+1}}(xy-1)^{s_i-k},
\end{align}
where $\wti{F}_\mu=(xy)^{\mu}-1, \mu>~0$.
 
On the other hand,
  
$\Delta_{Q_j} (x,y) = 
{\displaystyle
\sum_{\lambda_{k,q_j}}
}
(-1)^k b'_{j,r_1} b'_{j,r_2} 
\cdots b'_{j,r_k} (x-1)^k(y-1)^k F_{-\mu'_1} F_{-\mu'_2} 
\cdots F_{-\mu'_{k+1}}$, and hence we have:
 % (5.13)\hfil
\begin{align}
&\ \ \ \ \delti_{Q_j} (x,y)\nonumber\\
& = (xy)^{\lambda'}
(xy-1)^{q_j+1} \Delta_{Q_j}(x,y) \nonumber\\
&=-\sum_{\Lambda_{k,q_j}}
(-1)^k b'_{j,r_1} b'_{j,r_2} \cdots 
b'_{j,r_k} (x-1)^k(y-1)^k \wti{F}_{\mu'_1} \wti{F}_{\mu'_2} 
\cdots \wti{F}_{\mu'_{k+1}}(xy-1)^{q_j-k}.
\end{align}
 
\medskip
{\bf 5.3. Proof of (5.9) (II)}
 
To evaluate  $\Delta_B(x,y)$, we must split and smooth
at various crossings.
  We classify these operations into two types.
  
Type 1.  Split all crossings at every $d_i$ and $e_j$.

Type 2   Smooth some crossings at some $d_i$ and/or 
$e_j$.
 
From Type 1 operation, we obtain the following term in
$\delti_B(x,y)$ :
%(5.14)  
\begin{align}
A = &(-1)^m
d_1\cdots d_m (-1)^{m-1} 
e_1 \cdots e_{m-1}(x-1)^{2m-1} (y-1)^{2m-1}\nonumber\\
&\times
\prod_{i=1}^{m}\delti_{P_i} (x,y)
\prod_{j=1}^{m}\delti_{Q_j} (x,y).
\end{align}
 
Terms in $A$ with the max $y$-degree are obtained by\\
(1) taking $y^{2m-1}$ from $(y-1)^{2m-1}$, \\
(2) taking , in each $P_i$, $y^k$ from $(y-1)^k$, 
$(xy)^{\mu_i}$ from each $\wti{F}_{\mu_i}$ and 
$(xy)^{s_i -k}$ from 
$(xy-1)^{s_{i}-k}$, and\\
(3) taking, in each $Q_j, y^k$ from $(y-1)^k, 
(xy)^{\mu'_i}$ 
from $\wti{F}_{\mu'_i}$, and
$(xy)^{q_j -k}$ from 
$(xy-1)^{q_j -k}$.

Therefore, the max $y$-degree in $A$ is 

${\displaystyle
   2m-1+\sum_{i=1}^{m}(k+\mu_1+\cdots+\mu_{k+1}+s_i -k)+    
   \sum_{j=1}^{m}(k+\mu'_1 + \cdots \mu'_{k+1}+q_j- k)
   }$
   
${\displaystyle
= 2m-1+\sum_{i=1}^{m}(s_ i +\lambda_i ) +
   \sum_{j=1}^{m}(q_j + \lambda'_j)
   = 2m-1 +\sum_{i=1}^{m}s_i  + \sum_{j=1}^{m}q_j + \lambda
   }$.\\
While, the min $y$-deg in A is obviously $0$.  \\
Since $
{\displaystyle
\delti_B(x,y) = (xy)^{\sum_{i=1}^{m}\lambda'_i}
(xy-1)^{\sum_{i=1}^{m}(s_i+1)+\sum_{j=1}^{m}(q_j+1)}
\Delta_B(x,y)
}$, \\
the $y$-degree of
$\Delta_B(x,y)$ is at least
$2m-1 +\sum s_i+ \sum q_j + \lambda - 
(\sum s_i  +m +\sum q_j  + m) = \lambda - 1$,
that coincides with \Proposition \ref{prop:3.9}.  
Therefore, 
these terms are in fact the terms with maximal $y$-degree.

\medskip
{\bf 5.4. Proof of (5.9) (III)}
 
Next we show that Type 2
operation does not yield a term with 
max $y$-degree in $\delti_B(x,y)$.
To see this, we can assume without loss of generality that 
we smooth only crossings at $d_1$, but not others.  
Namely, we split at other crossings $d_i (i\neq1)$ and 
$e_j, 1\le j\le m-1$.
 
 \medskip
 
Case 1. Suppose $a_{1,s_1+1} > a'_{1,1}$.

By smoothing at $d_1$, we have a new canonical 
decomposition of the new continued fraction:
$\hat{S} = \{\hat{P}_1, \hat{c}_1, \hat{Q}_1, e_1, P_2, 
d_2, Q_2, e_2, \cdots, 
P_m, d_m, Q_m\}$ ,
where

 $\hat{P}_1=[[a_{1,1},b_{1,1},a_{1,2},b_{1,2}, \cdots, 
 a_{1,s_1},
 b_{1,s_1},a_{1,s_1 +1} - a'_{1,1}]], \ 
\hat{c}_1 = - b'_{1,1}$ and

 $\hat{Q}= [[-a'_{1,2},-b'_{1,2}, 
 \cdots,-a'_{1,q_1},-b'_{1,q_1},-a'_{1,q_1+1}]]$.
 
Consider
$\delti_B(x,y) = (xy)^{\lambda'}
(xy-1)^{\sum (s_i+1)+\sum(q_j+1)}
\Delta_B(x,y)$.
Using the previous argument, 
we can determine the terms of max 
$y$-degree of $\Delta(\hat{S})$ in $\delti_B(x,y)$.
Since the terms 
of max $y$-degree are obtained as those in each 
$P_i$ and $Q_j$, we will determine these terms for 
$\hat{P}_1$ and $\hat{Q}_1$.
For $\hat{P}_1$, the max $y$-degree is
%(5.15)
%\hfil
\begin{equation}
k + \hat{\mu}_1 +  \cdots  + \hat{\mu}_{k+1} + 
s_1 + 1 - (k + 1).
\end{equation}

Since $\hat{\mu}_1 +  \cdots  + \hat{\mu}_{k+1} = 
a_{1,1} + a_{1,2} + 
\cdots + a_{1,s_1 +1} - a'_{1,1} = \lambda_1 - a'_{1,1}$, 
it follows from (5.15) that the max $y$-degree is  
$\lambda_1 - a'_{1,1} + s_1$.
For $\hat{Q}_1$, the maximal terms are contained in\\
${\displaystyle
%\sum_{
%\genfrac{}{}{0pt}{}
%{2\le r_1<\cdots<r_k\leq q_1}
%{k=0,1, \ldots, q_{1}-1}
%}
\sum
(-1)^k b'_{1,r_1}\cdots b'_{1,r_k} (x-1)^k(y-1)^k 
\wti{F}_{\mu'_1}  \cdots 
\wti{F}_{\mu'_{k+1}} (xy-1)^{q_1-k-1} 
        (xy)^{a'_{1,1}} (xy-1),
        }$
where the sum is taken over
$2\le r_1<\cdots<r_k\leq q_1,
k=0,1, \ldots, q_{1}-1$.
        
Since the original multipliers 
$(xy)^{\lambda'}$
cannot be cancelled out in this case, $(xy)^{a'_{1,1}}$
remains.
Therefore, max $y$-degree in $\hat{Q}_1$ is
   $k + q_1 - k - 1 + \hat{\mu}'_1 +  \cdots  + 
   \hat{\mu}'_{k+1} + a'_{1,1} + 1 
  = q_1 + \lambda'_1 - a'_{1,1} + a'_{1,1} 
  = q_1 + \lambda'_1$,
since  $\hat{\mu'}_1 +  \cdots  + \hat{\mu'}_{k+1} 
= \lambda'_1 - a'_{1,1}$,  and hence, the 
max $y$-deg of $\delti_B(x,y)$ is
   $2m - 1 + (\lambda_1 - a'_{1,1}) 
   +s_1 + \lambda'_1 + q_1 +\sum_{i=2}^{m}(s_i +\lambda_i ) 
   +\sum_{j=2}^m (q_j + \lambda_j)
= 2m - 1 +\sum_{i=1}^{m}s_i +\sum_{j=1}^mq_j + \lambda - a'_{1,1}$.
Since $a'_{1,1}  > 0$, we cannot get a term of 
the max $y$-degree from $\Delta(\hat{S})$.
 
 \medskip
Case 2.  $a_{1,s_1 +1} < a'_{1,1}$ or  $a_{1,s_1 +1} = a'_{1,1}$. 
A similar argument works, and hence omit the details.
Therefore, to evaluate $\wti{f}_\nu (x)$, it suffices to consider
$\delti(P_i)$, since the treatment for $\delti(Q_j)$
is similar to $\delti(P_i)$.  
In other words, we will show the following:
 
\begin{prop}\label{prop:7.1}%7.2 June24
Let $S=[[u_1, v_1,u_2,v_2, \cdots, u_s,v_s,u_{s+1}]]$, where 
$u_i > 0,  1\le i\le s +~1$.
Write 
   $\delti_{B}(x,y) = f_{s+\lambda}(x) y^{s+\lambda} + 
   \cdots + f_0(x)$,
where $\lambda = \sum_{i=1}^{s+1}u_i$.  Then we can write
as follows, using some integer $\gamma_{s+\lambda-\rho}
\neq 0$.
%(5.16) 
\begin{equation}
f_{s+\lambda}(x) = 
\gamma_{s+\lambda-\rho}x^{s+\lambda-\rho}
 +  \cdots + \gamma_\zeta x^\zeta, \ 
 {\it for\ some}\ 
 \zeta \ge 0, s+\lambda-\rho>\zeta.
 \end{equation}

 \end{prop}

\medskip
{\bf 5.5. Auxiliary Lemmas}
 
Before we proceed to the proof of 
\Proposition \ref{prop:7.1}, we show the following 
two lemmas.

\begin{lem}\label{lem:7.2}  %7.3 June 24
Assume $n\ge k\ge 0$ and $n\ge m \ge 0$.  Then
%(5.17)
\begin{align}
&\ \ \ \bi{n}{k}-\bi{n-1}{k-1}\bi{m}{m-1}+\bi{n-2}{k-2}\bi{m}{m-2}- \nonumber\\
&\ \ \  \ \ \ 
\cdots + (-1)^\ell \bi{n-\ell}{k-\ell}\bi{m}{m-\ell}+
\cdots+(-1)^m
\bi{n-m}{k-m}\bi{m}{0}\nonumber\\
&=\bi{n-m}{k}
\end{align}
\end{lem}

Note.  In (5.17) we assume that $\bi{n}{k}=0$
if $n\le 0$ or $k\le 0$, and  $\bi{0}{0}=1$.
 
 \medskip
 
{\it Proof.}  
We prove (5.17) by induction on $n,k$ and $m$.
Direct calculations prove the validity of the first step.
Suppose (5.17) holds up to $n, k$ and $m-1$.
Then we see the following:
The LHS of (5.17) is
\begin{align*}
&\bi{n}{k}-\Bigl\{\bi{n-1}{k-1}\bi{m-1}{1}
+\bi{n-1}{k-1}\bi{m-1}{0}\Bigr\}\\
&\ \ \ 
+\Bigl\{\bi{n-2}{k-2}\bi{m-1}{2}
+\bi{n-2}{k-2}\bi{m-1}{1}
\Bigr\}-\cdots\\
&\ \ \ 
+(-1)^{m-1}
\Bigl\{\bi{n-m+1}{k-m+1}\bi{m-1}{m-1}+\bi{n-m+1}{k-m+1}
\bi{m-1}{m-2}\Bigr\}\\
&\ \ \ 
+(-1)^m\Bigl\{
\bi{n-m}{k-m}\bi{m-1}{m-1}\Bigr\}\\
&=\bi{n}{k}-\bi{n-1}{k-1}\bi{m-1}{1}+\bi{n-2}{k-2}\bi{m-1}{2}
+\cdots\\
&\ \ \  + (-1)^{m-1}\bi{n-m+1}{k-m+1}\bi{m-1}{m-1}\\
&\ \ \ \ 
-\Bigl\{
\bi{n-1}{k-1}\bi{m-1}{0}-\bi{n-2}{k-2}\bi{m-1}{1}+\cdots \\
&\ \ \ +
(-1)^{m-1}\bi{n-m}{k-m}\bi{m-1}{m-1}\Bigr\}\\
&=
\bi{n-(m-1)}{k}-\bi{n-m}{k-1}=\bi{n-m}{k}+\bi{n-m}{k-1}-
\bi{n-m}{k-1}=\bi{n-m}{k},
\end{align*} 

by induction hypothesis.
\qed

\begin{lem}\label{lem:7.3}
Let $n\ge k \ge0$.  
Then the following equality holds among integer polynomials 
in $n$ variables $x_1,x_2, \cdots,x_n$:
\begin{align*}
&\bi{n}{k}x_1\cdots x_n -\bi{n-1}{k}\sum_{\Lambda_{n-1}}x_{i_1}
\cdots x_{i_{n-1}}\\
&\ \ \ +\bi{n-2}{k}\sum_{\Lambda_{n-2}} x_{i_1}\cdots x_{i_{n-2}}+
\cdots 
+(-1)^{n-k} \bi{k}{k} \sum_{\Lambda_{k}} x_{i_1}\cdots
x_{i_k}\\
& 
=\bi{n}{k}(x_{1}-1)\cdots(x_{n}-1)+\bi{n-1}{k-1}
\sum_{\Lambda_{n-1}}(x_{i_1}-1)\cdots(x_{i_{n-1}}-1)\\
&\ \ \ 
+\bi{n-2}{k-2}\sum_{\Lambda_{n-2}}(x_{i_1}-1)\cdots
(x_{i_{n-2}}-1)+\cdots\\
&\ \ \ +\bi{n-k}{0}\sum_{\Lambda_{n-k}}
(x_{i_1}-1)\cdots(x_{i_{n-k}}-1),
\end{align*}
where the summation is taken over the set $\Lambda_j$ 
consisting of all indices $i_1,\cdots, i_j$ such that 
$1 \le i_1<\cdots< i_j\le n$.  
If $n-k=0$, then the last term on the right side is interpreted 
as~$1$.
\end{lem}

{\it Proof.} Since polynomials on both sides are 
symmetric polynomials over the symmetric group $S_n$, 
it is enough to compare the  coefficients of 
$x_1 x_2 \cdots x_m$, $1\le m\le n$.
For example, the constant term of the LHS is $0$ if $k>0$, 
while that of the RHS is

$(-1)^n \bi{n}{k} +(-1)^{n-1}\bi{n-1}{k-1}\bi{n}{n-1}
+\cdots+(-1)^{n-k}\bi{n-k}{0}\bi{n}{n-k}$\\
$=(-1)^n\sum_{i=0}^k (-1)^i \bi{n-i}{k-i}\bi{n}{n-i}$\\
$=(-1)^n\sum_{i=0}^k (-1)^i\bi{n}{k}\bi{k}{i}$\\
$=(-1)^n\bi{n}{k}\sum_{i=0}^k (-1)^i\bi{k}{i}$\\
$=0
$

First,  $x_1 x_2 \cdots x_n$ appears $\bi{n}{k}$ 
times in both sides.  Thus the 
formula is true for $x_1 x_2 \cdots x_n$.
Next, consider $x_1 x_2 \cdots x_{n-1}$. 
This appears  $-\bi{n-1}{k}$ times in the LHS, while 
it appears, in the RHS,
$-\bi{n}{k}+\bi{n-1}{k-1}=-\bi{n-1}{k}$ times.
Thus the formula is true.
In general,  $x_1 x_2 \cdots x_{n-r}$, $r\ge1$,  
appears $(-1)^r \bi{n-r}{k}$ times in the LHS, 
while in the RHS, it appears as many times as
\begin{align*}
&(-1)^r \Bigl\{
 \bi{n}{k}-\bi{n-1}{k-1}\bi{r}{r-1}+\bi{n-2}{k-2}
\bi{r}{r-2}- \cdots + (-1)^r \bi{n-r}{k-r}\Bigr\}\\
&=(-1)^r \bi{n-r}{k}\ {\rm (by}\ {\rm \Lemma}\ \ref{lem:7.2}).
\end{align*}
\qed

{\bf 5.6. Proof of \Proposition \ref{prop:7.1}.}
 
By (5.10), we can write
%(5.18)
%\hfil
\begin{align}
\delti_B(x,y) &= (xy-1)^{s+1}   \Delta_B(x,y)\nonumber\\
&=\sum_{p=0}^{s}\sum_{\Lambda_{p,s}}
         (-1)^p v_{i_1} v_{i_2} 
\cdots v_{i_p} (x-1)^p(y-1)^p \wti{F}_{\mu_1} \wti{F}_{\mu_2} 
        \cdots \wti{F}_{\mu_{p+1}}
        (xy-1)^{s-p},
\end{align}
where $\mu_1 + \cdots + \mu_{p+1} =
  \lambda = u_1 + u_2+ \cdots + u_{s+1}$.
 
In (5.18), terms with $y^{s+\lambda}$ are obtained as follows.
Let $B_k$ be the coefficient of the term 
$x^{s+\lambda-k}y^{s+\lambda}$.\\
(1) For $p = 0$,
 since we smooth all crossings at $v_i$, 
  we have only one term
$\wti{F}_{\mu_1}(xy-1)^s$.
Since $\mu_1 = \lambda$, we have one term 
$x^{\lambda+s}y^{\lambda+s}$.\\
(2) For $p=1$,
we have the following polynomial
\begin{equation*}
(-1)\sum_{i=1}^s v_i(x-1)(y-1) \widetilde{F}_{\mu_1}
 \widetilde{F}_{\mu_2} (xy-1)^{s-1}.
 \end{equation*}
Thus the contribution to $B_0$  by these polynomials is
$(-1)\sum_{i=1}^s v_i$.\\
(3) For general $p$,
the contribution to $B_0$ by the polynomials in
 
    $(-1)^p \sum v_{i_1} v_{i_2}
     \cdots v_{i_p} (x-1)^p(y-1)^p 
    \wti{F}_{\mu_1} \wti{F}_{\mu_2} \cdots 
    \wti{F}_{\mu_{p+1}}(xy-1)^{s-p}$\ 
is \

  ${\displaystyle
    (-1)^p \sum_{1\le i_1<\cdots<i_p\le s}
    v_{i_1} v_{i_2} \cdots v_{i_p}
    }$.
 
Therefore, by letting $n=s$ and $k=0$ in \Lemma \ref{lem:7.3},
we have:
%(5.19)  
\begin{align}
B_0 &= 1 - \sum_{i=1}^{s}v_i+\sum_{1\le i<j\le s}v_i v_j+ 
\cdots\nonumber\\
& \ \ \ +
(-1)^p\sum_{\Lambda_{p,s}}v_{i_1} v_{i_2} 
\cdots v_{i_p} + \cdots + (-1)^s v_1 v_2 
\cdots v_s\nonumber\\
&=(-1)^s( v_1 - 1)( v_2 - 1) \cdots (v_s - 1).
\end{align}
If $\rho = 0$, i.e., $v_j\neq 1$
for any $j$, then $B_0 \neq 0$, 
i.e. $x^{\lambda+s}y^{\lambda+s}$ does exist.
However, if $\rho >0$, then $B_0=0$, and hence
$x^{\lambda+s} y^{\lambda+s}$ does not exist.
Next we consider $B_1$.  The terms 
$x^{\lambda+s - 1 }y^{\lambda+s}$ are obtained as follows.\\
(1) If $p =0$, we do not get the term 
$x^{\lambda+s - 1 }y^{\lambda+s}$.\\
(2) Suppose $p \ge1$. Then in order to get 
$x^{\lambda+s - 1 }y^{\lambda+s}$, 
we must take every possible $y$-term of maximal degree.  
In other words, from each $\wti{F}_{\mu_j}$,
take $(xy)^{\mu_j}$ and 
$(xy)^{s - p}$ from $(xy - 1)^{s - p}$ and 
$y^p$ from $(y - 1)^p$.  
For the $x$-terms we take  
$(-1)^{}\bi{p}{p-1} x^{p-1}$ from $(x-1)^p$.\\
Therefore, we have, by \Lemma \ref{lem:7.3},
%(5.20) 
\begin{align}
B_1&=(-1)\sum_{1\le i_1\le s} v_{i}(-1)\bi{1}{0}+
(-1)^2\sum_{\Lambda_{2,s}} v_{i_1}v_{i_2}(-1)\bi{2}{1}\nonumber\\
&\ \ \ +(-1)^3\sum_{\Lambda_{3,s}}v_{i_1}v_{i_2}v_{i_3}(-1)
\bi{3}{2}+\cdots+
(-1)^{s}v_1\cdots v_{s}(-1)\bi{s}{s-1}\nonumber\\
&=(-1)\Bigl\{s(v_{1}-1)\cdots(v_{s}-1)+
\sum_{1\le i_{1}<\cdots<i_{s-1}\le s}(v_{i_1}-1)
\cdots(v_{i_{s-1}}-1)
\Bigr\}.
\end{align}
If $\rho = 1$, then the first term in the RHS is 
$0$, but one term in the second summation survives.  
Thus, $x^{\lambda+s - 1 }y^{\lambda+s}$ does exist, and

$B_1 =-(v_1 - 1)(v_2 - 1) \cdots (v_{t-1} - 1)
(v_{t+1} - 1) \cdots (v_s -1)$
for some $t$.\\
However, if $\rho \ge 2$, then $B_1 = 0$.
By the same argument, we can show;
%(5.21)
\begin{align}
B_r
&=(-1)^r\sum_{\Lambda_{r,s}} v_{i_1}\cdots v_{i_r}
(-1)^r\bi{r}{0}\nonumber
+(-1)^{r+1}\sum_{\Lambda_{r+1,s}}v_{i_1}\cdots v_{i_{r+1}}
(-1)^r\bi{r+1}{1}\nonumber\\
&\ \ \ +(-1)^{r+2}\sum_{\Lambda_{r+2,s}}v_{i_1}\cdots v_{i_{r+2}}
(-1)^r\bi{r+2}{2}+\cdots
+(-1)^s v_1\cdots v_{s}(-1)^r \bi{s}{s-r}\nonumber\\
&=(-1)^{s+r}\Bigl\{
\bi{s}{r}v_1 \cdots v_{s} - \bi{s-1}{r}
\sum_{\Lambda_{s-1,s}}v_{i_1}\cdots v_{i_{s-1}}\nonumber\\
&\ \ \ \ \ \ 
+\bi{s-2}{r}\sum_{\Lambda_{s-2,s}}v_{i_1}\cdots v_{i_{s-2}}
-\cdots\nonumber\\
&\ \ \ \ \ \ 
+(-1)^{s-r-1}\bi{r+1}{r}\sum_{\Lambda_{r+1,s}}v_{i_1}\cdots
v_{i_{r+1}}+(-1)^{s-r}\bi{r}{r}\sum_{\Lambda_{r,s}}v_{i_1}
\cdots v_{i_{r}}\Bigr\}\nonumber\\
&=(-1)^{s+r}\Bigl\{\bi{s}{r}(v_1-1)\cdots (v_{s}-1)
+\bi{s-1}{r-1}\sum_{\Lambda_{s-1,s}} (v_{i_{1}}-1)\cdots
(v_{i_{s-1}}-1)+\cdots\nonumber\\
&\ \ \ \ \ \ +\bi{s-r}{0}\sum_{\Lambda_{s-r,s}}(v_{i_1}-1)\cdots
(v_{i_{s-r}}-1)\Bigr\}
\end{align}
Thus, if  $\rho \ge r+1$, then $B_r = 0$.  However,
if $\rho =r$, say  $v_1=v_2 =\cdots=v_r=~1$, but 
$v_j\neq 1, j\ge r+1$, then only the last summation 
contains one non-zero term:  
$(v_{r+1} -1)\cdots (v_s -1) \neq 0$.  
Therefore, if $\rho = r$, then 
$B_0=B_1=\cdots=B_{\rho-1}=0$, but there exist  
$s - \rho$ integers  $v_{i_1}, v_{i_2},
\cdots, v_{i_{s-\rho}}$, 
each of which is not $1$, and
%(5.22) \ \ \ 
\begin{equation}
B_\rho = (-1)^{s+r}(v_{i_1} - 1)( v_{i_2} - 1) 
\cdots( v_{i_{s-\rho}} - 1) \neq 0.
\end{equation}
This proves \Proposition \ref{prop:7.1}.
\qed
 
\medskip
{\bf 5.7. Precise form of $\Delta_B(x,y)$}
 
Now we arrive at our final theorem of this section.

\begin{thm}\label{thm:7.4} %Theorem 7.4.
 Let S = $\{P_1, d_1, Q_1, e_1, P_2, d_2, Q_2, e_2, 
 \cdots, P_m, d_m, Q_m\}$  be the canonical 
 decomposition of the continued fraction of $\beta/2\alpha$.
  Let $\rho$ and $\lambda$ be the numbers defined in 
  Definition \ref{dfn:3.8}.  Write
  \begin{equation*} 
 \Delta_{B(2\alpha,\beta)}(x,y) = f_{\lambda - 1}(x) 
 y^{\lambda - 1} + \cdots + f_0(x),
 \end{equation*}
where  $f_{\lambda - 1}(x)  \neq 0$ and $f_0(x) \neq 0$, 
and $f_i(x), 0 \le i \le \lambda - 1$, are integer polynomials.  
Then we have:
%(5.23) &\\
\begin{align}
(1)\ &f_i(x^{-1}) x^{\lambda-1} =
 f_{\lambda-1-i}(x), 0\le  i\le \lambda - 1,\nonumber\\
(2)\ &f_{\lambda - 1}(x) =
\gamma_{\lambda-1,\lambda-1-\rho} 
x^{\lambda-1-\rho} + \cdots +\gamma_{\lambda-1,\zeta} x^{\zeta},\nonumber\\
 &f_{\lambda - 2}(x)  =
\gamma_{\lambda-2,\lambda-1-\rho+1} 
x^{\lambda-1-\rho+1} + \cdots,\nonumber\\
 & \ldots \nonumber\\
 & f_{\lambda - i - 1 }(x)  =
\gamma_{\lambda- i-1,\lambda-1-\rho+i} x^{\lambda-1-\rho+i} + 
\cdots \nonumber\\
 &  \ldots \nonumber\\
& f_{\lambda - 1 - \rho}(x)  =\gamma_{\lambda-1 - 
\rho,\lambda-1} x^{\lambda-1} + \cdots,\ {\it where\ } \gamma_{\lambda-1,\lambda-1-\rho} =
\gamma_{\lambda-1-\rho,\lambda-1} \neq 0, \nonumber\\ 
& {\it and\  hence,\ } \nonumber\\
 &\deg f_{\lambda-1}(x)=\lambda-1-\rho,
\deg f_{\lambda-1-i}(x) \le \lambda - 1 +i- \rho,  
1 \le i \le \rho - 1\ {\it  and\ }\nonumber\\
 &\deg f_{\lambda-1-\rho}(x) = \lambda - 1.\nonumber\\
(3)\ & \ {\it All\ \mbox{\it non-zero}\ leading\ coefficients\ of}\ 
f_i(x)\ {\it are\ of\ the\ same\ sign}.\nonumber\\
(4)\ & 
\prod_{i=1}^m d_i \prod_{j=1}^{m-1} e_j
\ {\it divides}\ 
\gamma_{\lambda - 1, \lambda - 1 - \rho}\
{\it (and}\ \gamma_{\lambda-1-\rho,\lambda-1})\nonumber\\
(5)\ &
\gamma_{\lambda - 1, \lambda - 1 - \rho}\
{\it (and}\ \gamma_{\lambda-1-\rho,\lambda-1}\ {\it )\
is\ equal\ to}\ \pm 1\ {\it if\ and\ only\ if}\nonumber \\
&
  \begin{tabular}{ll}
(i)&  all $d_i = \pm 1$, $1 \le i \le m$, and\\
 (ii)&  all $e_j = \pm 1$, $1\le j \le m-1$, and\\
(iii) &all $b_{i,k}$, $1 \le i \le m, 1 \le k \le s_i$ \\
 & and all $b'_{j,k},1\le  j\le m, 1 \le k \le q_j$
are either $1$ or $2$.
\end{tabular}
\end{align}
\end{thm}

{\it Proof.}
(1)  Since a 2-bridge link $B(2\alpha,\beta)$ is invertible, we have 
$\Delta_{ B(2\alpha,\beta)} (x^{-1},y^{-1}) 
x^{\lambda-1} y^{\lambda-1}= \Delta_ {B(2\alpha,\beta)}(x,y)$.
This implies:\\
$x^{\lambda-1} y^{\lambda-1} 
\Bigl\{f_{\lambda-1}
(x^{-1}) y^{-(\lambda-1)} + f_{\lambda-2}
(x^{-1}) y^{-(\lambda-2)} + \cdots + f_0(x^{-1})
\Bigr\}$
 
$= f_{\lambda-1}(x^{-1}) x^{\lambda-1} + 
f_{\lambda-2}(x^{-1}) x^{\lambda-1}y
 + \cdots + f_0(x^{-1}) x^{\lambda-1} y^{\lambda-1}$,
and hence, we have (1).

(2) \Proposition \ref{prop:7.1}
 shows that $f_{\lambda-1}(x)$ is a required form.  
 Since $B(2\alpha,\beta)$ is 
interchangeable, we see  $\Delta_B(x,y) = \Delta_B(y,x)$, 
and hence
$\gamma_{\lambda-1,\lambda-1-\rho} 
=\gamma_{\lambda-1-\rho,\lambda-1}$.
 
Next, to show that $\deg f_{\lambda-1-i}\le
\lambda-1+i-\rho, 1 \le i \le \rho-1$, 
we need the following easy lemma.
 
 \begin{lem}\label{lem:7.5}
 The number of terms of 
$\Delta_{B(2\alpha,\beta)}(x,y)$ is exactly $\alpha$.  
In other words, if we write
${
\displaystyle
\Delta_{B(2\alpha,\beta)}(x,y) =\sum_{0\le p,q}
c_{p,q} x^p y^q
}$,
then 
${\displaystyle
\sum_{0\le p,q}|c_{p,q}| = \alpha
}$.
\end{lem}

{\it Proof.}  
The group of $B(2\alpha,\beta)$ has the following 
Wirtinger presentation:\\
$\pi_1(S^3 - B(2\alpha,\beta)) =\langle x,y|R\rangle$,
 where $R=W x W^{-1}x^{-1}$, and \\
 $W=y^{\es_1} 
 x^{\es_2} y^{\es_3}\cdots
 y^{\es_{2\alpha-1}},  \es_i = \pm1$.  
 Therefore, the Alexander matrix  $M$ is of the form:
 
$M=\Bigl[\frac{\partial R}{\partial x}\ \frac{\partial R}
{\partial y}\Bigr]^\phi
=\Bigl[\frac{\partial W}{\partial x}(1-x)+W-1,
\frac{\partial W}{\partial y}(1-x)\Bigr]^\phi
=\Bigl[\frac{\partial W}{\partial y}(1-y)\ \frac{
\partial W}{\partial y}(1-x)\Bigr]^\phi
$
and hence,
$\Delta_B(x,y) =\det \bigl[
\frac{\partial W}{\partial y}\bigr]^\phi$.

Here $\det \bigl[
\frac{\partial W}{\partial y}\bigr]^\phi$
is the sum of $\alpha$ terms, while  
$|\Delta_{B(2\alpha,\beta)}(-1,-1)| = \alpha$, 
and hence no cancellation occurs among these 
$\alpha$ terms.
\qed

Now we return to the proof of (2).
Suppose  $\deg f_{\lambda-1-i}(x) > 
\lambda - 1 +i- \rho$.
Then
   $\deg f_{\lambda-1-i}(t) t^{\lambda-1-i} > 
   2( \lambda - 1)- \rho$.\\ 
Write 
$f_{\lambda-1-i}(x) =
\gamma_{\lambda - 1 -i, k} x^k 
+ \cdots +\gamma_{\lambda-1-i,r} x^r$,
where $k>\lambda - 1 + i- \rho$, and $k\ge r$.
Then by~(1),
 
$f_i(x) = f_{\lambda-1-i}(x^{-1})x^{\lambda-1} =
\gamma_{\lambda - 1 -i, r} x^{\lambda-1-r} + 
\cdots +\gamma_{\lambda-1-i,k} x^{\lambda-1-k}$.
 
Since $\lambda - 1 - r \ge \lambda - 1 - k$,  
$\Delta_{B}(t,t)$ contains the term with degree 
$\lambda - 1 - k + i$.
Since no cancellation occurs  when we set $x=y=t$, 
we see
 \begin{align*}
 \deg \Delta_{B(2\alpha,\beta)}(t,t) &>
2(\lambda - 1) - \rho - (\lambda - 1 -  k + i)\\
& =
 \lambda - 1- i - \rho + k  \\
 &> \lambda - 1- i - \rho + 
 \lambda - 1+ i - \rho \\
 &=2 \lambda - 2 - 2\rho.
 \end{align*}
This contradicts \Proposition \ref{prop:3.11}.  This proves (2).

(3) follows also from the fact that no cancellations 
occur when we set  $x=y=t$ in $\Delta_B(x,y)$.
(4) follows from (5.14).
(5) follows also from (5.14) and (5.22).
 
\Theorem \ref{thm:7.4} is now proved.
\qed
 
\begin{rem}\label{rem:7.6}
It is quite likely that
%(5.24)\hfil 
\begin{equation}
\deg f_{\lambda - 1- i }(x) =  
\lambda - 1+i - \rho,  1 \le i \le \rho - 1.
\end{equation}
\end{rem}

\section{Monic Alexander polynomials}%sec8  Feb 25 Sec 6

In this section, we determine 
when the Alexander polynomial of 
$K(2\alpha, \beta|r), r>0$, is monic. 
We use the results proved in the previous section.
In Subsection 6.1, we deal with the case $\ellkB \neq 0$,
using the
continued fraction of $\beta/2\alpha$.
However, if  $\ellkB = 0$, 
we cannot characterize $K(2\alpha, \beta|r)$ with
monic Alexander polynomials in terms of continued fractions.
We then deal with this case in Subsection 6.2. %%%
Let $\{P_1, d_1, Q_1, e_1, P_2, d_2, Q_2, e_2, 
  \cdots, P_m, d_m, Q_m\}$ be the canonical 
  decomposition of the continued fraction 
  of $\beta/2\alpha$.\\
  Write
 $P_i=[[a_{i,1},b_{ i,1},a_{i,2},b_{i,2}, 
 \cdots, a_{i,s_i}, b_{i,s_i},a_{i,s_i +1} ]]$, 
 $a_{i,j}>0$, and \\ 
 $Q_j= [[-a'_{j,1},-b'_{j,1}, 
 \cdots,-a'_{j,q_j},-b'_{j,q_j},-a'_{j,q_j+1}]]$,  
 $a'_{j,k}>0$.
 
 \medskip
{\bf 6.1. The case $\ellkB >0$.}\

The purpose of this subsection is to state 
algebraic conditions equivalent to that in Theorem \ref{thm:A}.
Namely, we prove the following:

\newpage

%%Thm 6.1 final
\begin{thm} \label{thm:8.4} %%\Thm 8.4.  %624 Thm 8.5 
Suppose $\ell  = \ellkB\neq0$.\\
(1)  Suppose $\ell =r=1$.  
Then,
$\Delta_{K(2\alpha, \beta|1)}(t)$ is monic 
       if and only if, for any $i, j, p, q$,
       
(a)  $d_i, e_j = \pm 1$ and
(b)  $b_{i,k}=b'_{j,p}=2$, \\
(2)   Suppose $\ell  \ge 2$.  
Then for any $r \ge1$,
$\Delta_{K(2\alpha, \beta|r)}(t)$ is monic
      if and only if
      
    (a)  $d_i, e_j = \pm1$ and
    (b) $ b_{i,k}$ and $b'_{j,p}$ are $1$ or $2$.
\end{thm}

Let  $\Delta_{B(2\alpha,\beta)}(x,y)$ be the 
Alexander polynomial of $B(2\alpha,\beta)$.  
Suppose $\ell=\ellkB>0$. 
Then by \Proposition \ref{prop:6.1}
 (1) and \Theorem \ref{thm:7.4}, 
the Alexander polynomial 
$\Delta_K(t)$ of $K=K(2\alpha,\beta|r), r >0$, 
is given by
%(6.1) \hfil %%%
\begin{equation}
\Delta_K(t) = \dfrac{1 - t}{1-t^\ell}
\bigl\{f_{\lambda-1}(t) t^{(\lambda-1)\ell r} + 
f_{\lambda-2}(t)t^{(\lambda-2)\ell r} + \cdots +
f_{0}(t)\bigr\},
\end{equation}
where $\lambda - 1$ is the maximal 
$y$-degree of $\Delta_{B(2\alpha,\beta)}(x,y)$.
 
First we determine the degree of  
$\Delta_{B(2\alpha, \beta)}(t, t^{\ell r})$.

\begin{prop}\label{prop:8.1+}
(1) The highest degree of 
$\Delta_{B(2\alpha, \beta)}(t, t^{\ell r})$ is
%(6.2)
\begin{equation}
\lambda-1-\rho+(\lambda-1)\ell r.
\end{equation}

(2) The lowest degree of 
$\Delta_{B(2\alpha, \beta)}(t, t^{\ell r})$ is $\rho$.
\end{prop}

\noindent
{\it Proof.}
(1)
We write
  $f_{\lambda-1}(x) =
  \gamma_{\lambda-1,\lambda-1-\rho}
  x^{\lambda-1-\rho} + \cdots + 
  \gamma_{\lambda-1,\zeta}x^\zeta,
  \gamma_{\lambda-1,\lambda-1-\rho}\neq0$.\\
We show that if $\ell r \ge 2$, then 
$\gamma_{\lambda-1,\lambda-1-\rho} 
t^{\lambda-1-\rho} t^{(\lambda-1)\ell r}$ 
is the only term with the highest degree in
$\Delta_{B(2\alpha, \beta)}(t,t^{\ell r})$.
In fact, by \Theorem \ref{thm:7.4}, we see that
for  $1\le i\le \rho-1$,
$\deg f_{\lambda-1-i}(x)\le \lambda-1+i-\rho$,
and hence, since $\ell r\ge 1$,  for $1\le i\le \rho-1$
we have:
%(6.3) \hfil %%%
\begin{equation}
\lambda-1-\rho+(\lambda-1)\ell r \ge
\lambda-1+i-\rho+(\lambda-1-i)\ell r.
\end{equation}
Moreover, obviously, $\deg f_j(x)\le \lambda-1$,
for $0\le j\le \lambda-2-\rho$, and hence,
if $\ell r\ge1$, then for $0\le j\le \lambda-2-\rho$, we see
%(6.4) \hfil %%%
\begin{equation}
\lambda-1-\rho+(\lambda-1)\ell r \ge \lambda-1+j\ell r.
\end{equation}

Combining (6.3) and (6.4), we have (1). %%% %%%
In particular, if $\ell r\ge 2$, the strict inequality
holds in (6.3), and therefore, %%%
$\gamma_{\lambda-1,\lambda-1-\rho}t^{\lambda-1-\rho}
t^{(\lambda-1)\ell r}$ is the only term with 
the highest degree in
$\Delta_{B(2\alpha,\beta)}(t,t^{\ell r})$.
However, if $\ell r=1$, then the equality holds in (6.3), and %%%
thus, the terms with the highest degree appear at least in 
$f_{\lambda-1}(t) t^{(\lambda-1)\ell r}$
and
$f_{\lambda-1-\rho}(t) t^{(\lambda-1-\rho) \ell r}$.

(2) First we note from \Proposition \ref{prop:3.11} 
that the lowest degree
of $\Delta_{B(2\alpha,\beta)}(t, t)$
is $\rho$, since the highest degree of 
$\Delta_{B(2\alpha, \beta)}(t,t)$ is $2(\lambda-1)-\rho$ by (1).

By \Theorem \ref{thm:7.4}, we see that
 $\deg f_{\lambda-1-i}(x) \le  \lambda -1+i - \rho,
 1 \le i\le  \rho - 1$ and, of course,
$\deg f_j \le  \lambda - 1$, for  
$0\le j\le \lambda - 2 - \rho$.

 Now, since $f_0(x) = f_{\lambda-1}(x^{-1}) 
 x^{\lambda-1}$, it follows that the lowest degree of 
 $f_{\lambda-1}(t) t^{(\lambda-1)\ell r}$ is 
 $\zeta + (\lambda-1)\ell r$, 
 and that of $f_0(t)$ is $\rho$.
Since $\rho\le\lambda-1$, 
$\min\{\zeta+(\lambda-1) \ell  r,\rho\}=\rho$,
and hence, if $\ell r\ge 1$, $f_0(t)$ contains
the term of degree $\rho$.
Furthermore, since the lowest degree of
$\Delta_{B(2\alpha, \beta)}(t,t)$
is $\rho$, the degree of any term in
$\Delta_{B(2\alpha, \beta)}(t,t^{\ell r})$ is
at least $\rho$. This proves (2).
\qed

\Proposition \ref{prop:8.1+} implies the following:

\begin{prop}\label{prop:8.1} %\Prop 8.1. %624 Prop 8.2 
If $\ell r \ge 1$, then
  $\deg \Delta_{K(r)}(t)=(\lambda-1)(\ell r+1)
  - 2\rho-(\ell -1)$.
In particular, if $\ell r=1$ (i.e. $\ell =r=1$), then
$\deg \Delta_{K(r)}(t)=2(\lambda -\rho-1)$.
\end{prop}

Using the above results, we can characterize the 
monic Alexander polynomial.
 
First, we see that if 
$\ell r\ge2$, then the leading coefficient of
$\Delta_K(t)$ is given by 
$\gamma_{\lambda-1,\lambda-1-\rho}$.
Therefore,$\Delta_K(t)$ is monic if and only if 
$\gamma_{\lambda-1,\lambda-1-\rho}$ is $\pm 1$, 
and hence \Theorem \ref{thm:7.4}(5) gives us 
immediately the following:

\begin{prop}\label{prop:8.2} %\Prop 8.2.  %624 Prop 8.3 
Suppose $\ell r\ge2$.  Then $\Delta_K(t)$ is 
monic if and only if the following conditions hold:

(1)  $d_i, e_j=\pm1$ for any $i,j$, and

(2)  $b_{i,k}= 1$ or $2$ and $b'_{j,p}= 1$ or $2$, 
for any $1\le i\le m$, $1\le k\le s_i+1$, 
and $1\le j\le m, 1\le p\le q_j +1$.
\end{prop}

Note that $a_{i,j}$ and
$a'_{j,k}$ are arbitrary.
 
If $\ell r=1$, then the following proposition holds.

\begin{prop}\label{prop:8.3} %\Prop 8.3.  %624 Prop 8.4 
Suppose $\ell =r=1$. Then $\Delta_{K(1)}(t)$
 is monic if and only if

(1)  $d_i, e_j = \pm1$ for any $i,j$, and

(2)  $b_{i,k}= b'_{j,p}=2$  for any $i,k,j,p$. 
(In particular, $\rho =0$.)
\end{prop}
 
{\it Proof.}  
%The \lq only if' part follows from Theorem \ref{thm:7.4}.
(1) Suppose that $d_i$ or $e_j$ is not $\pm 1$.
Then  $\Delta_{K(1)}(t)$ is not monic by Theorem \ref{thm:7.4}.
(2) Suppose $\rho \neq 0$. Then $\Delta_{B(2\alpha, \beta)}(x,y)$
contains at least two non-zero terms, 
$\gamma_{\lambda-1, \lambda-1-\rho}
x^{\lambda-1-\rho}y^{\lambda-1}$
and $\gamma_{\lambda-1-\rho, \lambda-1}
x^{\lambda-1}y^{\lambda-1-\rho}$.
Since $\gamma_{\lambda-1, \lambda-1-\rho}=
\gamma_{\lambda-1-\rho, \lambda-1}$ by Theorem \ref{thm:7.4}(2),
we see that $\Delta_{K(2\alpha, \beta|1)}(t)$ is not monic.
Further, as is proved in Subsection 5.6, 
we have all $b_{i,k}=b'_{j,k}=2$.
The converse follows from Theorem \ref{thm:7.4}.
\qed
%Since $\ell  r = 1$, (6.4) holds, %%%
%and hence, if $\Delta_K(t)$ is monic, 
%then $\rho =0$. Then \Proposition \ref{prop:8.3} 
%follows from \Theorem \ref{thm:7.4}. %thm 7.5
%\Proposition
%\ref{prop:8.2} (2).
 
 \medskip

By these results above, we obtain Theorem \ref{thm:8.4}. 
%  We collect  these results into one theorem.

\medskip

Finally, we note that 
we can prove the following
(c.f., \cite[Theorem 4.2]{KMS})
as a simple consequence of
\Proposition \ref{prop:8.1}.

\begin{prop}\label{prop:8.6++}
Suppose $\ell, r>0$.
Then $K(2\alpha, \beta|r)$ is unknotted if and only if
$(2\alpha, \beta)=(4,3)$ and $r=1$,
and hence $\ell=2$.
\end{prop}

{\it Proof.}
Since the \lq\lq if\rq\rq\ part is obvious, we only
consider the \lq\lq only if\rq\rq\ part.
Suppose $K(r) = K(2\alpha, \beta|r)$
is unknotted.
Then, by \Proposition \ref{prop:8.1},
we have:
%(6.5) \hfil %%%
\begin{equation}
(\lambda-1)(\ell r+1)-2\rho-(\ell-1)=0.
\end{equation}
Rewrite the LHS of (6.5) as %%%
%(6.6) \hfil %%%
\begin{equation}
2(\lambda-1-\rho)+(\lambda-1)(r-1)\ell+
(\lambda-2)(\ell-1)=0.
\end{equation}
Since $\lambda\ge 2$ and $\lambda-1\ge \rho$,
it follows that
each term of the LHS is non-negative, 
and hence, the equality holds only if 
we have:
(i) $\lambda-1=\rho, r=1$, and
(ii) $\lambda=2$ or $\ell=1$.
Now, since $\lambda-1=\rho$, 
we see that $(2\alpha, \beta)=(2\lambda, 2\lambda-1)$
and $\ellk B(2\alpha, \beta)=\lambda$.
Since $\lambda\ge2$, $\ell$ cannot be $1$,
and hence the conclusion follows.
\qed

%%%\section{Monic Alexander polynomials (II)} %sec 10 New 9

\medskip
{\bf 6.2. The case $\ellkB =0$.}
 
In this subsection, we characterize the
 monic Alexander polynomial  of $K(r)$ when 
$\ellkB=0$.
 To do this, first, we give a geometric interpretation of
$\Delta_{K(r)}(t)$.  We remind that the calculation of the 
Alexander polynomial of $K(r)$ from 
$\Delta_{B(2\alpha,\beta)}(x,y)$ is quite different for 
$\ellkB=0$, 
as is seen in \Proposition \ref{prop:6.1}.
  Let  $B(2\alpha,\beta)$ be a 2-bridge link 
  consisting of $K_1$ and $K_2$.  Suppose that 
  $\ell k(K_1,K_2)=0$.  
  Consider the infinite cyclic 
  cover $M^3$ of $S^3\setminus K_2$.  
  Since $K_2$ is unknotted, $M^3$ is an infinite cylinder
  $D^2 \times {\mathbb R}^1$.
  Let $\{\wti{K}_m, m=0, \pm1,\pm2,\ldots\}$ be the set of 
  lifts of $K_1$ in $M^3$, where 
  $\wti{K}_j=\psi^{j}(\wti{K}_0)$ with the covering 
  translation $\psi$.
  Since $\ell k(K_1,K_2)=0$, 
  each lift $\wti{K}_m$ 
  is a knot in $M^3$ with orientation inherited from that of 
  $K_1$.
Denote $c_j =\ell k(\wti{K}_0,\wti{K}_j), j\neq 0$, and let
%(6.7) \hfil %%%
\begin{equation}
\Gamma(t) = \sum_{-\infty<j<\infty} c_jt^j,\ 
{\rm where}\ c_0 = - \sum_{-\infty<j<\infty} c_j.
\end{equation}   
Note that $c_0$ is well defined, and also $c_j=c_{-j}$ 
for any $j$.  Then it is proved in \cite{MM}

\begin{prop}\label{prop:10.1} %\Proposition 10.1 %new 9.1
    $\Gamma(t) \doteq (t-1)\Bigl[\dfrac{\Delta_{B}(x,y)}
    {1-y}\Bigr]_{x=t,y=1}$
\end{prop}

%\noindent%\begin{rem}\label{rem:10.2} 
%{\bf Remark.} 
We note that
Gonzalez-Acu$\tilde{{\rm n}}$a 
also studied this polynomial
$\Gamma(t)$ in \cite{Gon}.

\medskip
  Using \Proposition \ref{prop:10.1}, 
we can estimate the 
maximal and minimal degree of  
$\Delta_{K(2\alpha,\beta|r)}(t)$.\\ 
Let $S = [[ u_1,v_1,u_2,v_2,\ldots,u_s,v_s,u_{s+1}]]$ 
be the continued fraction of $\beta/2\alpha$.
Let $G(S)$ be the graph of $S$.  
  This graph $G(S)$ will be used to estimate the degree of 
  $\Delta_{K(2\alpha,\beta|r)}(t)$.

 \begin{prop}\label{prop:10.2} %Proposition 10.2. %new 9.2
  Let $h$ and $q$ be the highest and lowest 
  $y$-coordinates of $G(S)$.  Then
 $\deg \Delta_{K(2\alpha,\beta|r)}(t) \le 2 \max\{h,|q|\}$.
 \end{prop}
 
Note that $h\ge 0$ and $q\le0$.
 
{\it Proof.}  
We span $K_2$ by a disk $D$ in such a way that 
$K_1$ intersects $D$ transversally 
at $\sum_{i=1}^{s+1}|u_i|$  points.  
Using $D$, we construct $M^3$.  
Then it is easy to evaluate the linking number 
between $\wti{K}_0$ and $\wti{K}_j$ for 
$j \ge1$ using the 
primitive disk for $K_1$. See Example \ref{ex:lzero}.

%in \Figure 9.1.
% \Fi{9.1}

Also we can easily determine $h$ and $q$ from the graph 
$G(S)$.  In fact, $h$ is the $y$-coordinate of the 
absolute maximal vertices of $G(S)$, and $q$ is the 
$y$-coordinate of the absolute minimal vertices.  
Let $V_{i,1}, V_{i,2},\ldots,V_{i,p}$ be the absolute maximal 
vertices and $V_{j,1},V_{j,2},\ldots,V_{j,s}$ 
be the absolute minimal vertices of $G(S)$.
 
 Let $w_{i,k}$ be the weight of $V_{i,k},
 k=1, 2, \cdots, p$,
 and $w_{j,n}$ the weight of $V_{j,n},
 n=1,2, \cdots, s$.
 Then we have:
% (6.8) & %%%
\begin{align}
&(1)\ {\rm If}\ h>|q|,\ {\rm then}\ 
\ellk(\wti{K}_0,\wti{K}_h)=
 -\frac{1}{2}\sum^p_{k=1}
 w_{i,k}.\nonumber\\
 &(2)\  {\rm If}\ h<|q|,\ {\rm then}\ 
\ellk(\wti{K}_0,\wti{K}_q)=
-\frac{1}{2}\sum^s_{n=1}
w_{j,n}.\nonumber\\
&(3)\ {\rm If}\ h=|q|,\ {\rm then}\ 
\ellk(\wti{K}_0,\wti{K}_h)=-\frac{1}{2}
 \left\{
 \sum^p_{k=1}w_{i,k}+\sum^s_{n=1}w_{j,n}
 \right\}.\nonumber\\
&(4)\ {\rm If}\ d>\max\{h,|q|\},\ {\rm then}\
\ellk(\wti{K}_0,\wti{K}_d)=0.
\end{align}

Therefore,
$\max \deg \Gamma(t) \le \max\{h, |q|\}$, 
              and
$\min \deg \Gamma(t) \ge - \max\{h,|q|\}$,
and hence,
$\deg \Delta_{K(2\alpha,\beta|r)}(t) 
\le 2 \max\{h,|q|\}$.
 \qed
 
Under the same notation used in the proof of \Proposition 
\ref{prop:10.2}, we have

\begin{cor}\label{cor:10.3} %\Corollary 10.3. %new 9.3
 $\biggl[\dfrac{\Delta_{K(2\alpha,\beta)}(x,y)}{1-y}
 \biggr]_{\genfrac{}{}{0pt}{}
{x=t}{y=1}}(1-t)$ 
 is monic  and its degree is equal to
 $2\max\{h,|q|\}$
 if and only if

(1) $\sum_{k=1}^{p} w_{i,k}=\pm 2$, when $h>|q|$,

(2) $\sum_{n=1}^s w_{j,n} = \pm 2$, when $h<|q|$,

(3) $\sum_{k=1}^p w_{i,k} +\sum_{n=1}^s w_{j,n} 
= \pm 2$,
when $h=|q|$.
\end{cor}

\begin{prop}\label{prop:10.4} %\Prop 10.4. %new 9.4
Suppose $\ellkB=0$.  
Then, for $r >1$,

(1)    $\Delta_{K(2\alpha,\beta|r)}(t)$ is 
either non-monic or  $\Delta_{K(2\alpha,\beta|r)}(t) = 1$.
 
(2) $\Delta_{K(2\alpha,\beta|1)}(t)$ is monic if and only if
 $\biggl[\dfrac{\Delta_{K(2\alpha,\beta)}(x,y)}{1-y}
 \biggr]_{\genfrac{}{}{0pt}{}
{x=t}{y=1}}$ 
 is monic.
\end{prop}

This is an immediate consequence of 
\Proposition \ref{prop:6.1}.

\medskip
The rest of this paper (except for the last three sections)
will be devoted to the proofs of our main theorems.

\section{Construction of a Seifert surface $F_1$ for $K_1$.} 
%%%
By Theorem \ref{thm:1}, we assume $\ellkB =\ell
\ge0, r>0$.
In section 3, %%%
we constructed a primitive spanning disk 
$F_D$ for $K_1$, which consists of disks and bands
corresponding to the edges 
and vertices in $G(S)$. 
Recall that $F_D$ intersects $K_2$ 
as many times as the number of edges in $G(S)$,
which is equal to $\lambda$.
In this section, we construct a Seifert surface
$F(r)$ for $K(2\alpha, \beta|r)=K(r)$.
First, 
using $G(S)$,
we construct a new Seifert surface $F_1$ for $K_1$
which intersects  $K_2$ exactly
$\ell \ge0$ times.
We call $F_1$ a {\it canonical surface} for $K_1$.

Let $S=\{P_1,d_1,Q_1,e_1,P_2, \ldots \}$ be the 
canonical decomposition of the continued fraction of 
$\beta/2\alpha$.  Let $G(S)$ be the graph of $S$,
which by definition is the graph $G(S^*)$ of
the modified continued fraction $S^*$ of $S$.
We construct a new surface $F_1$ by induction on 
$\nu(G(S))$, 
the total number of local maximal and local minimal 
vertices,
including the end vertices.

Case 1: $\nu(G)=2$

Since $\ell\ge 0$, 
$G$ is an ascending
line segment. See \Figure 7.1 (a).

\tena %\Fi{10.1}

In this case, $F_D$ itself is our new surface $F_1$,
which consists of two disks connected by a twisted band.

Case 2: $\nu(G)=3$

 There are two cases. See \Figure 7.1 (b1) and (b2).
 Note that the vertex $B$ may be on the $x$-axis.
 
 %\tenb %\Fi{10.2}
 
 For the first case (b1), $F_D$ consists of $p$, say,
 positive disks $D_1, D_2, \ldots, D_p$ followed by $q$,
 say, negative disks
 $D'_1, D'_2, \ldots, D'_q, p\ge q$, and $p+q-1$ bands
 $B_j, j=1,2, \ldots, p+q-1$, connecting
these disks. See \Figure 7.2 (a).

\tenc 
%\Fi{10.3}

Then replace two disks $D_p$ and $D'_1$ 
by a cylinder $R_1$, where 
$R_1\cap F_D=\partial R_1=\partial(D_p \cup D'_1)$.
The orientation of $R_1$ is naturally induced from 
those of disks.
Next we replace two disks $D_{p-1}$ and $D'_2$ 
by a cylinder $R_2$ that is inside of $R_1$ and 
$R_2\cap F_D=\partial R_2=\partial(D_{p-1} \cup D'_2)$. 
Repeat this operation for every pair of a positive 
disk $D_{p-i+1}$ and a negative disk $D'_i$, $1\le i\le q$, 
in the same manner, 
so that we have a sequence of cylinders.  
Positive disks and bands corresponding to the subgraph $OB'$ 
are untouched as in Case 1.
This untouched part of $F_D$ and the cylinders and all
bands form our new surface $F_1$.
See \Figure 7.2 (b).
For the second case (b2), $F_1$ is constructed in the same
manner, and it looks like a surface depicted in 
\Figure 7.2 (c).
We note that
in this construction, all bands are untouched, 
and therefore, 
we are only concerned with disks.

Case 3: $\nu(G)\ge 4$.

Subcase (i): The origin $O$ is local minimal.

Then a local maximal vertex $A$ is followed 
by a local minimal vertex $B$.
See \Figure 7.3 (a1) and (a2).

Subcase (ii): The origin $O$ is local maximal.

Then a local minimal vertex $A$ is followed by a local maximal vertex 
$B$. See \Figure 7.3 (b1) and (b2).

\tend
\medskip

In, Subcase (i) (a1),
first, apply the argument used in Case 2 (b1) 
on the subgraph $OA\cup AB$
of $G(S)$, and replace pairs of disks corresponding to the edges on
$B'A$ and $AB$ by cylinders. Secondly, delete the subgraph 
$B'A\cup AB$ from $G$ and then identify $B'$ and $B$ to obtain a new
graph $G'$. 
Since $\nu(G')=\nu(G)-2$, 
we can inductively construct a surface $F'_1$ in such a way that
all cylinders in $F'_1$ are inside the cylinders we previously
constructed.
Our new surface $F_1$ is the union of cylinders firstly 
constructed and $F'_1$, (and all bands).

As an example, in \Figure 7.4 below, 
we depict a sequence of 
modifications of $G$ and corresponding surfaces.
The last surface is the surface $F_1$ we sought.

In Subcase (i) (a2),
apply the argument used in Case 2 (b1) on the subgraph $OA\cup AB'$,
and replace pairs of disks by cylinders.
Then delete the subgraph $OA\cup AB'$ from $G(S)$ and
identify $O$ and $B'$ so that a new graph $G'$ is obtained. 
Since $\nu(G')=\nu(G)-1$, apply induction on $G'$.

In Subcase (ii) (b1) and (b2),
apply the argument used in Case 2 (b2) and repeat similar arguments
used in Subcase (i) (a1) and (a2).

%\Fi{10.6}
\tenf

\begin{ex} \label{ex:9.1} %Example 9.1. new 10.1
\Figures 7.5 and 7.6 depict graphs and 
corresponding canonical surfaces.
\end{ex}

\teng

\tenh

 To construct $F_1$, we start with a primitive spanning
 disk %(Definition \ref{dfn:primitivedisk}) Mar3 2006
 and the genus increases by $1$ each time
 we replace a pair of small disks by a cylinder. 
  We have
$(\lambda-\ell)/2$ pairs of disks to replace.
Therefore, we have the following:

\medskip

 \begin{prop}\label{prop:10.last} %10.2
 Let $\lambda$ be the number of edges in $G(S)$.
 The number of disks in $F_1$ constructed above is
 $\ell=\ellkB$,
 and the number of cylinders of $F_1$ is $\frac{1}{2}(\lambda-\ell)$,
 therefore,  $g(F_1)= \frac{1}{2}(\lambda-\ell)$.
 \end{prop}
 
Now we twist $F_1$ by $K_2$.
First, suppose $\ell=0$. 
Then just by twisting, we obtain a Seifert surface $F(r)$
for $K(r)$. In Section 10, (10.1),
we show that $g(K(r))=g(F(r))$, which is equal to the number
of cylinders, i.e., $(\lambda-\ell)/2$.

Next, suppose $\ell \neq 0$.
$F_1$ consists of $\ell$ disks and 
$(\lambda-\ell)/2$ cylinders, connected by bands,
and $g(F_1)= \frac{1}{2} (\lambda-\ell)$
(by \Proposition \ref{prop:10.last}).
If we twist $F_1$, $r$ times by $K_2$, 
we obtain a singular surface, 
in which cylinders penetrate the $\ell$ disks transversely.
Remove ribbon singularities by smoothing intersections 
in the standard way. 
Then we obtain a Seifert surface $F(r)$ for
$K(r)$.
See \Figure 7.7 for the case $r=1$.

For a technical reason, we need more specific description
(given at the end of this section) 
on the position of bands connecting  disks
and cylinders.

\fotnA

Each time we make a hole, 
the genus of the surface is increased by one.  
We see that there are exactly 
$\frac{1}{2}(\lambda-\ell)\ell r$ intersections.  
Furthermore, by twisting along $K_2$, the boundaries of $\ell$ disks 
form a torus link of 
type  $(\ell, \ell r)$. %whose genus is $\frac{1}{2}(\ell-1)\ell r$.  

Therefore, we have:
%(7.1)  \hfil
\begin{align}
g(F(r))&=\frac{1}{2}
(\lambda-\ell) + 
\frac{1}{2}(\lambda-\ell)\ell r + 
\frac{1}{2}(\ell-1)\ell r\nonumber\\ 
&=  \frac{1}{2}\bigl\{(\lambda-1)(\ell r +1) - (\ell-1)
         \bigr\}.
\end{align}
         
In the proof of Theorem \ref{thm:C}, in the next section,
we show  (i) if $\rho=0$, then $F(r)$ is a minimal genus
Seifert surface for $K(r)$, 
and (ii) if $\rho>0$, then $F(r)$ admits
compressions $\rho$ times and the result is a minimal genus
Seifert surface for $K(r)$, where $\rho$ is the deficiency
(Definition \ref{dfn:3.8}).

Now we give a precise description of the relative
position of bands connecting
disks and cylinders. 
Proposition \ref{prop:4.2} also implies
that, in $F_1$,  we have a freedom of relative positions
of the bands.
However this freedom
is lost when we have twisted $F_1$ to obtain $F(r)$.

A band $B$ in $F_1$ 
is of one of the five types below. See \Figure 7.8. 
Note that the boundary of 
($F_1 -$ bands) consists of circles. In Figure 7.8,
a circle with an arrow
 heading toward left (resp. right) 
corresponds to a rising (resp. falling) edge of $G$.
Each band corresponds to a vertex in $G$.

Type I: 
$B$ connects the boundary of an outermost cylinder.

Type II:
$B$ connects two stacked cylinders.

Type III:
$B$ connects two cylinders side by side 
(of the same sign
or the opposite sign). 

Type IV:
$B$ connects a disk and a cylinder. 
In this type, we have three subtypes
as in Figure 7.8.3.
%to the sign of the disks and cylinders connected;
%(i) a positive disk and a positive cylinder,
%(ii) a positive disk and a

Type V:
$B$ connects two positive disks. %\Figure 14.2 (e).

\constAa

\constAb

\constAc

\constAd

Note. 
In Figure 7.8.3 (c), there are no
disks above the depicted cylinder
by our construction of $F_1$. 
%In the right figure
%of Figure 7.8.3, the bands signed with $*$ are obtained
%in the process of twisting to obtain $F(r)$.

To obtain $F(r)$, we place the bands of Type V
as in Figure 7.8.4, so that they are placed close to
the bands emerged by twisting.
The bands of Types I, II and III should be
arranged as in Figure 7.9, where the horizontal disk
is the very top disk in $F_1$, i.e, it corresponds
to the edge in $G$ that is the first rising edge
after the last intersection of $G$ and the $x$-axis.
Bands of Type IV, where a positive cylinder is connected (Figure 7.8.3 (a), (b)),
should be arranged as in Figure 7.7. 
If $r\ge 2$,
then the first band (resp. the second band, if any)
is placed before (resp. after) the $r$-twists.
A band of Type IV, where a negative cylinder is connected (Figure 7.8.3 (c)) is similarly done
as Figure 7.8.3 (b). See Figure 8.5 for
a local picture for this type of band
after twisting by $K_2$.

Now we have constructed a Seifert surface
for $K(r)$.

\constB

\medskip

\secti{Proof of \Theorem \ref{thm:C}} %non zero genus
 
%In this section, we prove Theorem \ref{thm:C}.

Let $S = \{P_1,d_1,Q_1,e_1,P_2,d_2,Q_2,e_2,\ldots, P_m,d_m,Q_m\}$
 be the
canonical decomposition of  the continued fraction of 
$\beta/2\alpha$.  
Express each $P_i$ and $Q_j$ 
by modified continued fractions, 
and write,
\begin{align*}
&P_i =[[1,b_{i,1},1,b_{i,2},1,\ldots,1,b_{i.s_i},1]],   
   1\le i\le m,\nonumber\\     
&Q_j= [[-1,-b_{j,1}',-1,-b_{j,2}',-1,
   \ldots,-1,-b_{j,q_j}',-1]],  
   1\le j\le m,
\end{align*}
where $b_{i,k}$ and $b_{j,k}'$ are arbitrary, and may be $0$.
Since $\ell \neq 0$, we see from \Proposition \ref{prop:6.1},
 %(8.1) 
\begin{equation}
\Delta_{K(r)}(t) = \frac{1-t}{1-t^\ell}
\Delta_{B(2\alpha,\beta)}(t,t^{\ell r}),
\end{equation}
where $\ell >0$ and $r>0$.
Then, by \Proposition \ref{prop:8.1}, we have:
\begin{equation*}
\deg \Delta_{K(r)}(t) = (\lambda-1)(\ell r+1)-(\ell-1)-2\rho.
\end{equation*}

Recall that in Section 7, (7.1), 
we constructed a Seifert surface $F(r)$
for $K(r$) with $g(F(r))
=  \frac{1}{2}\bigl\{(\lambda-1)(\ell r +1) - (\ell-1) \bigr\}$.

If $\rho =0$, then $F(r)$ 
is a minimal genus Seifert surface for $K(r)$, 
since $g(F(r))= \frac{1}{2}\deg \Delta_{K(r)}(t)$.
Therefore, to prove \Theorem \ref{thm:C},
%\Proposition \ref{prop:15.1}, 
it suffices to confirm that we can compress $F(r)$
as many times as $\rho$ (Definition \ref{dfn:3.8}).
%$b_{i,k}=1$ and $b_{j,k}'=1$ in $S$, 
In \Proposition \ref{prop:15.2} below,
we demonstrate where to apply compression
corresponding to each $b_{i,k}=1$ and $b_{j,k}'=1$ in $S$.
Since we need it in the proof of Theorem \ref{thm:A},
we also show
where we can deplumb twisted annuli from
$F(r)$. But first, we apply compressions.

\begin{prop}\label{prop:15.2}   
We can compress $F(r)$ $\rho$ times, where
each compression 
corresponds to an occasion of
 $b_{i,k} = 1$ or $b'_{j,k} = 1$.
Furthermore, 
we can deplumb $\sum_{i=1}^{m} s_i +
\sum_{j=1}^{m} q_j + 2m-1-\rho$ unknotted,
twisted annuli from $F(r)$,  where
each deplumbing corresponds
to an occasion of
$b_{i,k} \neq 1, b'_{j,k} \neq 1, d_i $ and $e_i$.
\end{prop}

%%{\it Proof.}
%%Let $B$ a band corresponding to $b_{i,k}=1$ or
%%$b'_{j,k}=-1$, i.e., contributing to $\rho$.
%%Then there are two cases:\\
%%(i) $B$ connectes two cylinders in $F(r)$, and
%%(ii) $B$ connects two disks in $F(r)$.
%%
%%
%%
%%
%%corresponding to each 
%%$b_{i,k}=1$ and $b_{j,k}'=1$ in $S$, 
%%we can deduce the genus of $F(r)$ by one.
%%In \Proposition \ref{prop:15.2} below
%%we demonstrate where to apply compression.
%%Since we need it in the proof of Theorem \ref{thm:A},
%%we also show where we can deplumb twisted annuli from
%%$F(r)$.
%%
%%First we express each $P_i$ and $Q_j$ 
%%by modified continued fractions, 
%%and write,
%% 
%%   $P_i =[[1,b_{i,1},1,b_{i,2},1,\ldots,1,b_{i.s_i},1]],   
%%   1\le i\le m$,     
%% 
%%   $Q_j= [[-1,-b_{j,1}',-1,b_{j,2}',-1,\ldots,-1,b_{j,q_j}',-1]],  
%%   1\le j\le m$,
%% 
%%where $b_{i,k}$ and $b_{j,k}'$ are arbitrary, and may be $0$.
%% 
%
%
%%\begin{prop}\label{prop:15.2}   
%%Given the surface  $\hat{F}(r)$, 
%%we can construct a Seifert surface 
%%$F^*(r)$ such that the minimal genus Seifert surface
%%$K(r)$ is 
%%Murasugi sum of $F^*(r)$ and  $\sum_{i=1}^{m} s_i +
%%\sum_{j=1}^{m} q_j + 2m-1-\rho$ unknotted,
%%twisted annuli (that span elementary torus links).
%%\end{prop}

{\it Proof.}  
At each band connecting disks and cylinders,
we explicitly show how to apply either compression or
deplumbing of an annulus. 
Each deplumbing corresponds to removing a band,
and each compression corresponds to cutting
the surface along a properly embedded arc.

If a band is of Type I, it is obvious that
we can deplumb an unknotted annulus with
$d_i$ or $e_j$ full twists. Compressions never occur
for this type.
For a band of Type II, the relevant part of $F(r)$ is depicted in Figure 8.1.
%%Then we see that $F_1$ is obtained from $F_1- B$ 
%%by plumbing an
%%annulus $A$ with $b-1$ (or $-(b'-1)$) full-twists
%%(\Figure 14.2 (2) (i) or (ii)) along the curve given in 
%%\Figure 14.2 (1).
If $b\neq 1$ or $b' \neq 1$, then we can deplumb
a twisted annulus and thus remove a band.
See \Figure 8.1 (a).
In particular,
if $b$ or $b'$ is either $0$ or $2$,
then the annulus is a Hopf band.
However, if $b=1$ or $b'=1$, 
then the annulus yields a compressing disk, so we
do not deplumb a band, but apply compression
(see Figure 8.1 (b)). 

\fotnC

Take a band $B$ of Type III. If $B$ connects two
cylinders showing the same side (see Figure 8.2), 
we can deplumb a twisted annulus. 
In this case, compression never occurs.
In particular we can deplume a Hopf band if
$d=\pm1$ or $e=\pm1$.
On the other hand, if $B$ connects two cylinders 
showing the opposite sides (see Figure 8.3), 
then we can apply compression if $b$ or $b'$
equals $1$,
and otherwise deplumb a twisted annulus.

\fotnD
\fotnE

Take a band $B$ of Type IV (Figure 7.8.3). 
There are three subtypes according
to the feature of the corresponding vertex $v_B$ in $G$:
(a) $v_B$, not on the $x$-axis
is between two rising edges of $G$,
(b) $v_B$ is a local minimum, and
(c) $v_B$, on the $x$-axis, is between two rising edges.

For (a), see Figure 8.4. If $b=1$, then we can compress,
and otherwise we can remove the band by deplumbing an
annulus with $b-1$ full twists. (See Figure 7.7.)
For (b), we can deplumb a band with $e$ full twists,
and compression never occurs.
For (c), see Figure 8.5. 

\fotnF %8.4
\fotnG

For bands of Type V, 
deform $F(r)$ by isotopy as in Figure 8.6
so that each band of Type V is adjacent to
a band emerged by twisting $F_1$.
Then we can remove the band by deplumbing if $b \neq 1$,
and otherwise cancel it with its neighbor, which corresponds
to compressing $F(r)$.
Note that since $r\ge1$, even if we cancel all bands of
Type V, still there are bands connecting each 
pair of adjacent disks.
\qed

\fotnH

\medskip
By \Proposition \ref{prop:15.2}, the proof of
\Theorem \ref{thm:C} is now completed.

\medskip
%\secti{Reformulation of \Theorem \ref{thm:A}}%5
\secti{Proof of Theorem \ref{thm:A}}

By Theorem \ref{thm:1},
we assume,  $\ellkB=
\ell\ge0$ and $r>0$.

Let $\beta/2\alpha=[[c_1, c_2, \ldots, c_{2d+1}]]$ 
be a continued fraction of $\beta/2\alpha$, and
$S=$\\
$\{P_1, d_1, Q_1, e_1, P_2, \ldots, P_m, d_m, Q_m\}$
be the canonical decomposition of $S$.
Write
\begin{align*}
&P_i=[[a_{i,1},b_{ i,1},a_{i,2},b_{i,2}, 
 \cdots, a_{i,s_i}, b_{i,s_i},a_{i,s_i +1} ]], a_{i,j}>0,\
 {\rm  and}\\ 
&Q_j= [[-a'_{j,1},-b'_{j,1}, 
 \cdots,-a'_{j,q_j},-b'_{j,q_j},-a'_{j,q_j+1}]],
 a'_{j,k}>0.
 \end{align*}
 
{\bf 9.1. Reformation of Theorem \ref{thm:A}.}

In Theorem \ref{thm:8.4}, we have
characterized $K(2\alpha, \beta|r)$ with $\ell>0, r>0$
whose Alexander polynomial is monic.
Hence now Theorem \ref{thm:A} is equivalent
to Theorem \ref{thm:5.1} below.

%%Thm 9.1 final. 
\begin{thm}\label{thm:5.1}
Fibredness of $K(r)=K(2\alpha,\beta|r)$ with
$\ell=\ellkB \neq 0$ is determined as follows,
where we assume $\ell>0$ and $r>0$
by Theorem \ref{thm:1}. 
In each case below, 
$a_{i,j}$ and $a'_{i,j}$ are arbitrary.\\
Case 1. $\ell=r=1$.
$K(1)$ is fibred if and only if

(a) $d_i, e_j=\pm 1$, for any $i, j$, and

(b) in each $P_i$ and $Q_i$,
$b_{i,k}$ and $b'_{j,p}$ are $2$.\\
Case 2. $\ell r\ge2$.
$K(r)$ is fibred if and only if

(a) $d_i, e_j=\pm 1$, for any $i, j$, and

(b) in each $P_i$ and $Q_i$,
$b_{i,k}$ and $b'_{j,p}$ are $1$ or $2$.
\end{thm}

\medskip
%\secti{Proof of \Theorems \ref{thm:5.1} (I)}  %15
%
%
{\bf 9.2. Proof of Theorem \ref{thm:5.1}, Case 1.}
 
%In this subsection, we prove \Theorem \ref{thm:5.1} 
%Case 1.  %This case is slightly different from other cases.  
We assume in this subsection
%(9.1)    \hfil
\begin{equation}
\ell = \ellkB = 1\ {\rm and}\ r=1.
\end{equation}

First suppose $K(r)$ is fibred.  
Then $\Delta_{K(r)}(t)$ is monic, and hence, by 
\Proposition \ref{prop:8.3},
%(9.2) \hfil
\begin{align}
&(a)\ d_i, e_j = \pm 1,\ {\rm and}\nonumber\\
&(b)\ b_{i,k}=b_{j,p}'=2\ {\rm  for\ any}\ i,j,k,p.
\end{align}
 
This proves the \lq \lq only if" part of Case 1. 

Conversely, suppose (9.2) is satisfied.  
Rewrite the continued fraction as the
modified continued fraction.
Then some of new $b_{i,k}, b'_{j,p}$ may be
zero, but still (9.2) implies the
deficiency $\rho = 0$.
Now, by Proposition \ref{prop:15.2},
we see that (9.2) also implies the following:
Let $F^*$ be the surface obtained from
$F(r)$ (constructed in Section 7) by
removing all the bands connecting the disks
and cylinders. Then $F^*$ is
obtained from $F(r)$ by deplumbing Hopf bands.
Therefore, to prove that $K(r)$ is fibred,
it suffices to show that $F^*$ 
is a fibre surface.
%%We should remark that $b_{i,k}$ and
%%$b'_{j,k}$ in (9.2)
%%are elements of the standard continued fraction of 
%%$\beta/2\alpha$,
%%while in Section 8, we use the same letter $b_{i,k}$ and $b'_{j,k}$
%%for elements of the {\it modified} continued fraction.
%%Therefore, $b_{i,k}$ and $b'_{j,k}$
%%in (9.2) are not $0$, but $b_{i,k}$ and $b'_{j,k}$ in 
%%Section 8 can be $0$.
%%However, we proved in \Proposition \ref{prop:15.2} that if $b_{i,k}$ or 
%%$b'_{j,k} =0$, then the surface is also obtained from some surface by 
%%plumbing a Hopf band.
%%Therefore, these cases are implicitly included in (9.2).

The following lemma shows that $F^*$ is a 
fibre surface, and hence Theorem \ref{thm:A} Case 1
is proved.
\qed

\begin{lem}\label{lem:16.1} (Braided fibre surface)
Let $L_1$ and $L_2$ be (naturally) oriented closed braids
in a tubular neighbourhood 
$N(L)$ of a Hopf link $L$, where $L_1$ and $L_2$ are 
embedded in different components of $N(L)$. 
Suppose that $L_1$ is a positive closed braid.
Then $L_1$ is a fibred link and a fibre surface $S$ for $L_1$
is obtained by applying Seifert algorithm.
Now replace each component $L_{2,i} (i\le i\le \mu)$ of 
$L_2$ by an annulus
$B_i$ whose core is $K_{2,i}$, 
but the number of twists of $B_i$ is arbitrary.
We assume that $B_i$ intersects $S$ transversally
in ribbon singularities.
By smoothing all the
the ribbon singularities, we obtain a Seifert surface
$F$ for $L_1\cup\partial 
B_1\cup \cdots \cup\partial B_{\mu}$.
Then $F$ is a fibre surface.
\end{lem}

{\it Proof.}  
The surface $S$ consists of disks and 
bands connecting these disks.  
Since $L_1$ is a positive braid, 
each band has only a positive half twist.  
Since neighbouring two bands form a Hopf band, 
we can eliminate one of the bands by deplumbing. 
After all, 
we may assume that $L_1$ is a trivial knot, i.e. 
$S$ consists of $\nu$ disks and $\nu-1$ bands, 
and it suffices to show that $F$ constructed 
using this $S$ is a fibre surface. 
 See \Figure 9.1.
 
 \fif

Denote by $\mu$ the number of strings of the braid of $L_2$.
Now consider a sutured manifold $M=F \times I$.  
For the definition of sutured manifold and its decompositions,
see
\cite[pp.8--10 and Appendix A]{Ga1} and 
\cite[Section 1]{Ga2}.
$M$ is a solid ball with $\nu \times \mu$ holes and $\nu \times \mu$ 
$1$-handles attached.
Applying a series of $C$-product decompositions, first we fill
these holes by $2$-handles, and obtain a  ball 
$M'$ with 
$\nu \times \mu$ 1-handles attached,
where the suture on $M'$ is the equator, and
each $1$-handle has exactly one suture, which is parallel to a co-core.
Since each of the  $1$-handles
connects the north hemisphere
and the south hemisphere of $M'$
without a local knotting,
we can arrange the $1$-handles by sliding their feet
so that they are attached to $M'$ trivially.
Then, by a $C$-product decomposition, we can amalgamate
a pair of $1$-handles, and eventually, we have
a solid torus whose sutures are two meridians.
Applying one more $C$-product decomposition, we have
a ball with a single suture.
Therefore, the original surface $F$ is a fibre
 surface.
\qed

{\bf 9.3. Proof of \Theorem \ref{thm:5.1} Case 2.}

% 
%In this subsection, we prove 
%\Theorem \ref{thm:5.1}, Case 2. So w
In this subsection, we assume that
%(9.3) \hfil
\begin{equation}
 \ell r\ge 2.
 \end{equation}
 
First, we note that \Proposition \ref{prop:8.2} 
proves the \lq \lq only if" part.
Therefore, suppose that
the continued fraction of $\beta/2\alpha$ satisfies
%(9.4) \hfil
\begin{align}
&(a)\  d_i, e_j = \pm1\ {\rm   and}\nonumber\\ 
&(b)\ b_{i,k}\ {\rm and}\ 
b_{j,p}' = 1\ {\rm  or}\ 2\ {\rm  for\ any}\ i,j,k,p.
\end{align}

Again rewrite the continued fraction as a 
modified continued fraction. Then in (9.4) (b)
we have \lq $0, 1$ or $2$', in stead of 
\lq $1$ or $2$'.

By Proposition \ref{prop:15.2},
we can apply compressions
corresponding to all $b_{i,k}, b'_{j,p} =1$, and
remove all bands with $d_i, e_j = \pm1$ 
and $b_{i,k} = 
b_{j,p}' = 0$ or $2$ by deplumbing Hopf bands.
Denote the
resulting surface by $\wti{F}$. 
To complete the proof, it suffices to show
$\wti{F}$ is a fibre surface.
Recall that in the 
proof of Proposition \ref{prop:15.2},
each compression corresponds to cutting the surface
along a properly embedded arc.

%%Then we obtain a surface like $F^*$, 
%%where all the bands are just removed
%%with (9.4) are just removed, and hence
%%\Lemma \ref{lem:16.1} is applicable.
%%Therefore, to prove that $F^{**}(r)$ is a fibre surface, 
%%we first fill the space 
%%that occurs as the result of compression
%%by a suitable Hopf band, 
%%and show that the new surface thus obtained
%%is of the form depicted in \Figure 15.1 (4).
%%Then \Lemma \ref{lem:16.1} shows that $F^*(r)$ 
%%is a fibre surface.
In the following, 
we depict how to undo each of the cuts 
by plumbing a Hopf band.  
To do this, the assumption that 
$\ell r \ge 2$ is essential.
In fact, to undo the cuts, 
we use two consecutive holes 
that occur as intersections of a 
cylinder and the disk(s).

It suffices to consider each of Types II, III and IV in 
\Proposition \ref{prop:15.2}.

 For Type II (resp. III),
we plumb a band $B$ along the  curve depicted in 
\Figure 9.2 (resp. 9.3).  
Then we can undo the cut by
plumbing a Hopf band.

\stnA
\stnB

%\stnAB

Type IV is a bit complicated. Since the argument is similar, we only prove for the case where the band connects
the positive disk and the positive cylinder
(in Figure 7.8.3 (a)). The cut made by compression is
in Figure 8.4.
In \Figure 9.4,
 the band $B'$ gives a compressing disk so that the
 result of compression is as in \Figure 8.4.
 
  \stnC

 To undo the cut, we consider three subcases:
 Let $A$ be the annulus to which $B'$ is connected.
% and $F$ is the original surface for $K_1$ before twisting along $K_2$.
 Subcase (i): There are some disks below $A$ in 
 $\wti{F}$ (\Figure 9,4 (a)).
 Subcase (ii): There are more than one disk
 above $A$ but no disks below $A$ in 
 $\wti{F}$ (\Figure 9.4 (b)), 
 Subcase (iii), There are only one disk above $A$
 and no disks below $A$ in $\wti{F}$
 (\Figure 9.4 (c)). Note that in Subcase (iii),
 $r\ge 2$ by assumption that $\ell r>1$.
 
In each subcase, using the arc depicted in \Figure 9.4,
we can add a band $B$ by
plumbing a Hopf band after compressing at $B'$.
Then, by sliding $B$ along $A$
to the cite of compression,
we can undo the cut. 
This fact is also understood by seeing that the bands 
$B$ and $B'$ cancel each other.

%
% 
% First add a half-twisted band $B'$ as is shown in \Figure 16.4 (1).
% But this addition of a band $B'$ is in fact a plumbing of a Hopf band 
% along the dotted line in  \Figure 16.4 (1) and (2), and this Hopf plumbing
% fills a space that was created in \Figure 14.6.
% 
% %                 \Fi{16.4}
% 
Now all cuts are undone by Hopf plumbings.
Then as in Case 1, we can further deplumb Hopf bands
and apply
\Lemma \ref{lem:16.1}.
 This proves \Theorem \ref{thm:5.1}, Case~2.
\qed

\medskip

\secti{Proof  of \Theorem \ref{thm:D}} %11

In this section, 
we determine the genus of a knot $K(r)$ for 
the case $\ellkB=\ell=0$, and thus prove Theorem
\ref{thm:D}.

In Section 7, we span $K_1$ by a canonical Seifert surface $F_1$.  
By applying Dehn twists on $F_1$ along $K_2$, 
we obtain a Seifert surface $F(r)=F$ for $K(r)$.  
We will show that $F$ is of minimal genus.
Since $F$ and $F_1$ have the same genus, 
we show in fact that
$g(K(r))$ is equal to the number of cylinders in $F_1$.
However, $g(F)$ is much 
larger than one half of 
the degree of the Alexander polynomial of $K(r)$,
c.f., \Proposition \ref{prop:10.2}.
Therefore, in order to show that $F$ is of minimal genus, 
we use geometry.
 First, deplumbing a twisted annulus from $F$ does not
hurt the genus-minimality by the additivity of genus under
the Murasugi sum. 
So we remove all bands connecting
the boundaries of the same cylinder.
Then our main tool is the sutured manifold hierarchies.
%For the definition of sutured manifold and its decompositions,
%see
%\cite[pp.8--10 and Appendix A]{Ga1} and 
%\cite[Section 1]{Ga2}.
 As a special case of general results
of sutured manifold hierarchies, we have the following
(see \cite[Corollary 1.29]{Ga1}):

\begin{prop}\label{prop:std}
Let $(M, \gamma)=(F\times I, \partial F \times I)$
be the sutured manifold obtained from a Seifert surface $F$.
Apply complementary disk- (annulus-) decompositions to
$(M, \gamma)$ and suppose we obtain
$(V, \delta)$ where $V$ is a standard solid torus and
each suture is a loop running longitudinally once and
meridionally non-zero times. 
Then $F$ is of minimal genus.
\end{prop}

Throughout the rest of this section,
we omit the adjective \lq complementary'
for complementary sutured manifold decompositions,
since we only deal with such decompositions and
no confusions are expected.

Let $[2u_1,2v_1,2u_2,2v_2,\ldots, 2u_m,2v_m, 2u_{m+1}]$ 
be the continued fraction of $\beta/2\alpha$.   
Suppose that $\ellkB=0$.  
Then to prove Theorem \ref{thm:D}, we show
the following:
%(10.1) \hfil
\begin{equation}
g(K(r))= \dfrac{1}{2}\sum_{i=1}^{m+1}|u_i|
=\dfrac{\lambda}{2}
(= \#\{{\rm cylinders\ of\ } F_1\}).
\end{equation}

 Note that since $\ellkB=0, \sum_{i=1}^{m+1}|u_i|=
\lambda=\#\{{\rm edges\ of\ } G(S)\}$
is even.
 
 {\it Proof.}
 We prove (10.1) by induction on $\lambda$.
 If $\lambda=2$, then (10.1) is obvious since
 $F(r)$ is a plumbing of two twisted annuli.
 Suppose $\lambda\ge4$.
 First we deplumb all bands corresponding to
 proper local maximal, or minimal vertices,
 i.e., those connecting the two boundaries of a cylinder.
 Denote by $\hat{F}$ the resulting surface.
 We inductively reduce the graph $G(S)$ and accordingly
amalgamate the solid tori in 
$(\hat{F}\times I, \partial \hat{F}\times I)$
by disk- (annulus-) decompositions,
until we have only one torus where each of the sutures
run longitudinally once and
meridionally non-zero times. 
After that, we will see that all such deplumbing and
amalgamations commute
with Dehn twists along $K_2$, and hence
by \Proposition \ref{prop:std}, we have (10.1).

\elebb %11.2

Case 1: There is a vertex in $G(S)$ incident to
two consecutive rising edges and two consecutive
falling edges, as in \Figure 10.1 (a1) or (b1), 
where the white 
vertices may be a terminal vertex
or a non-terminal vertex, and they may lie on the $x$-axes.

 As in \Figure 10.2, we amalgamate the top solid torus
with the second top one.
 There are several cases according to the
 number of  twists in the two bands connected.
 
  \eleb
  
  Case 1.1. The two bands are twisted in the opposite directions:
 See \Figure 10.2 (a). First apply a 
disk decomposition
using the disk with shadow, then apply a product
disk decomposition. Note that if the two bands
are both only half-twisted, then the first disk decomposition
is also a product disk decomposition.

 Case 1.2. Two bands are twisted in the same direction:
 See \Figure 10.2 (b) and (c).
As before, we apply a disk decomposition and a product disk
decomposition. \Figure 10.2 (b) depicts the case where both bands
are more than half-twisted. In this case, we have two extra
sutures, but they do not affect the following inductive amalgamations.
\Figure 10.2 (c) depicts the case at least one of the bands is
half-twisted. Note that if both bands are half-twisted,
then the first disk decomposition is a product disk decomposition.

Case 2: There are no subgraph considered in Case 1, or all such subgraphs
are removed.
Now it suffices to find a subgraph as in \Figure 10.3 or 10.4 and
amalgamate a solid torus.

\elecc

\eledd

Since the other cases are similar, we only consider subgraphs in
\Figure 10.3 (a1) and \Figure 10.4.

%it suffices by induction to show that we can amalgamate solid tori
%according to the reduction as depicted in \Figure 11.4 or 11.5.
By construction,
a rising edge above (resp. below) 
the $x$-axes is paired with the
a falling edge on its right (resp. left) so that
the pair corresponds to a cylinder in $F_1$.
Sutured manifold decompositions amalgamate
the solid tori as respectively depicted in \Figures 10.5 and 10.6.
In \Figure  10.5, corresponding to \Figure 10.3 (a1),
two cylinders showing the same side are
connected by an even-twisted band.
In \Figure 10.6, corresponding to \Figure 10.4,
two cylinders showing the opposite sides are connected
by an odd-twisted band.
However, as seen in \Figure 10.2, longitudinal sutures may
have accumulated.
In \Figures 10.5 and 10.6, we first
apply an annulus decomposition and then a disk decomposition.
Now we have amalgamated all the tori into one and see
that we may apply the Dehn twist along $K_2$ beforehand.
Therefore, by \Proposition \ref{prop:std},
(10.1) is obtained.
\qed

The proof of \Theorem \ref{thm:D} is now completed.

\eleee
\medskip

\eleff

\medskip
As an immediate consequence of \Theorem \ref{thm:D}, 
we have:

\begin{cor}\label{cor:11.last}
Suppose $\ellkB=0$ and $\alpha\ge2$.
Then, for any $r>0$, $K(2\alpha, \beta|r)$
is never unknotted.
\end{cor}

\medskip
\medskip

\secti{Proof of \Theorem \ref{thm:B}} %B (I)} %sec12

{\it Proof of \Theorem \ref{thm:B} (a). }%B(2)(a).}
%%%Suppose $K_0$ is unknotted.  
%%%Then, it follows from \Corollary \ref{cor:11.last}
%%%that 
%%%we may assume that $K(2\alpha,\beta|r; K_0)$
%%%is knotted for $r \ge 1$ and $\alpha \ge 2$, and 
%%%$\ellkB=0$.
Write
   $\Delta_{B(2\alpha,\beta)}(x,y) = (x-1)(y-1)f(x,y)$.
Then by \Proposition \ref{prop:6.1}, we can write
$\Delta_{K(r)}(t) = r(t-1)^2 f(t,1) +\es t^k$,
where $\es = \pm 1$ and $k$ is chosen so that 
$\Delta_{K(r)}(t)$ is symmetric.
Suppose $r\ge 2$.
If $f(t,1) \neq 0$, then $\Delta_{K(r)}(t) $
is not monic and hence $K(r)$ is not 
fibred.  
Suppose $f(t,1)=0$.  Then $\Delta_{K(r)}(t)=1$.  
However, since $\alpha \ge 2$,  $K(2\alpha,\beta|r)$ 
is 
not trivial by Corollary \ref{cor:11.last}, 
and therefore $K(r)$ is not 
fibred for 
$r \ge 2$.
\qed

\medskip
The rest of this section is devoted to the
proof of Theorem \ref{thm:B} (b).
As we remarked in subsection 3.4, 
we may assume $r>0$ and $\beta>0$.
To prove  \Theorem \ref{thm:B} (b), 
we need the following two propositions.\\
Let $\{P_1,d_1,Q_1,e_1,P_2,d_2,Q_2,e_2,
  \ldots, P_m,d_m,Q_m\}$ be the
canonical decomposition of $\beta/2\alpha$.
Using \Theorem \ref{thm:D}, first we prove the following.
 
\begin{prop}\label{prop:12.1}
 
Suppose $\ellkB=0$.  
Then for any $r\ge 1,
        g\bigl(K(r)\bigr)=\frac{1}{2}deg \Delta_{K(r)}(t)$
 if and only if $m=1$, i.e. $\{P_1,d_1,Q_1\}$ is the 
canonical decomposition of  $\beta/2\alpha$.
\end{prop}

{\it Proof.}  By \Theorem \ref{thm:D}, for any $r\ge1$,
${\displaystyle
2g\bigl(K(r)\bigr)=   \sum_{i=1}^{m}\Bigl\{
\sum_{k=1}^{s_i+1}|a_{i,k}| +\sum_{k=1}^{q_i +1}|a'_{i,k}|\Bigr\}
}$.\\
On the other hand,
       $\deg \Delta_{K(r)}(t) \le 2 max\{h, |q|\}$,
where $h$ (and $q$) is the $y$-coordinate of
the absolute maximal 
(and minimal) vertices in $G(S)$ 
(Proposition \ref{prop:10.2}).
Therefore, if there exist local minimal vertices (not end vertices), 
we see easily that
   $2g\bigl(K(r)\bigr) > 2 \max\{h, |q|\} \ge \deg 
   \Delta_{K(r)}(t)$,
and hence, there is only one local (and hence, absolute) 
maximal vertex, and therefore, $m=1$.
\qed
 
\begin{prop}\label{prop:12.2}
Suppose $\ellkB=0$.   
Suppose further, $\{P_1,d_1,Q_1\}$ is the canonical 
decomposition of  $\beta/2\alpha$.  Then
$\Delta_{K(1)}(t)$ is monic if and only if 
   $d_1 = \pm1$.
\end{prop}

{\it Proof.}  
This follows from \Corollary \ref{cor:10.3}.
\qed

 \medskip
 
Now we proceed to the proof of \Theorem \ref{thm:B}(b).

 {\it Proof of the \lq\lq if\rq\rq\ part.}
Suppose 
that
the modified continued fraction of $\beta/2\alpha$
is of the form
$S = [[1,b_1,1,b_2,
\ldots,1,b_{p-1},1,d_{1},-1,-b'_{p-1},-1,\ldots, -b'_{1},
-1]]$,
where $b_i$ and $b_i'$, $1\le i \le p-1$, are $0$ or $1$ 
and 
$d_1 = \pm1$.
What is to show is that $K(1)=K(2\alpha, \beta|1)$ is fibred.
Let $F_1$ be the canonical Seifert surface for
$K_1$. (See \Figure 7.6 (a).)
By twisting $F_1$ once by $K_2$, we obtain a
Seifert surface for $K$.
Since $d_1=\pm 1$, the band corresponding to the maximal
vertex, is a Hopf band, and hence, we may remove it
by deplumbing. Denote the resulting surface by $\hat{F}$.
Now, since $b_i$ and $b_i' (1\le i\le p-1)$ are either $0$ or $1$,
every band in $\hat{F}$ is only half-twisted.
Therefore, every disk decompositions employed in the proof of
(10.1) Case 1
is in fact a product disk decomposition,
and hence $\hat{F}$ is a fibre surface by \cite[Theorem 1.9]{Ga2},
and $K(1)$ is a fibred knot.

{\it Proof of the \lq\lq only if\rq\rq\ part.}
Suppose $K(1)=K(2\alpha,\beta|1)$ is a 
fibred knot.  
Then by 
\Proposition \ref{prop:12.1}, the continued fraction $S$ 
for $\beta/2\alpha$ must be 
$\{P_1,d_1,Q_1\}$.  
Therefore, the modified continued 
fraction is of the form:

 $\beta/2\alpha=[[1,b_1,1,b_2, \ldots, 
 1,b_{p-1},1,d_1,-1,-b_{p-1}',-1, \ldots, -1,-b'_1,-1]]$.
 
Consider the canonical Seifert surface $F_1$ for $K_1$,
as in \Figure 7.6 (a).
A Seifert surface $F$ for $K(1)$, consisting of $p$ cylinders
$A_1, A_2, \ldots, A_p$ and $2p-1$ bands,
 is obtained from $F_1$
 by applying  a Dehn twist once along $K_2$.  
 In Section 10, we have seen that $F$ is of minimal genus,
 and hence a fiber surface for $K(1)$.
 By \Proposition \ref{prop:12.2}, $d_1 = \pm1$
 and hence we remove the band on the top annulus $A_p$
 by deplumbing a Hopf band, and denote by
 $\hat{F}$ the resulting surface.
Then  the 
 inclusion map below must be onto:
 % (11.1)    \hfil
%\begin{equation}
\hfill
$\phi : \pi_1(\hat{F})  \longrightarrow \pi_1(S^3 - 
 \hat{F})$. \hfill (11.1)\\
%\end{equation}
\setcounter{equation}{1}

Now, we show that
all $b_j$ and $b_j'$ are $0$ or $1$ and hence that
if at least one of $b_j$ and
$b_j'$ is neither $0$ nor $1$,
then $\phi$ is not onto and hence $K(1)$
 is not fibred.

If $b_1$ and $b_1'$ are $0$ or $1$,
then we can \lq remove\rq\ the bottom annulus
$A_1$ by two product decompositions
as in the proof of (10.1) Case 1.
Therefore by \cite[Lemma 2.2]{Ga2}, we may 
assume without loss of generality that 
$b_1'\neq 0,1$, i.e., the band $B'_1$ is more than
half-twisted.
To show that $\phi$ is not onto, we 
need explicit presentations of the groups 
$\pi_1(\hat{F})$ and $\pi_1(S^3-\hat{F})$.
To do that, we deform $\hat{F}$ as depicted in
\Figure 11.1 (where $p=3$).
Note that both $\pi_1(\hat{F})$ and 
$\pi_1(S^3-\hat{F})$ are free of rank $2p-1$.

 \Figure 11.1 also depicts
 a base point $**$ and the generators of
 $\pi_1(S^3-\hat{F})$ 
denoted by $x_1, x_2, \ldots,x_p$ and 
$a_1,a_2,\ldots,a_{p-1}$.
Take a base point $*$ for $\pi_1(\hat{F})$ on 
$A_1$ as in \Figure 11.1.
 The generators for 
$\pi_1(\hat{F})$ are denoted by 
$\alpha_1, \alpha_2, \ldots, \alpha_p$ and 
$\beta_1, \beta_2,\ldots, \beta_{p-1}$.
A loop $\alpha_i$ starts at $*$ moving toward on $A_i$ 
through bands
$B_1,\ldots,B_{i-1}$ and circle once around $A_i$ 
counter-clockwise, 
and then return to $*$ through $B_{i-1}, \ldots, B_1$.  
A loop $\beta_i$ starts at $*$ moving toward $A_{i+1}$ through 
$B_1,B_2, \ldots, B_i$ 
and returns to $*$ passing through first $B_i'$ and then 
$B_{i-1}, B_{i-2}, \ldots, B_1$.  

%\Fi{13.1}
\thir

We must express $\phi(\alpha_i), \phi(\beta_j)$ in terms of 
$x_i , a_j$.   
Let $D=x_1 x_2 \ldots x_p$.
For simplicity, we use $\alpha_i$, (or $\beta_j$) instead of  $\phi(\alpha_i)$ 
(or $\phi(\beta_j)$).
Then we have the following:% \lq \lq equations":
% (11.2)
\begin{align}
&\alpha_1= a_1D,  \nonumber\\
&\alpha_2= u_1a_2D u_1^{-1}, \nonumber\\
&\alpha_3= u_1u_2a_3 D u_2^{-1} u_1^{-1},  \nonumber\\
&\indent \vdots\nonumber\\
& \alpha_{p-1}= u_1 u_2\ldots u_{p-2}a_{p-1}D u_{p-2}^{-1}
 \cdots u_1^{-1},  \nonumber\\
& \alpha_p= u_1u_2\cdots u_{p-1} D u_{p-1}^{-1}\ldots u_1^{-1}.\\
%&\nonumber\\
%(11.3)
&\beta_1 = u_1 w_1,\nonumber\\
&\beta_2 = u_1(u_2 w_2) u_1^{-1},\nonumber\\
&\beta_3 = u_1u_2(u_3 w_3) u_2^{-1} u_1^{-1},\nonumber\\
& \indent
 \vdots \nonumber\\
& \beta_{p-1} = u_1u_2\ldots u_{p-2}(u_{p-1} w_{p-1}) 
 u_{p-2}^{-1}\ldots u_1^{-1},
\end{align}
where $u_i=a_{i}^{b_{i}}$
and $w_i = x_{i+1}^{-1}(x_i^{-1} \ldots x_1^{-1} D x_1 
\ldots x_i D^{-1} a_i^{-1} )^{b_i'},  1 \le i\le p-1$.

Denote $h_i=x_i^{-1} \ldots x_1^{-1} D x_1 
\ldots x_i D^{-1} a_i^{-1}$.
By assumption, $b'_1 \neq 0,1$.
Let $H = {\rm Im} \phi$ and $G= \pi_1(S^3-\hat{F})$. 
We want to show that $H \neq G$.
First we may suppose that $a_1 \in H$ and  $x_1 x_2 \in H$.  
Otherwise, obviously, $H \neq G$, and we are done.
Now since $a_1\in H$, by (11.2), we have $D \in H$, and hence, 
$a_2 \in H$, since $H \ni\alpha_2=u_1a_2D u_1^{-1}$, and 
$H \ni    u_1=a_1^{b_1}$ and $H \ni D$.
Repeat this process to obtain
%(11.4) \hfil
\begin{equation}
a_1, a_2, \ldots, a_{p-1} \in H\ {\rm  and}\ D\in H.
\end{equation}
Therefore, $H$ is generated by
%(11.5)  \hfil
\begin{equation}
\{ a_1, a_2, \ldots, a_{p-1}, D, x_1x_2, 
   \alpha_1, \ldots, \alpha_p, 
   \beta_1,\ldots, \beta_{p-1}\}.
\end{equation}   
However, since $\alpha_i$ is written in terms of $a_i$ and $D$, 
we can eliminate these $\alpha_i$ from the set of generators (11.5), 
and hence $H$ is generated by
 %(11.6)  \hfil
\begin{equation}
\{ a_1, a_2, \ldots, a_{p-1}, D, x_1x_2, \beta_1,\ldots, \beta_{p-1}\}.
\end{equation}
Since $u_i \in H, 1\le i\le p-1$, we introduce new generators 
$\gamma_j$, replacing $\beta_j$, as
%(11.7) \hfil
\begin{equation}
\gamma_1 = w_1,\gamma_2 = w_2, \ldots, \gamma_{p-1} = w_{p-1}.
\end{equation} 
Then $H$ is generated by
%(11.8) \hfil
\begin{equation}
\{ a_1, a_2, \ldots, a_{p-1}, D, x_1x_2, 
\gamma_1,\ldots, \gamma_{p-1}\}.
 \end{equation}
Since $h_2^{b_2'} = 
(x_2^{-1} x_1^{-1} D x_1 x_2 D^{-1} a_2^{-1})^{b_2'} \in H$ 
and
$H\ni\gamma_2=x_3^{-1} h_2^{b_2'}$, it follows that 
${x_3}^{-1} \in H$.  
%By the similar reasoning, 
Similarly,
using $\gamma_3, \ldots, \gamma_{p-1}$, 
we have ${x_4}^{-1},\ldots, {x_p}^{-1} \in H$.
Therefore, we can replace $x_{i+1}$ by $\gamma_i, i\ge2$, 
and we see that $H$ is generated by
%(11.9) \hfil 
\begin{equation}
\{ a_1, a_2, \ldots, a_{p-1}, D, x_1x_2, x_3, \ldots, x_p, \gamma_1\}.
\end{equation} 
Since $D=x_1x_2\cdots x_p$, we can eliminate $D$ from the set of generators of $H$, 
and $H$ is generated by $2p-1$ elements
%(11.10)   \hfil
\begin{equation}
\{ a_1, a_2, \ldots, a_{p-1}, x_1x_2, x_3, \ldots, x_p, \gamma_1\}.
\end{equation} 
Since $H$ is free of rank $2p-1$, the above set is a set of free generators of $H$.\\ 
On the other hand, $G$ is freely generated by
$\{ a_1, a_2, \ldots, a_{p-1}, x_1,x_2, x_3, \ldots, x_p\}$.
 
Now introduce a new generator $z_2 = x_1x_2$ and 
replace $x_2$ by $z_2$. Then we have:
%\begin{tabular}{ll}(11.11) &
\begin{align}
&(1)\  G \ {\rm is\ generated\ (freely)\ by}\
\{ a_1, a_2, \ldots, a_{p-1}, x_1,z_2, x_3, \ldots, x_p\},\ {\rm  and}\nonumber\\
&(2)\  H\  {\rm is\ generated\ (freely)\ by}\
\{ a_1, a_2, \ldots, a_{p-1}, z_2, x_3, \ldots, x_p, \gamma_1\},\nonumber\\
&\ \ \ \ {\rm where} \gamma_1=w_1=x_2^{-1}(x_1^{-1} D x_1D^{-1} a_1^{-1})^{b_1'}.
\end{align}
 
Therefore, if $H=G$, then $x_1\in H$. 
In other words, $x_1$ can be written as a word on 
$a_i, 1\le i \le p-1,z_2,x_j, 3\le j\le p$, and $\gamma_1$.  
(Note that $H$ is a free group of rank $2p-1$.)
 
Recall  $b'_1\neq 0,1$.

Case 1: $ b'_1 > 0$,   i.e. $b'_1\ge2$.  
Then,
\begin{align*}
\gamma_1 &=
x_2^{-1}(x_1^{-1} D x_1D^{-1} a_1^{-1})(x_1^{-1} D x_1D^{-1}
a_1^{-1})^{b'_1-1}\\
&=
z_2^{-1} D x_1D^{-1} a_1^{-1}(x_1^{-1} D x_1D^{-1} 
a_1^{-1})^{b'_1-1}.
\end{align*}
Since $z_2^{-1}$, $D$ and $D^{-1}a_1^{-1}$ are in $H$, 
we can replace $\gamma_1$ by
\begin{equation*}
\delta_1= x_1(D^{-1} a_1^{-1})(x_1^{-1} D x_1D^{-1} 
 a_1^{-1})^{b'_1-2}(x_1^{-1} D x_1).
 \end{equation*}
Since  $D=z_2 x_3 \cdots x_p$, $D$ does not involve $x_1$ and hence 
$\delta_1$ is of the form:\\
 $\delta_1= x_1 W_1(x_1^{-1} W_2 x_1W_1)^{b'_1-2}(x_1^{-1} W_2 x_1)$,
where $W_1 =D^{-1} a_1^{-1}$ and $W_2 = D$, 
none of which involves $x_1$.
Therefore, $\delta_1$ is a reduced word.  
Since $b'_1-2\ge0, \delta_1$ involves $x_1$ 
at least three times, and $\delta_1$ is of the form:
$x_1 U x_1$.  
Therefore, we cannot write $x_1$ in terms of $a_i, 1
\le i\le p-1, z_2,x_3,\ldots,x_p$ and $\delta_1$.
 
Case 2:  $b'_1<0$.  Write $b'_1=-q, q\ge1$.  Then

 $\gamma_1=x_2^{-1}(a_1 D x_1^{-1}D^{-1} x_1)^q
=z_2^{-1}x_1(a_1 D x_1^{-1}D^{-1} x_1)^q$.
Again, since $z_2 \in H$, we can replace $\gamma_1$ by $\delta_1'$,
$\delta_1'=x_1(a_1 D x_1^{-1}D^{-1} x_1)^q
=x_1(W_1^{-1}x_1^{-1} W_2^{-1} x_1)^q$.

Since $q>0, \delta_1' $is a reduced word and $\delta_1'$
 is of the form $x_1Vx_1$, and $V$ contains $x_1$ at least once.  
 Therefore $x_1$ cannot be written in terms of 
 $\delta_1', a_i,
 1\le i\le p-1, z_2, x_3, \ldots, x_p$.
It proves that $H\neq G$, and hence $\phi$ is not onto. 
\Theorem \ref{thm:B} is now
 completely proved.
\qed

\section{Examples}
In this section, we discuss several examples
that illustrate our main theorems.

\begin{ex}\label{ex:lzero}
All $2$-bridge links in this example have linking number
$0$. 
\Theorem \ref{thm:D} determines 
$g(K(2\alpha,\beta|r)$
and the fibredness is checked by  \Theorem 
\ref{thm:B}.

(1) Consider $K(48,31|r)$.
Since $31/48=[2,2,-4,2,2]$,
the genus is $2$, and it is not fibred for any $r>0$.
The graph is given in \Figure 12.1 (1).
The lifts $\{\wti{K}_j\}$ of $K_1$ in the infinite cyclic cover
$M^3$ of $S^3\setminus K_2$ are depicted in 
\Figure 12.2 (1).
From it, we see that $\Delta_{K(48,31|1)}(t)$ is not monic.
This can also be checked by evaluating 
$\Delta_{B(48,31)}(x,y)$:\\
$\Delta_{B(48,31)}(x,y)=(1-x)(1-y)\{
1-(x+y)-xy(x+y)+x^2y^2\}$,
and hence, $\Delta_{K(48,31|r)}(t)=2r(1-t)^2+t$.
Thus, we see that
$\Delta_{K(48,31|1)}(t)$ is not monic.

(2) Consider $K(64,41|r)$.
Since $41/64=[2,2,-4,-2,2]$, the genus is $2$,
and it is not fibred for any $r>0$.
See \Figure 12.1 (2) for its graph.
The lifts $\{\wti{K}_j\}$ of $K_1$ in $M^3$
is depicted in \Figure 12.2 (2).
From it, we see that $\Delta_{K(64,41|r)}(t)=1$.
The same result is also obtained using 
$\Delta_{B(64,41)}(x,y)$.
$\Delta_{B(64,41)}(x,y)=(1-x)^2(1-y)^2(1+xy)$,
and hence $\Delta_{K(64,41|r)}(t)=1$ 
for any $r>0$.

\etna

\etnb

(3) Consider $K(40,11|r)$.
Since $11/40=[4,2,-2,-2,-2]$,
the genus is $2$ and it is fibred only for $r=1$.
See \Figure 12.1(3) for its graph.
The lifts are depicted in \Figure 12.2~(3).
Also, we have:
$\Delta_{B(40,11)}(x,y)=(1-x)(1-y)\{(x+y)-xy+
xy(x+y)\}$ and hence,
$\Delta_{K(40,11|r)}(t)=r(1-t-t^3+t^4)+t^2$.

(4) Consider $K(112,71|r)$.
Since $71/112=[2,2,-2,2,2,2,-2]$,
the genus is $2$ and it is not fibred for any $r>0$.
The lifts $\{\wti{K}_j\}$ of $K_1$ are depicted in \Figure 12.2 (4).
Also, we have:
$\Delta_{B(112,71)}(x,y)=
(1-x)(1-y)\{1-2(x+y)+(x+y)^2-2xy(x+y)+x^2y^2\}$,
and hence, $\Delta_{K(112,71|r)}(t)=2r(1-t)^2+t$.
\end{ex}

\begin{ex} %18.2
Each of the first two $2$-bridge links has linking 
number $1$, and the other two links have
linking number $2$.
The graphs are depicted in \Figure 12.3.
We use \Theorem \ref{thm:A},
\Proposition \ref{prop:8.1} and
\Theorem \ref{thm:8.4}.\\
(1) Consider $K(18,13|r)$.
Since $13/18=[2,2,2,-2,-2]$, we see
$P_1=[[1,1,1]], Q_1=[[-1]]$ and $d_1=-1$.
Since $b_{11}=1$, it follows from \Theorem \ref{thm:8.4}
that $\Delta_{K(18,13|r)}(t)$ is not monic for $r=1$,
but it is monic for $r>1$,
and hence $K(18,13|r)$ is fibred for $r>1$.
Since $\lambda=3$ and $\rho=1$, the degree of
$\Delta_{K(18,13|r)}(t)$ is $2r$, and hence,
the genus is $r$ by \Proposition \ref{prop:8.1}
and \Theorem \ref{thm:A}.
These facts are also confirmed
by evaluating
$\Delta_{B(18,13)}(x,y)$:
$\Delta_{B(18,13)}(x,y)=(x+y)-(x^2+3xy+y^2)
+xy(x+y)$,
and use \Proposition \ref{prop:6.1} (1).\\
(2)
 Consider $K(482,381|r)$.  
Since $381/482=[2,2,2,2,-4,-2,-2,2,2,2,2]$, 
we see that $P_1=[[1,1,1]] , Q_1=[[-2,-1,-1]], 
P_2=[[1,1,1]], d_1=1$ and $e_1=1$.
Therefore, we see that $\lambda = 7$ and  
$\rho =3$, and hence, by 
\Proposition \ref{prop:8.1}, the degree of  
$\Delta_{K(482,381)}(t) = 6r$ and the genus is $3r$. 
Further, it follows from \Theorem \ref{thm:8.4}(1) 
that it is not fibred for $r=1$, but it is fibred for $r >1$.   
These facts are also confirmed by evaluating the 
Alexander polynomial of $B(482,381)$, and
using \Proposition \ref{prop:6.1} (1).
$\Delta_{B(482,381)}(x,y) =(-x^3 + 2x^2 -x)y^6 + 
(-3x^4 + 8x^3 -8x^2 +4x -1)y^5 +
 (-3x^5 + 12x^4 -17x^3 +14x^2 -8x +2)y^4 + 
 (-x^6 + 8x^5 -17x^4 +21x^3 -17x^2 +8x -1)y^3 +
 (2x^6 - 8x^5 + 14x^4 -17x^3 +12x^2-3x)y^2 + 
 (- x^6 + 4x^5 - 8x^4 +8x^3 -3x^2)y -x^5 + 2x^4 -x^3$.\\
(3) Consider $K(60,47|r)$.
Note that $\ellk B(60,47)=2$.
Since $47/60=[2,2,2,2,$ $-2,-2,2]$,
we see that the genus is $3r$ and
$\Delta_{K(60,47|r)}(t)$ is monic
and hence is fibred for any $r>0$.
Note that $\Delta_{B(60,47)}(x,y)=
(x+y)-(2x^2+3xy+2y^2)+
(x^3+5x^2y+5xy^2+y^3)-xy(2x^2+3xy+2y^2)
+x^2y^2(x+y)$.\\

(4)
Consider $K(1732,-671)$.
Since $-671/1732=[-2,2,4,2,-2,2,2,4,2]$,
we see by \Theorem \ref{thm:8.4} (2) that 
$\Delta_{K(r)}(t)$ is monic for $r>0$
and hence $K(r)$ is fibred for any $r>0$.
Also, since $\lambda=6$ and $\rho=0$, we have
$g(K(r))=5r+2$.
Note that
$\Delta_{K(1732,-671)}(x,y)=
(x^5-5x^4+9x^3-7x^2+2x)y^5
-(5x^5-23x^4+44x^3-42x^2+18x-2)y^4+
(9x^5-44x^4+87x^3-86x^2+42x-7)y^3-
(7x^5-42x^4+86x^3-87x^2+44x-9)y^2+
(2x^5-18x^4+42x^3-44x^2+23x-5)y+
2x^4-7x^3+9x^2-5x+1$.

%(4) Consider $K(76,21|r)$.
%Since $21/76=[4,2,-2,-2,-2]$, the genus is $3r+1$
%and $\Delta_{K(76,21|r)}(t)$
%is monic for $r>0$ and hence it is fibred for $r>0$.
%Note that $\Delta_{B(76,221)}(x,y)=
%1-2(x+y)+xy(x^2+7xy+y^2)-
%2x^2y^2(x+y)+x^3y^3$.
\end{ex}

\etnc

\section{Fibred satellite knots of tunnel number one}
In this section, 
we determine the genera of the satellite knots of 
tunnel number one, and also solve the question of 
when it is fibred.

Let $\widehat{K}$ be a satellite knot of tunnel number one.  
According to \cite{MS}, 
the companion of $\widehat{K}$ is a torus knot of type $(a, b)$,  
say, and $|a|, |b| >1$, 
and the pattern is the torti-rational knot 
$K(2\alpha, \beta|ab)$ for $|\alpha|>1$.  
To be more precise, 
let $B(2\alpha,\beta)=K_1 \cup K_2$ be a $2$-bridge link.  
Then by construction, a torti-rational knot, 
$K(2\alpha,\beta|r)$ is a knot in 
the interior of a (unknotted) solid torus $V$, 
where $V$ is homeomorphic to 
$S^3 \setminus N(K_2)$, 
$N(K_2)$ being a tubular neighbourhood of $K_2$.  
Let $m$ be a meridian of $V$. 
Then the pattern  is a pair $(V, K(2\alpha, \beta|r))$.  
Generally, if the pattern is $(V, K(2\alpha, \beta|r))$, 
then our technique can be applied on any satellite knot 
with fibred companion $K_C$.  
Therefore, 
in this section we can prove slightly more general 
theorems as follows.

\begin{thm} \label{thm:17.1}
Let $K_C$ be a non-trivial fibred knot and 
$\widehat{K}$ be the satellite knot 
with companion $K_C$ and the pattern 
$(V, K(2\alpha, \beta|r)), r \ne 0$. 
Suppose $\ell = lk(K_1,K_2) \ne 0$. 
Then the following hold:\\
(1)	the genus of $\widehat{K}$ is exactly half of the 
degree of the Alexander polynomial of $\widehat{K}$.\\
(2)	$\widehat{K}$ is fibred 
if (and only if) the Alexander polynomial of 
$\widehat{K}$ (and hence, that of 
$K(2\alpha, \beta|r)$ is monic.
\end{thm}

\begin{thm}\label{thm:17.2}
Under the same notation of Theorem \ref{thm:17.1}, 
suppose $r \ne 0$ and  $\ell =0$. Then
(1)	$\widehat{K}$ is not fibred for any $r\ne 0$.
\cite[Proposition 4.15]{BZ} 
(2)	Let $[2c_1,2c_2, \ldots, 2c_m]$ be the continued fraction of 
$\beta/2\alpha$. Then $g(\widehat{K}) = \sum_{odd\ j} |c_j|$.
\end{thm}

{\it Proof of Theorem \ref{thm:17.1}.}
(1) Let $\phi$ be a faithful homeomorphism of a solid torus 
$V$ in which 
$K(2\alpha,\beta|r)$ is embedded to a tubular
neighbourhood 
$N(K_C)$ of $K_C$ in $S^3$. 
The minimal genus Seifert surface $F$ for $K=K(2\alpha, \beta|r)$
we had in Section 9
%the proof of 
%\Proposition \ref{prop:15.2} 
intersects $\partial V$ at 
$\ell$ parallel longitudes, 
where $\ell=\ell k(K_1,K_2)>0$.  
Since the image of each longitude under $\phi$
spans a fiber surface $S_C$ of genus 
$g(K_C)$ in $S^3-\phi(V)$, we can construct
a Seifert surface $\widehat{F}$ for $\widehat{K}$ by
capping off $\phi(V\cap F)$ by $\ell$ copies of $S_C$.
Hence we have
$g(\widehat{K})\le
        g(K) + 
        \ell g(K_C)$.
Combining with Schubert's inequality
\cite[Proposition 2.10]{BZ}, we have:
%(13.1)  \hfil
\begin{equation}
g(\widehat{K})=
g(K) + \ell g(K_C).
\end{equation} 
However, since $K_C$ is a fibred knot, we see:
%(13.2)  \hfil
\begin{equation}
g(K_C) = \frac{1}{2}deg \Delta_{K_C}(t).
\end{equation} 
With the fact that $g(K)= 
\frac{1}{2} \deg 
\Delta_{K}(t)$ and 
Seifert's theorem \cite[Proposition 8.23 (b)]{BZ}, 
we obtain
%(13.3)  \hfil 
\begin{equation}
2g(\widehat{K})=
deg \Delta_{K}(t) + 
\ell  \deg \Delta_{K_C}(t) =  
 \deg \Delta_{\widehat{K}}(t).
\end{equation}
 
This proves (1).

Next, we prove (2).
Suppose $\Delta_{K(2\alpha,\beta|r)}(t)$ is monic.
Then construct a fibre surface $F$ for 
$K=K(2\alpha,\beta|r)$ as in Section 9. As we did in the proof of
(1), construct a minimal genus Seifert surface $\widehat{F}$
for $\widehat{K}$.
If we needed compressions to obtain $F$ from
the surface $F(r)$, then 
we undo the cuts in $\widehat{F}$
caused by compressions by plumbing Hopf bands. 
Note this is possible since
the plumbings of Hope band in the proof of Theorem \ref{thm:5.1}
can be locally done in $V$. 
Denote by $F'$ the resulting Seifert surface.
To complete the proof, it suffices to show that
$F'$ is a fibre surface.
Now $F'$ looks like as in Figure 7.7, where
horizontal parts are understood as parallel copies of
$S_C$. We can remove the bands by deplumbing Hopf bands
as we did in the proof of Proposition \ref{prop:15.2},
until we have a new Seifert surface 
$F''$ \lq consisting\rq\
of the annuli, and $\ell$ copies of $S_C$ and $\ell-1$ half-twisted band, 
where each pair of adjacent copies of $S_C$ is connected by
a half-twisted band. 
We show that $F''$ is a fibre surface by using
sutured manifolds.
Let $(F''\times I, \partial F'' \times I)$ be the sutured manifold
obtained from $F''$.
 Apply a $C$-product disk decomposition corresponding
to each site of ribbon holes arising from the intersection 
of the annuli and copies of $S_C$. 
Then we have $1$-handles each with a meridional suture. 
Actually, we have $\ell(\lambda-\ell)/2$ $1$-handles.
See Figure 13.1 (a).
 We can remove, by $C$-product decompositions, 
all such $1$-handles coming from the annuli (Figure 13.1 (b)).
Apply a $C$-product decomposition 
between a pair of
parallel copies of $S_C$'s. Then the complementary sutured manifold
is separated into two pieces: one of which is a product sutured manifold
between two copies of $S_C$ and hence we can disregard it
(Figure 13.1 (c)).
Repeating this, we have a sutured manifold obtained from $S_C$ 
(Figure 13.1 (d)). 
Since $S_C$ is a fibre surface, the complementary sutured manifold
is a product sutured manifold. Therefore, 
$F''$ is a fibre surface.
Theorem \ref{thm:17.1} (2) is now proved.
\qed

\satel

\medskip
{\it Proof of Theorem \ref{thm:17.2}.}
For any $r>0$ and any fibred companion $K_C$, we see
$g(\widehat{K})\ge g(K(2\alpha, \beta|r))$.
However, since $\ell=0$, from the construction
of $F_1$ in Section 7, we have $g(\widehat{K})=
g(K(2\alpha, \beta|r))$, and hence (2) follows immediately
from \Theorem \ref{thm:D}.
\qed

%Remark 17.3
\begin{rem}\label{rem:17.3}
When $r = 0$, $K(2\alpha, \beta|0)$ is a trivial knot.  
If we consider the satellite knot $\widehat{K}$
with $K(2\alpha,\beta|0)$ as a pattern knot, 
this satellite knot gives rise to a very interesting problem.  
%As is pointed out by D.Silver \cite{Dan}, 
As is known to some specialists,
even if the pattern and the companion of 
$\widehat{K}$ are both fibred, 
$\widehat{K}$ may not be fibred \cite{EN}.  
Therefore, the determination of the genus and fibredness 
of a satellite knot $\widehat{K}$ with a fibred companion and 
the pattern $K(2\alpha, \beta|0)$ is not straightforward. 
We leave this problem untouched.  
To the reader who are interested in this problem, 
we refer to \cite{Dan}.% as the best source. 
\end{rem}

\section{Genus one knots} %18

In this final section, 
we determine torti-rational knots
$K(2\alpha, \beta|r)$ of genus one (Theorem \ref{thm:gtOne}).
Recall that torti-rational knots have tunnel number one.
(We denote the tunnel number of $K$ by $t(K)$.)

Recently Scharlemann \cite{S}
proved a conjecture by Goda and Teragaito
which states that if
a knot $K$ has $g(K)=t(K)=1$,
then $K$ is a $2$-bridge knot or satellite knot.
Before that, Goda and Teragaito \cite{GT} had classified
satellite knots $K$ with $g(K)=t(K)=1$.

\begin{prop}\cite[Proposition 18.1]{GT} \label{prop:GT}
Let $K$ be a satellite knot with $g(K)=t(K)=1$.
Then the pattern knot is genus one $2$-bridge knot.
More precisely, the pattern knot is the torti-rational
knot $K(8d, 4d+1 |pq)$ and the companion knot is a
torus knot $T(p,q)$. 
\end{prop}

Note that $(4d+1)/8d = [2, 2d, -2]$, and hence
the pattern knot is a $2$-bridge knot whose
associated continued fraction is $[2d,2pq]$.
In Proposition \ref{prop:GT}, it is evident that
the associated $2$-bridge link has linking number zero,
and that the pattern knot has genus one.
Hence Proposition \ref{prop:GT} immediately follows from
Theorem \ref{thm:gtOne} below.

In Theorem \ref{thm:gtOne}, we have a family of
torti-rational knots $K$ with \mbox{$g(K)=t(K)=1$}, 
which, at a glance, looked like
counter-examples to the Goda-Teragaito conjecture
(See Example \ref{ex:19.1}).

\begin{thm}\label{thm:gtOne}
A torti-rational knot
 $K=K(2\alpha, \beta|r)$ with
 $\ell=\ellk B(2\alpha, \beta)\ge0, r>0$ has
$g(K)=1$ if and only if one of the following is
satisfied:

Case A: $\ell=0$.\\
A1: $\beta/2\alpha=\pm [2, 2d, -2], d\neq 0$ and $r$ is arbitrary.

Case B: $\ell>0$.\\
B1: $\beta/2\alpha=[2,2,2], r=2$\\
B2: $\beta/2\alpha=[2,2,2,2,2], r=1$\\
B3: $\beta/2\alpha=[2,2d,2], d\neq 1, r=1$. 
Note possibly, $d=0$.\\
B4: $\beta/2\alpha=\pm [
  \underbrace{2,2,\ldots,2}_{2a+1},
  2b,
  \underbrace{-2,-2,\ldots,-2}_{2a-1}]$,
  or
 $\pm [
  \underbrace{-2,-2,\ldots,-2}_{2a-1},
  2b,
  \underbrace{2,2,\ldots,2}_{2a+1}]$, $a \ge1, b\neq0, r=1$.
\end{thm}

\begin{rem}\label{rem:tunnel1}
In Theorem \ref{thm:gtOne},
$K=K(2\alpha, \beta|r)$ in A1, or in B4 with $a=1$,
is a $2$-bridge knot.
If $K$ is in B1 or B2, then $K$ is a trefoil knot.
If $K$ is in B3,
then $K$ is a twist knot.
If $K$ is in B4 with $a\ge 2$, then $K$ is a satellite knot
with its companion a torus knot 
$T(a,a+1)$ and the
pattern knot $B(4ab(a+1)-1,2a(a+1))$, 
where the Alexander polynomial is  
$\Delta_{K(1)}(t) = ab(a+1)(t-1)^2 + t$.
\end{rem}

As a direct consequence, we have:

\begin{cor}\label{cor:tunnel2}
A torti-rational knot $K$ of genus one is
a satellite knot if and only if $K$ is in 
Case B4 with $a\ge 2$ of Theorem \ref{thm:gtOne}.
\end{cor}

\begin{rem}\label{cex}
(1) Our knots in B4 with $a\ge 2$
has genus one and hence is prime and not a cable knot.
Therefore, they are negative examples to the question
posed in \cite{AM}.

(2) We can also prove that if $\beta/2\alpha=[4,\pm2,\pm 4]$,
or $[4,\pm3, \pm4]$, then $K(2\alpha, \beta|\pm1)$
cannot be a prime satellite knot 
(c.f. \cite[Theorem 1.6 (2)]{GHS}).
\end{rem}

\begin{ex} \label{ex:19.1} %Example 19.1.
Consider a 2-bridge link $B(46,39)$, where 
$39/46=
$\\
$
[2,2,2,2,2,2,-2,-2,-2]$. 
(This is the case $(a, b)=(2,1)$ 
in Case B4). 
The diagram of $K(r), r=1$, is depicted in \Figure 14.1(a)
together with a compressible Seifert surface $F$.

 Compressing $F$ three times, we obtain
 a minimal genus Seifert surface $F'$ of
 genus $1$, isotopic
 to those depicted in \Figures 14.1 (b), (c) and (d).
Then we see that  $F'$ is a Seifert surface 
for  the  satellite 
knot with its companion 
$T(2,3)$ and the pattern knot a 
$2$-bridge knot $B(23,12)$.%, a twist knot.   
\end{ex}

%\Fi{18.1(1)-(5)}
\eighteen

%\begin{prop}\label{prop:18.1}
%Suppose a knot $K=K(2\alpha, \beta|r), 
%r > 0$ and $\alpha\ge2$ has the genus one.
%Then $K(2\alpha, \beta|r)$ is either
%(a) a $2$-bridge knot or (b) a satellite knot
%with its companion a torus knot 
%$T(a,a+1)$ and the
%pattern knot $B(4ab(a+1)-1,2a(a+1))$, 
%where the Alexander polynomial is  
%     $\Delta_{K(1)}(t) = ab(a+1)(t-1)^2 + t$.
%\end{prop}

%{\it Proof.}
%  (I)  Assume that $K_0$ is not trivial. 
%  Then
% $g(K)\ge|\ell| g(K_0) + 
%  g\bigl(K(2\alpha,\beta|r)\bigr)$,
%where  $\ell =\ellkB$.
%Since $K_0$ is not trivial, $K(2\alpha, \beta|r; K_0)$ is a 
%satellite knot of tunnel number one, and hence, 
%by \cite{MS} the companion is 
%a torus knot $T(a,b)$ 
%for some $|a|, |b| \ge 2$, and the pattern knot is 
%$K(2\alpha, \beta|ab)$.
%If  $\ell \neq0$, then the fact that $g(K)$, 
%$r=ab$, equals one implies 
%that $|\ell|=1$ and $g(K_0) = 1$, and 
%$g\bigl(K(2\alpha, \beta|r)\bigr) =0$.  
%This is impossible
%by \Proposition \ref{prop:8.6++}.
%Therefore, $\ell = 0$.
%Then $g(K) = 
%g\bigl(K(2\alpha, \beta|r)\bigr)$.
%Now, by (11.1), 
%$g(K) = 1$ if and only if 
%$\beta/2\alpha=[2,2d_1,-2], d_1\neq 0$.  
%Therefore, the pattern knot of 
%$K$ is 
%a $2$-bridge knot $B(2pq-1,2p)$ for some $p$ and $q$.
%This proves (I).
% 
%
%(II)  
{\it Proof of Theorem \ref{thm:gtOne}.}
Case A: By (10.1), 
$g(K) = 1$ if and only if 
$\beta/2\alpha=\pm[2,2d,-2], d\neq 0$.  
Therefore, $K$ is 
a $2$-bridge knot associated to
the continued fraction $[\pm2d,2r]$.

Case B:
By \Proposition \ref{prop:8.1},
we see
%(14.1)
\begin{equation}
2(\lambda-1-\rho)+(\lambda-1)(r-1)\ell+
(\lambda-2)(\ell -1)=2.
\end{equation}
Since
each term of the LHS in (14.1) is non-negative,
we have two cases: (i) $\lambda-1-\rho=0$
and (ii) $\lambda-1-\rho=1$.
For case (i), as is seen in the proof of 
\Proposition \ref{prop:8.6++}, we have
$(2\alpha, \beta)=(2\lambda, 2\lambda-1)$
and $\lambda=\ell$,
and hence, either
(a) $\lambda=\ell=2$ and $r=2$, or
(b) $\lambda=\ell=3$ and $r=1$.
Therefore, we have either B1 or B2.
For case (ii), we have either
(c) $\lambda=2$ and hence $r=1$ and $\rho=0$,  or
(d) $\lambda>2, \ell=r=1$ and hence $\rho=\lambda-2$.
Therefore, we have B3 or B4.
This proves Theorem \ref{thm:gtOne}.
\qed

\medskip

\noindent
{\bf Acknowledgements.}
The first author is
partially supported by MEXT, Grant-in-Aid for
Young Scientists (B) 18740035,
and the second author is
partially supported by NSERC Grant~A~4034.
The authors express their appreciation to 
K.Morimoto, M.Sakuma, T.Kobayashi, T.Kanenobu, M.Kidwell 
and D. Silver
for their invaluable comments.
%The first author was 
%partially supported by MEXT,
%Grant-in-Aid for Young 
%Scientists (B) 14740048.
%The second author was
%partially supported by
%NSERC Grant A 4034.

\end{document}